\documentclass[11pt]{article}
\usepackage[margin=1in]{geometry}

\usepackage{amsmath}
\usepackage{color}
\usepackage{mathtools}
\mathtoolsset{showonlyrefs}
\usepackage{amsfonts,amstext, amssymb,amsthm, amsmath, rotating, latexsym,wasysym}
\usepackage{graphicx,bm,bbm,bigints}% ,
\usepackage{adjustbox}
\usepackage{footnote}
\makesavenoteenv{tabular}
\makesavenoteenv{table}

\usepackage{multirow}
\usepackage{booktabs}
\usepackage{bigstrut}
\usepackage{enumitem}
\DeclareRobustCommand{\rchi}{{\mathpalette\irchi\relax}}
\newcommand{\irchi}[2]{\raisebox{\depth}{$#1\chi$}}
\RequirePackage[colorlinks,citecolor=blue,urlcolor=blue]{hyperref} 
\newtheorem{thm}{Theorem}[section]
\newtheorem{lemma}{Lemma}[section]
\newtheorem{cor}{Corollary}[section]
\newtheorem{prop}{Proposition}[section]
\newtheorem{remark}{Remark}[section]

\makeatletter
\def\namedlabel#1#2{\begingroup
    #2%
    \def\@currentlabel{#2}%
    \phantomsection\label{#1}\endgroup
}
\makeatother

\DeclareMathOperator*{\argmin}{arg\,min}

\newcommand{\R}{\mathbb{R}}
\newcommand{\F}{\mathcal{F}}

\newcommand{\p}{\mathbb{P}}

\newcommand{\E}{\mathbb{E}}
\newcommand{\M}{\mathcal{M}}

\newcommand{\B}{\mathfrak{B}}

\definecolor{lgray}{gray}{0.70}
\usepackage{tikz}

\newcommand{\clr}{\color{red}}

\providecommand{\norm}[1]{\left\lVert#1\right\rVert}

\usepackage[numbers]{natbib}

\usepackage{appendix,authblk,libertine}
\renewcommand{\clr}{\color{black}}

\begin{document}

 \title{\bf On Least Squares Estimation Under Heteroscedastic and Heavy-Tailed Errors}
  % \runtitle{LSE under Heavy-tails}

\author[1]{Arun K. Kuchibhotla\thanks{Email: {\tt             arunku@cmu.edu}.}}
 \author[2]{ Rohit K. Patra\thanks{Email: {\tt rohitpatra@ufl.edu}.}} 
\affil[1]{Carnegie Mellon University }
 \affil[2]{University of Florida}
    \date{}
    \maketitle
\begin{abstract}
 We consider least squares estimation in a general nonparametric regression model where the error is allowed to depend on the covariates. The rate of convergence of the least squares estimator (LSE) for the unknown regression function is well studied when the errors are sub-Gaussian. We find upper bounds on the rates of convergence of the LSE  when the error has a  uniformly bounded conditional variance and has only finitely many moments. Our upper bound on the rate of convergence of the LSE depends on the moment assumptions on the error, the metric entropy of the class of functions involved, and the ``local'' structure of the function class around the truth.  We find sufficient conditions on the error distribution under which the rate of the LSE matches the rate of the LSE under sub-Gaussian error. Our results are finite sample and allow for heteroscedastic and heavy-tailed errors.
 \end{abstract}
\noindent%

\section{Introduction} % (fold)
\label{sec:introduction}
% \todo[inline]{
%   1. don't say minimax\\
%   2. don't say optimal\\
% }

% section introduction (end)
% Suppose we have regression data with i.i.d. observations
% $(X_i, Y_i), 1\le i\le n$ with  covariates $X_i$ in a metric space
% $\rchi$ and with a real-valued response $Y_i$. The object of
% interest for us is $f_0(X) = E[Y|X]$, the conditional mean of Y 
% given X. Formally, define the  error as $\epsilon = Y - f_0(X). $
% This error has no special  property except for  $E[\epsilon|X] = 0.$ In particular, epsilon can be arbitrarily dependent on X. 
{\clr Consider the following least squares regression problem: we observe $n$ i.i.d. pairs $(X_i,Y_i)\in   \rchi\times \R, 1\le i\le n$,} where $X_i$ belongs to a metric space $\rchi$ and $Y_i\in \R$.
 The object of interest is  $f_0(x) := \E(Y|X=x)$, the conditional mean of $Y$ given $X=x.$ It is common to define the error to be 
 \begin{equation}\label{eq:error}
 \epsilon := Y- \E(Y|X) = Y- f_0(X).
 \end{equation}
 Note that $\{(X_i, \epsilon_i), 1\le i\le n\}$ are i.i.d.~and $\E(\epsilon_i|X_i)\equiv 0$, by definition. We stress that we do not assume independence between $\epsilon$ and $X$. 
 % under the constraint that $f_0 \in \F$, where $\F$ denotes a class of real-valued functions on $\rchi$.  Formally, the LSE is defined as
% Suppose  we have $n$ i.i.d.~observations $\{(X_i,Y_i)\in   \rchi\times \R, 1\le i\le n\}$
% from the nonparametric regression model
% \begin{equation}\label{eq:LSE_model}
% Y = f_0(X) +\epsilon,
% \end{equation}
% where $\rchi$ is a metric space,  $f_0:\rchi \to \R$ is an unknown measurable function, and $\epsilon$ satisfies $\E(\epsilon|X) = 0$ almost everywhere $P_X$, the distribution of $X$. We allow for $\epsilon$ to be arbitrarily \textit{dependent} on $X$. 
We consider the least squares estimator (LSE) for  $f_0$ under the constraint that $f_0 \in \F$, where $\F$ is a class of real-valued functions on $\rchi$.  The LSE is defined as
\begin{equation}\label{eq:L_2Loss}
\widehat{f}:= \argmin_{f\in \F}\sum_{i=1}^{n} (Y_i - f(X_i))^2.
\end{equation}

 The two most widely used metrics for assessing the error in estimation are the empirical loss ($\|\widehat{f}- f_0\|_n$) and the population loss ($\|\widehat{f} -f_0\|$), where for any function $g: \rchi \to \R$,
\begin{equation}\label{eq:Losses_def}
\|g\|^2_n := \frac{1}{n}\sum_{i=1}^{n} {g}^2(X_i)\; \text{ and }\; \|g \|^2 := {\int_{\rchi} g^2(x) dP_X(x)},
\end{equation}
and $P_X$ denotes the distribution of $X.$
We say that $\widehat{f}$ converges to $f_0$ at a rate $\delta_n$ if $\|\widehat f- f_0\|=O_p(\delta_n)$; $\delta_n$ is also called the  rate of convergence of  the LSE.
 In this paper, we find upper bounds on $\delta_n$  and  the tail probability of $\delta_n^{-1} \|\widehat f- f_0\|$.
{\clr Our goal in this work is to provide some general sufficient conditions on $\F$ and $\epsilon$ under which the LSE is ``rate optimal.'' For instance, LSE for $\gamma$-H{\"o}lder continuous functions (in $d$ dimensions) converges at the minimax optimal rate of  $n^{-\gamma/(2\gamma + d)}$ when $\epsilon$ is sub-Gaussian. In this work we relax the sub-Gaussian assumption to a finite moment assumption, i.e.,  how many finite moments  of $\epsilon$ are required for the LSE to converge at $n^{-\gamma/(2\gamma + d)}$ rate?  We answer this and a general version of the question by studying the LSE under various entropy conditions and  heavy-tailed heteroscedastic noise. 

Informally stating, our results show the following: if the error $\epsilon$ has $q$ moments, then the LSE can attain the minimax rate of convergence whenever $q \ge q^*$ for an explicit threshold $q^* \ge 2$ depending on the complexity of $\mathcal{F}$ and some properties on the locality of $f_0$.
It is worth mentioning that although the LSE under heavy-tailed noise can attain the minimax rate of convergence, its tail behavior suffers when compared to the sub-Gaussian case. The tail probability (i.e.,  $\p(\delta_n^{-1} \|\widehat f- f_0\|\ge t)$) under  heavy-tailed noise  decays polynomially as opposed to a sub-Gaussian decay, $\exp({-c t^2})$, under sub-Gaussian errors. The results of this paper should be seen from the viewpoint of understanding the theoretical behavior of the widely used LSE under realistic assumptions of heavy-tailed heteroscedastic noise. We are not arguing for the universal use of LSE and acknowledge the existence of estimators that outperform the LSE in various aspects. 
% can   for this work is to understand the behavior of the LSE under
  }

%   $n^{-2/5}$ rate? In the above two examples, these rates are known to be minimax; see Sections~\ref{sub:univConvex} and~\ref{sub:multi_index}. 
% % For a general discussion on the connection of minimax rates and entropy numbers used in this paper, see Section~\ref{sub:minimax_rate_and_entropy_numbers}.
%  Our results significantly expand the scenarios under which LSE can be proven to be  rate optimal 
%  % or ``be a safe choice"
%   for the function class $\F$ at hand. It is known 
 \cite{van1996consistency} have established necessary and sufficient conditions on $\F$ and $\epsilon$ for consistency of $\widehat{f}${\clr, in the random design setting with empirical norm}. 
Theorem 3.2.5 of~\cite{VdVW96} implies that the LSE $\widehat{f}$ defined on $\mathcal{F}$ satisfies $\|\widehat{f} - f_0\| =O_p(\delta_n)$ for any $\delta_n$ such that
\begin{equation}\label{eq:Characterization}
\mathbb{E}\left[\sup_{f\in\mathcal{F}:\,\|f - f_0\| \le \delta_n}\big|\mathbb{G}_n[2\epsilon(f - f_0)(X) - (f - f_0)^2(X)]\big|\right] \le C \sqrt{n}\delta_n^2,
\end{equation}
where $C$ denotes a  constant. By ``constant'' we will always mean a quantity that does not depend on $n$ but might depend on the various parameters introduced in our assumptions; we specify these parameters in each occurrence. In the rest of this paper, we make the convention that the constant $C$ is not necessarily the same on each occurrence. \cite{van2017concentration} show that the rate of convergence of the LSE is characterized  by the empirical process above; hence, sharp bounds on the expectation in~\eqref{eq:Characterization} lead to sharp rates for the LSE. Assuming that the functions in $\mathcal{F}$ are uniformly bounded by $\Phi < \infty$, the expectation in~\eqref{eq:Characterization} can be bounded using symmetrization and contraction (Theorem 3.1.21 and Corollary 3.2.2 of~\cite{Gine16}, respectively) by
\begin{equation}\label{eq:NecessaryExpectation}
\mathbb{E}\left[\sup_{f\in\mathcal{F}:\,\|f - f_0\| \le \delta_n}\big|\mathbb{G}_n\left[(|\epsilon| + \Phi)(f - f_0)(X)\right]\big|\right].
\end{equation}
{Inequalities leading to bounds on~\eqref{eq:NecessaryExpectation} are called maximal inequalities.} The path-breaking works by the authors of~\cite{birge1998minimum,MR762984,VANG,VdVW96} have provided sharp maximal inequalities to bound the expectation in~\eqref{eq:NecessaryExpectation}. However, the assumptions are often strong and might not be necessary:~\cite{birge1993rates,VanDeGeer90,VANG,VdVW96} assume restrictive conditions (such as boundedness or sub-exponential tails) on the distribution of $\epsilon$ and \cite{han2018robustness,han2017sharp} assume that $\epsilon$ is independent of $X$.  The study of the LSE in specific examples~\cite{audibert2011robust,shen1994convergence,zhang2002risk} has shown that such conditions are not necessary in general.

The uniform boundedness assumption of $\F$ plays crucial role in proving that~\eqref{eq:NecessaryExpectation} is an upper bound on the left side of~\eqref{eq:Characterization}. The boundedness assumption is widely used in the nonparametric regression setting~\cite{han2017isotonic,guntuboyina2018nonparametric,2016arXiv160804167G,han2017isotonic,kur2019optimality,rakhlin2017empirical,gaillard2015chaining,MR2829871}. However, as pointed out by~\cite{lecue2013learning,mendelson2015local}, the ``gap'' between~\eqref{eq:NecessaryExpectation} and the left side of~\eqref{eq:Characterization} can be large when the noise variance goes to zero with sample size. When $\F$ is sub-Gaussian, Mendelson and co-authors~\cite{lecue2013learning,mendelson2015local} provide provably (see Theorem 1.12 of~\cite{mendelson2015local}) tight bounds for~\eqref{eq:Characterization} even when the noise variance goes to zero. Although sub-Gaussian classes can be unbounded and accommodate vanishing noise variance, a wide range of uniformly bounded nonparametric function classes are not sub-Gaussian; see Section 5.2 and Proposition 3 of~\cite{han2017convergence} for details. For this reason,  in this paper, we focus on uniformly bounded nonparametric function classes and only consider noise distributions with a non-vanishing variance.

% However, in this paper, we focus only on the large noise scenarios, i.e., we assume $\sigma$ is a fixed non-zero constant.{\clr Such an assumption is quite common in statistical models considered in this paper.}

%  A notable departure from this is~\cite{mendelson2014learning,mendelson2018learning}, where $\F$ is allowed to be unbounded. 

 % ~\cite{mendelson2015local,mendelson2016upper} make minimal assumptions on $\epsilon$ but make structural assumptions on $\F$. 

  % In the nonparametric regression setting, $\F$ is commonly assumed to be uniformly bounded~\cite{?}. A notable departure from this is~\cite{mendelson2014learning,mendelson2018learning}, where $\F$ is allowed to be unbounded. 

\subsection*{Organization} % (fold)
\label{sub:organization}

% subsection organization (end)
The paper is organized as follows. In Section \ref{sec:assumptions_and_contributions}, we describe our framework, motivate our assumptions, and list our contributions. In Sections~\ref{sec:L2all},~\ref{sec:MajorGeneral}, and~\ref{sec:rate_of_convergence_for_vc_type_classes}, we find the rate of convergence of the LSE under the three main complexity measures on $\F$ described in Section~\ref{sec:assumptions_and_contributions}.
% In Section~\ref{sec:rate_of_convergence_for_vc_type_classes}, we find the rate of convergence of the LSE when $\F$ satisfies~\eqref{eq:Unif_entr}.
Each section ends with an example, and in each of these examples we show (for the first time) that sub-Gaussian errors are not needed for the LSE to be minimax rate optimal.   In Section~\ref{sec:oracle_ineq}, we briefly comment on the rate of the LSE under misspecification. In Section~\ref{sec:conclusion}, we summarize the contributions of the paper and briefly discuss some future research directions.  In Appendix~\ref{sec:auxiliary_results}, we state three new interpolation inequalities used in our examples. In Appendix~\ref{sec:a_new_maximal_inequality_for_finite_maximums}, we state a new maximal inequality for maximums over finite sets and discuss an application that is of
 independent interest. In Appendix~\ref{sec:PeelingResult}, we state our peeling result. The proofs of all the results in the paper are given in the supplementary file. All the sections, lemmas, and remarks in the supplementary file have the prefix ``S.''

\section{Assumptions and contributions}\label{sec:assumptions_and_contributions} In this section, we describe and discuss our main assumptions. We focus on uniformly bounded function classes $\F$ and relax the assumptions on $\epsilon$ and $\F$ when providing maximal inequalities to bound~\eqref{eq:NecessaryExpectation}. This, in turn, helps us establish the rate of convergence of the LSE under weaker assumptions.  We argue that there  are three properties concerning $\epsilon$ and $\F$ that play a crucial role when finding the rate of convergence of the LSE: (1) the tail behavior of $\epsilon$; (2) the ``complexity'' of $\F$; and (3) the ``local'' structure of $\F$ in the neighborhood of $f_0$. In the following, we discuss these three aspects in detail and state our main assumptions.

\subsection[11]{Assumptions on $\epsilon$}\label{sec:epsilon_assumptions}
 In this work, we  assume that  there exists a $\sigma>0$ such that
\begin{equation}\label{eq:cvar}
\E(\epsilon^2 |X) \le \sigma^2 \text{  almost everywhere (a.e.) }P_X, \tag{CVar}
\end{equation}
and there exists a finite $q\ge 2$  and $K_q < \infty$ such that 
\begin{equation}\label{eq:MomentCondition} \E(|\epsilon|^q)  \le K_q^q. \tag{$\mathcal{E}_q$}
\end{equation}
Note that~\eqref{eq:cvar} allows for heteroscedastic errors that can depend on the covariates arbitrarily and~\eqref{eq:MomentCondition} allows for heavy-tailed errors. Of course, we are not the first to consider heavy-tailed errors (i.e., $\epsilon$ with only finitely many moments). Both~\cite{han2017sharp} and \cite{mendelson2016upper} allow for heavy-tailed errors, but require $\epsilon$ to be independent of $X$ and $\F$ to be sub-Gaussian, respectively; see Section~\ref{sub:our_contributions} for more details on this.~\cite{shen1994convergence} and \cite{chen1998sieve} also allow for heavy-tailed errors, but their results do not directly relate the rate of convergence of the LSE to the moment assumptions on $\epsilon.$ As discussed earlier, we focus only on settings where $\sigma$ is bounded away from zero.
% We discuss these papers and their connections to our results in more detail in the upcoming sections.  
% {\cln Mendelson and his co-authors~\cite{lecue2013learning,mendelson2014learning,mendelson2015local} have studied the behavior of the LSE in heavy tailed regression setting when $\sigma$ is allowed to converge to zero.  In this paper, we focus only on the large noise scenarios, i.e., we assume $\sigma$ is a fixed non-zero constant.} {\clr Such an assumption is quite common in statistical models considered in this paper.}

\subsection[1]{Complexity of $\F$} \label{sec:Complexity}

{\clr Bounds} on~\eqref{eq:NecessaryExpectation} depend on the ``effective'' number of elements in the supremum. The effective number is given by number of functions that are essentially ``different.'' This number is usually described in terms of metric entropy numbers. In the following sections, we use three of the most widely used entropy numbers. For any $\zeta > 0$, function class $\mathcal{F}$, and metric $d(\cdot, \cdot)$ on $\mathcal{F}\times\mathcal{F}$, let $N(\zeta, \mathcal{F}, d)$ be the minimum $m\ge1$ for which there exist functions $\{g_i\}_{i=1}^m$  such that for every $f\in\F$ there exists a $j\le m$ such that $d(f, g_j)\le \zeta$.  We only use metrics $d(\cdot, \cdot)$ of the form $d(f, g) = D(f - g)$ for some norm $D(\cdot)$ and for these forms, we write $N(\zeta, \mathcal{F}, d)\equiv N(\zeta, \mathcal{F}, D)$. $N(\zeta, \mathcal{F}, d)$  and $\log N(\zeta, \mathcal{F}, d)$ are called the $\zeta$-covering number and the $\zeta$-metric entropy of $\F$ with respect to the metric $d$, respectively. For any $f:\rchi \to \R$, define $\|f\|_{\infty}:= \sup_{x\in \rchi} |f(x)|.$ In Section~\ref{sec:MajorGeneral}, we assume 
\begin{equation}\label{eq:inf}
\log N(\zeta,\F, \|\cdot\|_{\infty}) \le A \zeta^{-\alpha}, \quad \text{for some } A>0 \text{ and } \alpha \in [0,2)\tag{$L_\infty$},
\end{equation}
We call $\log N(\zeta,\F, \|\cdot\|_{\infty})$ the $L_\infty$-entropy.  In Section \ref{sec:rate_of_convergence_for_vc_type_classes}, we assume  there exists an $A>0$ such that  $\F$ satisfies
\begin{equation}\label{eq:Unif_entr}
\sup_Q\sup_{\mathcal{F}'\subseteq\mathcal{F} - f_0}\,\log N(\zeta\|F'\|, \mathcal{F}', \|\cdot\|_{2,Q}) \le\frac{ A}{\zeta^{\alpha}} \log^{\beta} \left(\frac{1}{\zeta}\right), \tag{VC($f_0$)} 
\end{equation}
for some $\alpha, \beta \ge 0$,
where  $\F-f_0:= \{f-f_0: f \in \F\}$, $F'(x) := \sup_{g\in\mathcal{F}'}|g(x)|$, the supremum in $Q$ is taken over all finitely supported discrete measures on $\rchi$,  and $\|\cdot\|_{2,Q}$ denotes the $L_2$-norm with respect to the measure $Q.$  If $\F$ satisfies~\eqref{eq:Unif_entr}, then $\F-f_0$ is said to be a uniform VC-type class.

The third entropy considered in the paper is the bracketing entropy. In contrast to covering numbers, the bracketing number $N_{[\,]}(\zeta, \mathcal{F}, d)$ is the smallest $m\ge1$ such that there exist pairs of functions $(g_1^L, g_1^U), \ldots, (g_m^L, g_m^U)$ that satisfy $d(g_j^U, g_j^L) \le \zeta$ for all $j\le m$ and for any $f\in \F$ there exists a $j \le m$ such that $g_j^{L}(x) \le f(x) \le g_j^U(x)$ for every $x\in\rchi.$
In Section~\ref{sec:L2all}, we study the LSE when  $\F$ satisfies 
\begin{equation}\label{eq:Brack}
\log N_{[\,]}(\zeta,\F, \|\cdot\|) \le A \zeta^{-\alpha},  \quad \text{for some } A>0 \text{ and } \alpha \in [0,2),\tag{$L_2$}
\end{equation}
  where $\|\cdot\|$ is the $L_2$-norm  with respect to $P_X$.

In~\eqref {eq:inf},~\eqref{eq:Brack}, and~\eqref{eq:Unif_entr}, $\alpha$ is known as the complexity parameter. For ``simple'' classes of functions, $\alpha$ is small, while a larger $\alpha$ corresponds to more ``complex'' $\F$. For example, when $\F$ is the class of real valued $\gamma$-H\"{o}lder functions on $[0,1]^d$, then $\alpha = d/\gamma$ for the $L_\infty$-entropy; see~\cite[Page 350]{Gine16}.  Note that for H\"{o}lder classes,  larger $\gamma$ (more smoothness) are ``simpler'' classes, thus $\alpha$ is inversely proportional to $\gamma.$ See Table~\ref{tab:Example_comp} for more examples. 

{\clr The entropy conditions~\eqref{eq:inf} and~\eqref{eq:Brack} are ``global'', while~\eqref{eq:Unif_entr} is ``local'', in the sense that~\eqref{eq:inf} and~\eqref{eq:Brack} do not depend on $f_0\in\mathcal{F}$, but~\eqref{eq:Unif_entr} depends on $f_0\in\mathcal{F}$. In particular, condition~\eqref{eq:Unif_entr} may hold for some functions $f_0\in\mathcal{F}$ but might fail to hold for other functions in $\mathcal{F}$; note that~\eqref{eq:Unif_entr} is an entropy condition for $\mathcal{F} - f_0$ not $\mathcal{F}$. For example, the class of all monotone functions satisfies~\eqref{eq:Unif_entr} with $f_0\equiv 0$ (identically zero), but does not satisfy~\eqref{eq:Unif_entr} with $f_0$ is strictly monotone.}

The above metric entropy numbers will be used to bound the expected maxima of the empirical process~\eqref{eq:NecessaryExpectation} via entropy integrals~\cite[Section 2.5]{VdVW96}. In contrast  to the Talagrand’s $\gamma$-functionals~\cite{MR1411488,MR3184689},  the entropy integral based upper bounds can be sub-optimal. However, we use metric {\clr and bracketing} entropies to bound~\eqref{eq:NecessaryExpectation}, because they are  well-understood and are easier to compute for a wide variety of nonparametric classes~\cite{guntuboyina2018nonparametric,VdVW96,VANG,mendelson2008weakly}.

\subsection[2]{Local structure of $\F$} % (fold)
\label{sub:local_structure_of_}
 The expression~\eqref{eq:NecessaryExpectation} can be rewritten as 
\begin{equation}\label{eq:mathcalF_delta_def}
\mathbb{E}\left[\sup_{g\in\mathcal{F}_{\delta_n}}\big|\mathbb{G}_n\left[(|\epsilon| + \Phi)g(X)\right]\big|\right],
\end{equation}
where for any $\delta>0$, $\F_\delta := \{f-f_0 : f \in \F\text{ and }\|f - f_0\| \le \delta\}$. Because the supremum is over functions in $\F_\delta$, if the local ball in $\mathcal{F}$ (centered at $f_0$) is nicely behaved, then bounds on expectations that take into account the local structure will lead to sharper rate bounds. We account for the local structure via the following envelope function
\begin{equation}\label{eq:F_delta_def}
  F_{\delta}(x) := \sup_{f\in\mathcal{F}:\|f - f_0\| \le \delta}|(f - f_0)(x)|.
  \end{equation}
%   the proof of Theorem~\ref{cor:L2_classical} is 
% \begin{equation}\label{eq:calF_delta_def}
% \F_\delta := \{f-f_0: f\in  \F\text{ and }\|f-f_0\| \le \delta\}
% \end{equation} 
We call $F_\delta$ the \textit{local envelope} at $f_0$.\footnote{{\clr Another way to account for local structure is via local entropy bounds, i.e., the entropy of $\mathcal{F}_{\delta}$; see Section~\ref{sub:minimax_rate_and_entropy_numbers} and~\cite[Pages 122, 152]{VANG},~\cite[Section 7]{yang1999information}}.} {\clr The motivation for using the local envelope when bounding~\eqref{eq:mathcalF_delta_def} stems from the results of~\cite[Section 3]{gine1983central}. Theorem 3.3 of~\cite{gine1983central} implies that
\[
\mathbb{E}\left[\sup_{g\in\mathcal{F}_{\delta_n}}\big|\mathbb{G}_n\left[(|\epsilon| + \Phi)g(X)\right]\big|\right] ~\gtrsim~ \frac{1}{\sqrt{n}}\mathbb{E}\left[\max_{1\le i\le n}\,(|\epsilon_i| + \Phi)F_{\delta_n}(X_i)\right].
\]
Also see Theorem 1.4.4 and Remark 1.4.6 of~\cite{MR1666908}. This shows that local envelope is an crucial quantity to consider when bounding the expected value in~\eqref{eq:mathcalF_delta_def}. See Section~\ref{sec:discussion_on_the_local_envelope} of the supplementary file for further discussion.}

Note that $F_\delta$ can depend on $f_0$, however, upper bounds for the integral norms of $F_\delta$ might not depend on $f_0$. Because of this, and notational convenience, we have suppressed the dependence of $F_\delta$ on $f_0$ in our notation.   The local envelope gives us an insight into the worst case behavior of functions in a $\delta$-neighborhood of $f_0$ in $\F$.  If $\sup_{f\in \F} \|f\|_{\infty} \le \Phi$, it is clear that $\|F_\delta\|_{\infty} \le 2 \Phi$, however, if the functions in $\F$ are smooth (e.g., uniformly Lipschitz) then $2 \Phi$ is a conservative bound. In fact, if $\rchi$ is a bounded interval in $\R$ and the functions in $\F$ are uniformly Lipschitz with Lipschitz constant $L$ then $\|F_\delta\|_{\infty} \le 2 L^{1/3} \delta^{2/3}$; see Lemma~2 of~\cite{chen1998sieve} for a proof of this.  In Figure~\ref{fig:LipschitzEnvelope} below, we plot $\F_\delta$ (left panel) and $F_\delta$ (right panel) when $\F$ is the class of $1$-Lipschitz functions and $f_0(x)=x$.

\begin{figure}[h!]
\centering
\includegraphics[width=.9\textwidth]{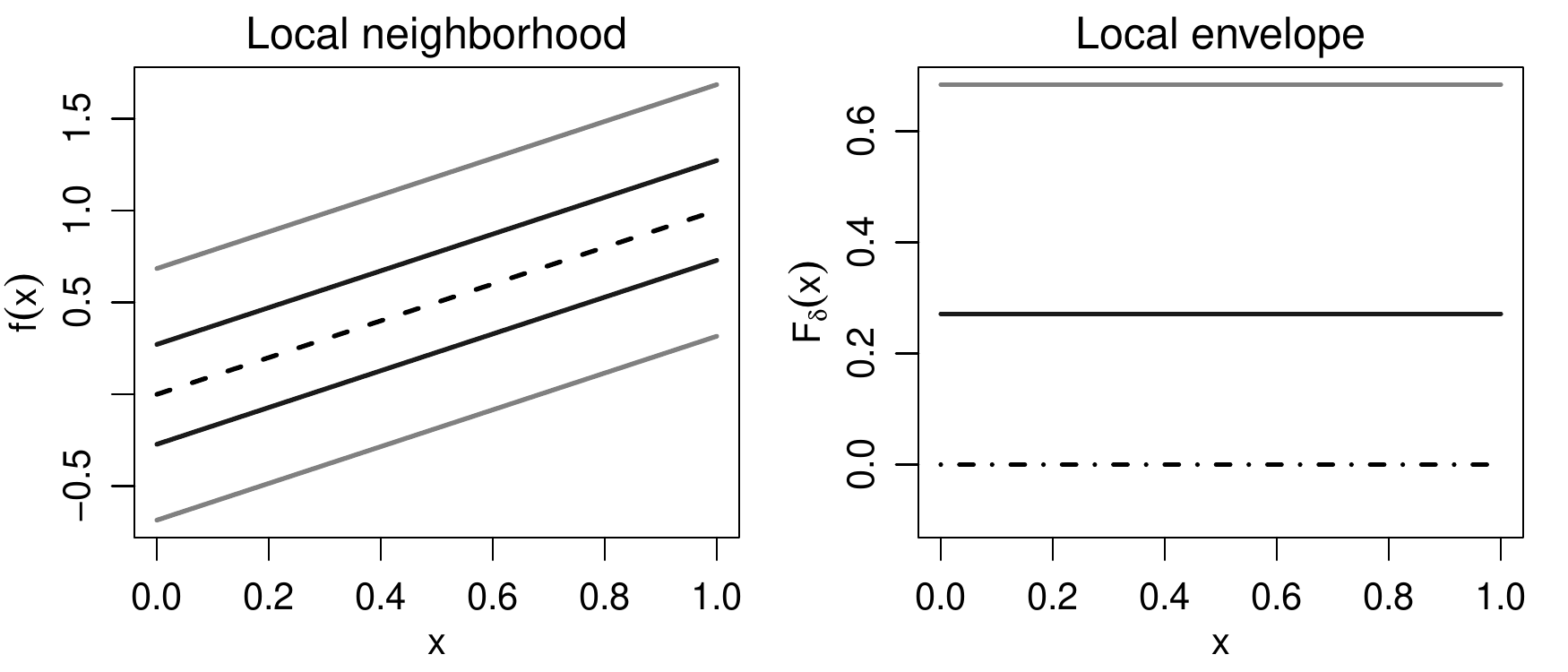}
%% scale=.8
  \caption[]{Illustration of $\F_\delta$ (left panel) and $F_\delta$ (right panel) when $\F:= \{f :[0,1] \to \R | |f(x)-f(y)| \le |x-y|\}$ and $f_0(x)= x$ for  $\delta=.2$ (solid gray) and $\delta=.05$ (solid black).  Any $1$-Lipschitz function $f$ that satisfies $\|f-f_0\| \le .2$ lies in the ``band'' created by the two solid  gray  lines. Here $s= 2/3.$ The dashed line in the left panel is $f_0$.  }
  \label{fig:LipschitzEnvelope}
\end{figure}

 In general, if $\sup_{f\in \F} \|f\|_{\infty} \le \Phi$, then one can invoke the rich theory of interpolation inequalities \cite{Agmon10,kolmogorov1949inequalities,nirenberg2011elliptic} to show that 
\begin{equation}\label{eq:Fdelta}
\|F_\delta\|_{\infty} \le C \Phi^{1-s}\delta^s \quad \text{for some } 0 \le s \le 1,
\end{equation} 
where $C$ denotes a constant independent of $\delta$. {\clr We will refer to conditions of the form~\eqref{eq:Fdelta} as ``envelope growth conditions'' and to the value $s$ as the ``envelope growth parameter.'' Because the local envelope $F_{\delta}$ depends on $f_0$, $s$ can also depend on $f_0$. However, for notational convenience, we will suppress the dependence of $f_0$ on $s$.} Note that for  a uniformly bounded class of functions, $s$ can only vary between $0$ and $1$. If $F_\delta$ does not shrink with $\delta$ (with respect to $\|\cdot\|_\infty$-norm) then $s\approx 0$. For a class of non-smooth functions $s$ will be small, and for the class of infinitely differentiable functions $s=1$; see Table~\ref{tab:interpolation} for examples. Intuitively, smaller values of $s$ correspond to more ``complex'' models.

% If $X\in [0,1]^d$ and $\F = \{\beta^\top x: \beta \in \R^d \}$ then we can easily check that $\F_\delta$ satisfies~\eqref{eq:Fdelta} with $s=1$ for every $f_0$. {\clr Example of $s=0$?}. 
% Many general interpolation inequalities can be used to find $s$ for general smoothness classes; 
Note that the entropy conditions (\eqref{eq:Brack},\eqref{eq:inf}, or~\eqref{eq:Unif_entr}) and~\eqref{eq:Fdelta} complement each other in the sense that the entropy conditions give control over the ``global'' behavior of $\F$ and~\eqref{eq:Fdelta} provides control over the ``local'' behavior of $\F$, i.e., the behavior of $\F_\delta$.  It should be noted that~\cite{chen1998sieve,shen1994convergence,guntuboyina2015global,MR2243881} have implicitly used the property~\eqref{eq:Fdelta} when studying the LSE for certain specific examples. However, their results do not lead to a general relationship between the envelope growth parameter and the rate of convergence of the LSE.  \cite{han2018robustness} use {\clr a similar envelope growth condition} for finding the rate of convergence when $\epsilon$ is independent of $X$ and $\F$ satisfies~\eqref{eq:Unif_entr}. 

 \begin{table}
 % \vspace{-.2in}
 \center
 
      \caption[Different choices of $\F$ and their corresponding parameters] {\label{tab:interpolation}The value of $s$ for widely used choices of $\F$. $P_X$ is the Uniform distribution over $\rchi$. For all of the classes $s$ can be computed from Lemma 2 of~\cite{chen1998sieve}; also see Lemma 4 of~\cite{2020arXiv200810979G}.}
      \centering
  \begin{tabular}{lcl}
  \toprule
  \multicolumn{1}{c}{$\rchi$} &$\F$ & $s$ \\
  \midrule
  $[0,1]^d$ &$\gamma$-H\"older class & $2\gamma/(2\gamma +d)$\\ 
  $[0,1]^d$ & $\gamma$-Sobolev class
  % \footnote{See this \cite[Page 19]{VANG} for a definition.}
   & $(2\gamma-1)/(2\gamma+d-1)$\\
  $[0,1]$ & Uniformly Lipschitz functions & $2/3$\\
 % $[0,1]$ & {\clr Uniformly bounded Convex functions}& $1/2$\\ 
  \bottomrule
  \end{tabular}
  \end{table}

\subsection[minimax rate]{Minimax rates under entropy conditions} % (fold)
\label{sub:minimax_rate_and_entropy_numbers}
In certain cases, it has been shown that the LSE can be rate sub-optimal when the errors have few ($\sim 2$) moments~\cite{han2017sharp,han2018robustness,brownlees2015empirical}. In this paper, our motivation is to understand the required number of moments on the errors $\epsilon$ in order for the LSE to attain the minimax rate of convergence.
% {\clg In other words,  when can the LSE be proven to be  rate optimal or be a ``safe choice''? }
{\clr By the minimax rate, we refer to the \emph{global} minimax rate, i.e.,
\begin{equation}\label{eq:global-minimax-rate}
\inf_{\widehat{g}}\sup_{f_0\in\mathcal{F}}\,\mathbb{E}\left[\|\widehat{g} - f_0\|^2\right],
\end{equation}
where the infimum is over all estimators. This means that we consider the worst case rate over all $f_0\in\mathcal{F}$.}
% Our assumptions~\eqref{eq:Brack},~\eqref{eq:Unif_entr}, and~\eqref{eq:inf} are also global entropy conditions.
%Local rates are discussed in Section~\ref{sec:rate_of_convergence_for_vc_type_classes} under the local entropy condition~\eqref{eq:Unif_entr}.} 
In the regression setup, the relationship between the minimax rate of convergence for an estimator of $f_0$ and {\clr global} entropy numbers of $\F$ is studied extensively by several authors; see e.g.,~\cite{birge1993rates,yang1999information,lecue2013learning}. Theorem B of~\cite{lecue2013learning} (based on the results of~\cite{yang1999information})
 % \todo{Even under extra conditions? Theorem B is the worst case rate for a function $f$.}
 shows that the minimax (rate) lower bound for estimating $f_0$ is $n^{-1/(2 + \alpha)}$, if for all small enough $\eta$, the function class satisfies
\begin{equation}\label{eq:entrop_lowe}
\log N(\eta/2, \mathcal{F}_{\eta}, L_2(P)) ~\gtrsim~ \eta^{-\alpha}.
\end{equation} Furthermore, in addition to~\eqref{eq:entrop_lowe}, if either $\log N_{[\,]}(\eta, \mathcal{F}, L_2(P)) \lesssim \eta^{-\alpha}$ or $\log N(\eta, \mathcal{F}, \|\cdot\|_{\infty}) \lesssim \eta^{-\alpha}$ for all $\eta$ small enough, then~\cite[Theorem~1]{birge1993rates} implies that the minimax rate (both upper and lower bounds) of convergence becomes $n^{-1/(2 + \alpha)}$. It also follows that the best rate of convergence of any estimator of $f_0$ under~\eqref{eq:Brack} or~\eqref{eq:inf} is $n^{-1/(2+\alpha)}$, in the sense that there exists a function space $\mathcal{F}$ (e.g., H\"older classes) satisfying~\eqref{eq:Brack} and~\eqref{eq:inf} for which the minimax rate is no better than $n^{-1/(2+\alpha)}$. Hence under~\eqref{eq:Brack} and~\eqref{eq:inf}, our question becomes: \emph{How many moments on the errors are required before the rate of convergence of the LSE becomes $n^{-1/(2+\alpha)}$?} We will answer this question partly in Sections~\ref{sec:L2all} and~\ref{sec:MajorGeneral} by providing an upper bound on the {\clr minimum} number of moments of $\epsilon$ {\clr required} so that the LSE  has an $n^{-1/(2 + \alpha)}$ rate of convergence. We do a similar study of the LSE under~\eqref{eq:Unif_entr} in Section~\ref{sec:rate_of_convergence_for_vc_type_classes}, see Theorems~\ref{thm:Vc-type} and~\ref{thm:Lowerbnd} and related discussion for more details. 

\paragraph{\clr Local vs global} {\clr In the literature, several authors have considered local versions of entropy conditions and minimax rates of convergence; see~\cite[Section 7]{yang1999information},~\cite[Section 4]{birge1989grenader},~\cite[Section 5]{chatterjee2015risk} and \cite[Section 5]{MR3405621}. The local minimax rate of convergence as opposed to the global minimax rate~\eqref{eq:global-minimax-rate} is given by
\begin{equation}\label{eq:local-minimax-rate}
\inf_{\widehat{g}}\sup_{f\in\mathcal{F}_{\eta_n}}\mathbb{E}[\|\widehat{g} - f\|^2],
\end{equation}
for a sequence $\eta_n$ converging to zero at some rate and $\mathcal{F}_{\delta} = \{f\in\mathcal{F}:\,\|f - f_0\| \le \delta\}$. This is similar to the local asymptotic minimaxity (LAM) considered for parametric inference~\cite[Section 8.7]{vanderVaart98}. In this paper, we mostly restrict ourselves to conditions under which the LSE attains the global minimax rate of convergence. Only in Section~\ref{sec:rate_of_convergence_for_vc_type_classes}, we consider conditions under which the LSE attains the local minimax rate of convergence for shape constrained classes. Also, considering the fact that the rate of convergence of the LSE depends on ``local'' supremums in~\eqref{eq:Characterization} and~\eqref{eq:NecessaryExpectation}, it suffices to consider the local complexity of the function space $\mathcal{F}$. More precisely, it suffices to consider assumptions on local entropy $\log N_{[\,]}(\zeta, \mathcal{F}_{\eta}, \|\cdot\|_2)$, $\log N(\zeta, \mathcal{F}_{\eta}, \|\cdot\|_{\infty})$. Because the local entropy can grow at a slower rate than the global entropy (i.e., lower $\alpha$), one can obtain faster rates of convergence for the LSE using the local entropy conditions; see~\cite{yang1999information} and~\cite[Chapter 7.5]{VANG}. 
{This is not very common and in most non-parametric examples, both the local and global entropies grow at the same rate.}
For this reason, we will restrict ourselves to the global entropy conditions~\eqref{eq:Brack},~\eqref{eq:inf}, and~\eqref{eq:Unif_entr}. Note that we do consider local structure of the function space $\mathcal{F}$ using the local envelope as discussed in Section~\ref{sub:local_structure_of_}. 
% It would be interesting to consider the effect of both the local envelope and the local entropy on the rate of convergence of the LSE. However, in most non-parametric examples, the local entropy at the same rate as the global entropy and hence, 
}
% Theorem 7 of~\cite{rakhlin2017empirical} proves that under the uniform entropy assumption~\eqref{eq:Unif_entr} (alone), the rate of convergence for estimating $f_0$ cannot be better than $n^{-1/(2+\alpha)}$. This means that there exists a function space satisfying~\eqref{eq:Unif_entr} and an $f_0\in\mathcal{F}$ such that the minimax rate for estimating $f_0$ is $n^{-1/(2+\alpha)}$. However, our results in {}Section~\ref{sec:rate_of_convergence_for_vc_type_classes} assume a growth condition~\eqref{eq:GeneralEnvelopeCondition_VC} on the local envelope, in addition to~\eqref{eq:Unif_entr}. Under these two conditions, the minimax rate is no longer $n^{-1/(2+\alpha)}$. The details for this case will be given in the next section.

  % subsection minimax_rate_and_entropy_numbers (end)
% subsection local_structure_of_ (end)
\subsection{Our contributions} % (fold)
\label{sub:our_contributions}
When $\F$ satisfies \eqref{eq:Brack} or~\eqref{eq:inf} with complexity parameter $\alpha\in [0,2]$ and $\epsilon$ is uniformly sub-Gaussian (i.e,  $ \E(|\epsilon|^q|X)^{1/q} \le C\sqrt{q}$ almost everywhere $P_X$ for all $q\ge2$), then the LSE is known to be minimax rate optimal and it converges at an $n^{-1/(2+\alpha)}$ rate~\cite[Chapter 3.4]{VdVW96}. In this paper, we show that for a wide variety of examples,  the uniform sub-Gaussianity of $\epsilon$ is not necessary for the $n^{-1/(2+\alpha)}$ rate of convergence of the LSE. We further provide tail probability  bounds for $\delta_n^{-1} \|\widehat{f}-f_0\|$ that decay as a polynomial of degree $q$ (as opposed to the Gaussian decay under sub-Gaussian errors).  Our framework is closely related to the works~\cite{han2017sharp,han2018robustness} and~\cite{mendelson2016upper}, but with the following important differences:
\begin{enumerate}
  \item \cite{han2017sharp,han2018robustness} study the rate of convergence of LSE under $q$ moments on errors, but assume independence between errors and covariates. In this paper, we do \emph{not} assume independence and allow for arbitrary dependence between the errors and covariates (except possibly for conditional moment assumptions). In addition to the~\eqref{eq:Brack} and~\eqref{eq:Unif_entr} assumptions in~\cite{han2017sharp}, we also consider function classes satisfying~\eqref{eq:inf}. Furthermore, we also study the impact of {\clr $s$, the envelope growth parameter,} on the rate of convergence of LSE under~\eqref{eq:Brack},~\eqref{eq:inf}, and~\eqref{eq:Unif_entr}, while~\cite{han2018robustness} considers the effect of the local envelope under~\eqref{eq:Unif_entr} only when $\epsilon$ and $X$ are independent.
  \item \cite{mendelson2016upper,mendelson2015local} allows for $q$ moments on errors as well as arbitrary dependence between errors and covariates. However, the authors require the function class $\mathcal{F}$ to be sub-Gaussian i.e., $(\mathbb{E}[|f(X) - g(X)|^p])^{1/p} \lesssim \sqrt{p}\|f - g\|$ for all $p\ge2$, $f,g \in \F$; a closely related relaxation is~\cite[Definition 1.7]{mendelson2016upper}. The sub-Gaussian condition implies the small-ball condition~\cite{mendelson2014learning,van2014higher}, which is not satisfied for several function classes (such as $\gamma$-H{\"o}lder continuous, monotone, and convex functions) we consider; see Section 5.2 and Proposition 3 of~\cite{han2017convergence} for details.
\end{enumerate}
In Sections~\ref{sec:L2all}, \ref{sec:MajorGeneral}, and \ref{sec:rate_of_convergence_for_vc_type_classes}, we relate the rate of convergence of the LSE to the behavior of $F_\delta$ and the moment assumptions on $\epsilon$ when $\F$ satisfies the global conditions~\eqref{eq:Brack},~\eqref{eq:inf}, and~\eqref{eq:Unif_entr}, respectively; our results are   summarized in  Table~\ref{tab:Rates_obtained}. Although we show that the rate of convergence of the LSE under heavy tails and sub-Gaussian error match for some choices of $q$, the tail behaviors of $\delta_n^{-1}\|\widehat f-f_0\|$ differ for every $q<\infty.$ In fact, we show that under~\eqref{eq:MomentCondition}, $\p(\delta_n^{-1}\|\widehat f-f_0\|\ge D)\lesssim D^{-q}$; see~\cite[Theorem 5.1]{MR2829871} for the tail behavior under bounded errors. {\clr This tail behavior is optimal for the LSE under~\eqref{eq:MomentCondition}; see~\cite[Proposition 1.5]{lecue2016performance}. However, there do exist several robust estimators that are  minimax rate optimal and have sub-Gaussian tails even under heavy-tailed errors~\citep{lugosi2016risk,mendelson2019unrestricted}.}

% In this regard, heavy tailed errors affect the behavior of the LSE i

 \begin{table}[!ht]
 \centering
      \caption[Rate of convergence of $\|\widehat{f}-f_0\|$ under entropy and {\clr envelope growth} assumptions.]{\label{tab:Rates_obtained}Rate of convergence of $\|\widehat{f}-f_0\|$ under entropy and {\clr envelope growth} assumptions when $\Phi=\sup_{f\in \F} \|f\|_{\infty} $ and $A$ in~\eqref{eq:inf},~\eqref{eq:Brack}, and~\eqref{eq:Unif_entr} do not change with $n$.}

      \resizebox{\textwidth}{!}{
  \begin{tabular}{lccc}
  % \multicolumn{3}{c}{Assumpti}&\multicolumn{2}{c}{ Moments needed for an $n^{-1/(2+\alpha)}$ rate}\\
  % \cmidrule(r){1-4} \cmidrule(l){5-6}
  \toprule
  Entropy  &  Envelope growth assumption & Moment assumption & Rate of convergence\\
  \midrule
  \eqref{eq:Brack}  &$\big\|(|\epsilon| +\Phi) F_{\delta}(X)\big\|_q \le C \Phi^2 \delta^{s}$ & $q\ge 2/s$ & $n^{-1/(2+\alpha)}$\\
  \eqref{eq:inf} &$\|F_\delta\|_\infty\le  C \Phi^{1-s} \delta^{s}$ & $q\ge\frac{ 2 + \alpha(1-s)}{s+ \alpha(1-s)}$ & $n^{-1/(2+\alpha)}$\\
  \eqref{eq:Unif_entr}& $\big\|  F_{\delta}\big\| \le C  \Phi^{1-s} \delta^{s}
$ &$q\ge 2$& $n^{-1/(2(2-s))}$\\
  \bottomrule
  \end{tabular}}
  \end{table}

% We will now briefly describe our results in Sections~\ref{sec:L2all},~\ref{sec:MajorGeneral}, and~\ref{sec:rate_of_convergence_for_vc_type_classes}.  

In Section~\ref{sec:L2all}, we consider classes of functions that satisfy~\eqref{eq:Brack}. In Theorem~\ref{thm:SecondMomentRate}, we show that if $F_\delta$ satisfies an $L_q$ version of~\eqref{eq:Fdelta} with {\clr envelope growth} parameter $s$, then the LSE converges at an $n^{-1/(2+\alpha)}$ rate if $\epsilon$ has at least $2/s$ (conditional) moments.
% ~\cite{han2017sharp} show that if $\F$ satisfies~\eqref{eq:Brack},  then the LSE converges at an $n^{-1/(2+\alpha)}$ rate if $\epsilon$ is independent of $X$ and has $1+2/\alpha$ moments no matter the local structure of $\F$ (i.e., even if $s=0$). 
In Section~\ref{sub:univConvex}, we apply Theorem~\ref{thm:SecondMomentRate} to show that the convex LSE converges at a near minimax rate if $\E(|\epsilon|^3|X) \le C<\infty$ a.e. $P_X$ and $f_0$ is bounded. Previously, minimaxity of the convex LSE with heteroscedastic errors was known only under uniformly sub-Gaussian errors, i.e,  $ \E(|\epsilon|^q|X)^{1/q} \le C\sqrt{q}$ for all $q\ge2$ a.e.~$P_X$.

% In both cases, our results show that the LSE converges faster if $\F$ satisfies the additional ``smoothness'' condition~\eqref{eq:Fdelta} with $s>0$. In Section~\ref{sec:L2all}, we actually find the rates of convergence of the LSE under a weaker version ($\|\cdot\|$-norm instead $\|\cdot\|_\infty$-norm) of the assumption in~\eqref{eq:Fdelta}; see~\eqref{eq:GeneralEnvelopeCondition_L2}.

In Section~\ref{sec:MajorGeneral}, we show that if $\F$ satisfies~\eqref{eq:inf}, then the LSE converges at an $n^{-1/(2+\alpha)}$ rate if $\epsilon$ has at least $1+ 2/\alpha$ moments. However, only $(2 + \alpha(1-s))/(s+ \alpha(1-s))$ ($< 1 + 2/ \alpha $) many moments for $\epsilon$ are enough if $\F$ satisfies~\eqref{eq:Fdelta} with envelope growth parameter $s>0$; see Theorem~\ref{thm:L_inf_smooth}. This is useful since classes with low  complexity $\alpha$ often have high {\clr envelope growth parameter} $s$.
 % A simple counter example is the class of indicators of closed intervals on $[0,1]$, i.e., $\{\mathbf{1}_{[a,b]}(\cdot): 0\le a \le b \le 1\}$. 
% For classes with complexity $\alpha$ and local smoothness $s>0$, $(2 + \alpha(1-s))/(s+ \alpha(1-s)) < 1+ 2/\alpha.$ 
In Section~\ref{sub:multi_index}, we apply Theorem~\ref{thm:L_inf_smooth} to find moment conditions on $\epsilon$ under which the LSE is minimax rate optimal for $d$-dimensional H\"{o}lder regression and some low dimensional submodels. 

In Section~\ref{sec:rate_of_convergence_for_vc_type_classes}, we consider classes of functions that satisfy~\eqref{eq:Unif_entr} and a $L_2$ version of~\eqref{eq:Fdelta}.  
% However, in many scenarios it turns out that $\F_\delta$ is a Uniform VC-type class. For example, if $\rchi=\R$, $\F$ is the class of monotone (or convex) functions on $\R$ and $f_0\equiv 0$ ($f_0(t) =t$) then $\F_\delta$ is uniform VC-type class.
 % In Section~\ref{sec:rate_of_convergence_for_vc_type_classes}, w
 We show that the LSE converges at a rate of $n^{-1/(2(2-s))}$ for any $\alpha<2$ when $\epsilon$ has just two moments. Theorem~\ref{thm:Vc-type} also allows for non-Donsker $\F$, i.e., $\alpha\ge 2$. Theorem~\ref{thm:Vc-type} is especially useful in proving adaptive properties for shape constrained LSEs; see~Section~\ref{sub:example_3_univariate_isotonic_regression} and Remark~\ref{rem:Extension}. In Theorem~\ref{thm:Lowerbnd}, we  show {\clr (via an example)} that the upper bound in Theorem~\ref{thm:Vc-type} cannot be improved (up to $\log n$ factors). 
 % It proves that under~\eqref{eq:Unif_entr} and~\eqref{eq:cvar}, the local envelope drives the rate of convergence of the LSE. 
 % In the supplementary file, we give a heuristic  argument that suggests that the local envelope has role in the rate of convergence for the LSE in general.

% {\clr For example, when $\F$ is class of monotone functions ($\F$ satisfies~\eqref{eq:Brack} with~$\alpha=1$), then~\cite{zhang2002risk} showed that the LSE converges at the rate $n^{-1/3}$ as long as $\epsilon$ has just  two moments.  But the proof of this uses of specific (closed form) characterization of the LSE and geometric properties of the class of monotone functions.  In contrast, we use the local and global structure of $\F$ to find the rate of convergence of $\widehat{f}$. Thus the approach is applicable to other similar  scenarios.} 

 The first step in proving the results discussed above is to find a tight upper bound on~\eqref{eq:NecessaryExpectation}. 
% {\clr Such bounds are known as maximal inequalities.} 
Then one uses the upper bound in conjunction with a peeling argument~\cite[Theorem 3.2.5]{VdVW96} to bound the tail probability of the LSE. However, most existing maximal inequalities require  $|\epsilon| (f-f_0)(X)$ to be bounded or to have exponential tails.  It should be noted that, while several exceptions such as~\cite[Lemma 6.12]{MR1915446},~\cite{van2011local}, and~\cite[Theorem 1.9]{mendelson2016upper} exist,  they are either not sharp enough or do not apply in settings considered here. In order to deal with heavy-tailed errors, in this paper, we introduce a peeling argument in Theorem~\ref{thm:GeneralPeeling}. This peeling argument uses a truncation device to split the bound on the tail probability into two parts; see~\cite{shen1994convergence,mendelson2008weakly} for other truncation based arguments.  We use new (Proposition~\ref{lem:bernouli_maximal}) and some existing (\cite[Lemma 3.4.2]{VdVW96} and Lemma~\ref{lem:UniformCoveringPhi}) maximal inequalities to bound the maximum of the truncated empirical process and the Markov inequality to control the unbounded remainder. Then we optimize over the truncation scale to find the rate of convergence; see proof of Theorem~\ref{thm:GeneralPeeling}. 

% Our results in Sections~\ref{sec:L2all}--\ref{sec:rate_of_convergence_for_vc_type_classes} should be thought of as the ``worst-case'' rates. The results are trying to find the best rates for the LSE under \textit{only} the entropy and {\clr envelope growth} conditions.

    % \item The moment condition on $\epsilon$ in Theorem~\ref{thm:rate_m_theta_LLSE} is slightly weaker  than Theorem 4 of~\cite{han2017sharp}. We require $\|\epsilon\|_p < \infty$ and while~\cite{han2017sharp} requires that $\|\epsilon\|_{p,1} < \infty$ (note that $\|\epsilon\|_{p} < \infty$ implies that $\|\epsilon\|_{p,1} < \infty$).
   
%\end{remark}

% subsection our_contributions (end)

\section[Rate of convergence under bracketing entropy]{Rates of convergence of the LSE using bracketing $L_2(P_X)$-entropy} % (fold)
\label{sec:L2all}
Assumption~\eqref{eq:Brack} is the most widely used entropy condition to study the rate of convergence of the LSE~\cite{Gine16,han2017sharp,VdVW96}.    The following theorem (proved in Section~\ref{sec:proofSecondMomentRate} of the supplementary file) finds an upper bound on the  rate of convergence of the LSE when $\epsilon$ is heavy-tailed and heteroscedastic. 
% The only assumption we make on the relationship between $\epsilon$ and $X$ is that the conditional variance of $\epsilon$ given $X$ is bounded almost every $X$. 

% Suppose $\F$ satisfies~\eqref{eq:inf}, $\epsilon$ satisfies~\eqref{eq:cvar} and~\eqref{eq:MomentCondition}, and $f_0\in \F$.  Let $\Phi:= \sup_{f\in \F} \|f\|_{\infty}$.  If there exists a constant $C>0$such that 
\begin{thm}\label{thm:SecondMomentRate}
Suppose $\F$ satisfies~\eqref{eq:Brack}, $\epsilon$ satisfies~\eqref{eq:cvar},  and $f_0\in \F$. Let $\Phi:= \sup_{f\in \F} \|f\|_{\infty}$. Suppose there exists a constant $C>0$  such that 
\begin{equation}\label{eq:GeneralEnvelopeCondition_L2}
\big\|(|\epsilon| +\Phi) F_{\delta}(X)\big\|_q \le C \Phi^2 \delta^{s},
\end{equation}
for some $s\in[0,1]$, and let
\begin{equation}\label{eq:rn_SecondMoment}
r_n := \min\left\{\frac{(nA^{-1})^{1/(2+\alpha)}}{(\sigma+\Phi)^{2/(2+\alpha)}}, 
% \frac{n^{1/2}}{\sigma + \Phi},
\frac{n^{(q-1)/(q(2-s))}}{\Phi^{2/(2-s)}},  
\frac{n^{1/(2 + \alpha + (2-qs)/(q-1))}}{(A^{q-1}\Phi^{2q})^{1/(2-qs+(2+\alpha)(q-1))}} \right\}.
\end{equation}
Then,  there exists a constant $C$ depending only on $\alpha$, $s$, and $q$ such that
\begin{equation}\label{eq:Tail_boundLSE}
\mathbb{P}\left(r_n\|\widehat{f} - f_0\| \ge D\right) 
% \le   \frac{C}{D^{2q/(2q-1)}}
\le C D^{-q + \mathbf{1}\{s=1\}/10}
\end{equation}
for any $n \ge 1$ and $D> 0$. 
% where $C$ depends only on $\alpha$, $s$, and $q$, and 
% \begin{equation}\label{eq:rn_SecondMoment}
% r_n := \min\left\{\frac{(nA^{-1})^{1/(2+\alpha)}}{(\sigma+\Phi)^{2/(2+\alpha)}}, 
% % \frac{n^{1/2}}{\sigma + \Phi},
% \frac{n^{(q-1)/(q(2-s))}}{\Phi^{2/(2-s)}},  
% \frac{n^{1/(2 + \alpha + (2-qs)/(q-1))}}{(A^{q-1}\Phi^{2q})^{1/(2-qs+(2+\alpha)(q-1))}} \right\}.
% \end{equation}
Hence, the rate of convergence of the LSE is $r_n^{-1}.$
% \todo[inline]{The rate above needs to be fixed because we changed the bound in~\eqref{eq:GeneralEnvelopeCondition_L2}.}
\end{thm}

\begin{remark}[Assumptions in Theorem~\ref{thm:SecondMomentRate}]\label{rem:assumptions}
We now make some observations on the assumptions of Theorem~\ref{thm:SecondMomentRate}. 
\begin{enumerate}
  \item  The covariate space $\rchi$ is not restricted to be Euclidean.  The only assumption on $\rchi$ is that it be a metric space. This comment applies to all the results of the paper. 

  \item  Observe that~\eqref{eq:MomentCondition} and~\eqref{eq:Fdelta} together imply~\eqref{eq:GeneralEnvelopeCondition_L2}, i.e., if $\|F_{\delta}\|_{\infty} \le C\Phi\delta^s$ and $\E(|\epsilon|^q) \le K_q^q$ then
\[\|(|\epsilon| + \Phi)F_{\delta}(X)\|_{q} \le (K_q + \Phi)C\Phi^{1-s}\delta^s \le C \Phi^2\delta^s, \; \text{where }C  := 1 + K_q/\Phi.\]

\item There are cases when $\|F_\delta\|_{\infty}\asymp 1$ but the above upper bound can still hold when $\|F_{\delta}\|_{q} \le C\Phi\delta^s$ and $\E(|\epsilon|^q|X) \le K_q^q$ a.e.~$P_X$; see  Section~\ref{sub:univConvex}.

 \item  The uniform boundedness assumption on $\mathcal{F}$ can be easily  relaxed to only $\|\widehat{f}\|_\infty =O_p(1)$; see Section~\ref{sub:univConvex} for a detailed argument. Also see~\cite{han2017sharp,han2018robustness,2017arXiv170800145K} for further examples. 

 \item {\clr We refer to condition~\eqref{eq:GeneralEnvelopeCondition_L2} as ``$L_2$-envelope growth condition.'' This growth condition can be relaxed to accommodate extra log factors.} For example, if $\big\|(|\epsilon| +\Phi) F_{\delta}(X)\big\|_q \le C \Phi^2 \delta^{s} \log^{\gamma}(1/\delta)$ then $r_n$ will increase by additional $\log n$ factors; where the power of $\log n$ will depend on~$\gamma, \alpha$ and $q$. This dependence is computed  explicitly in~\eqref{eq:ep_n_extra_log} in Section~\ref{sub:additional_log_factors} of the supplementary file. \label{item:extraLogfactor}
\end{enumerate}
\end{remark}

% If $s<1$ then the~\eqref{eq:Tail_boundLSE} implies that $\E(r_n^\eta \|\widehat{f} - f_0\|^\eta)$ 
\begin{remark}[Conclusions of Theorem~\ref{thm:SecondMomentRate}]\label{rem:Tailbound}
{\clr Note that the rate of convergence $r_n^{-1}$ is a function of $f_0$ because both the local envelope ($F_{\delta}$) and envelope growth parameter ($s$) depend on $f_0$.}
 The tail bound in~\eqref{eq:Tail_boundLSE} is a finite sample result and holds for all $n\ge 1$ {\clr and hence one can take the supremum over all $f_0\in\mathcal{F}$ with a fixed value of $s$ on the left hand side of~\eqref{eq:Tail_boundLSE}}. When $s=1$,~\eqref{eq:Tail_boundLSE} implies that the tail probability decays at a polynomial rate with an exponent of $-q+1/10$. The $1/10$ in the exponent is meant to represent a small constant. In fact when $s=1$, we show that $\mathbb{P}(r_n\|\widehat{f} - f_0\| \ge D) 
\le C D^{-\eta},$ for any $\eta <q$; see~\eqref{eq:TrueTailBound_L2} in Section~\ref{sec:proofSecondMomentRate} for a proof of this. Here the constant $C$ depends on $q, \alpha$, and $\eta$ only.  Because $q\ge 2$, the tail probability bound in~\eqref{eq:Tail_boundLSE} implies that $\E(r_n \|\widehat{f} - f_0\|)\le C$ for all $n\ge 1$ and some constant $C$. {\clr It should be noted that when the errors are sub-Gaussian, the tail probability of the LSE decays like $\exp(-D^2/c)$ for some constant $c$; see e.g., Theorem 5.1 of \cite{MR2829871}. Also see \cite{lugosi2016risk} for estimators that have a sub-Gaussian tail even under heavy-tailed noise. }
\end{remark} 

Because all the assumptions and results are finite sample,  $\F$, $A$, $\Phi$, and the distribution of $(X, Y)$ are all allowed to depend on $n$.  However, in most applications, these do not change with $n$, and hence the dependence of $r_n$ on $A$ and $\Phi$ in~\eqref{eq:rn_SecondMoment} can be ignored. Furthermore, if $\Phi < \infty$ and $\epsilon$ satisfies~\eqref{eq:MomentCondition},  then  \textit{every} uniformly bounded function class $\F$ satisfies~\eqref{eq:GeneralEnvelopeCondition_L2} with $s=0.$ The following corollary finds the rate of the LSE if $\F$ does not satisfy the {\clr envelope growth} assumption of~Section~\ref{sub:local_structure_of_}, i.e., $s=0$.

  % Another common scenario is that $F_\delta$ does not shrink with $\delta$, i.e., $s=0$. If both of the above scenarios hold true, then the result of Theorem~\ref{thm:SecondMomentRate} is greatly simplified. 

\begin{cor}\label{cor:L2_classical}
Suppose $f_0\in \F$ and $\Phi:= \sup_{f\in \F} \|f\|_{\infty}$ is a constant. Moreover, suppose  $\epsilon$ satisfies~\eqref{eq:cvar} and~\eqref{eq:MomentCondition} and $\F$ satisfies~\eqref{eq:Brack}.  Then for any $n \ge 1$ and $D > 0$, we have 
\begin{equation}\label{eq:rn_Classical}
\mathbb{P}\left(n^{{1}/{(\alpha+ 2 q/(q-1)})} \|\widehat{f} - f_0\| \ge D\right)\le  {C}{D^{-q}},
\end{equation}
where $C$ is a constant depending only on $A, \Phi, K_q, q, \sigma,$ and $\alpha$.
\end{cor}
  % \item There are estimators that 
% \end{enumerate}
% \end{remark}

The above result is a direct application of Theorem~\ref{thm:SecondMomentRate} with $s=0$. If $\epsilon$ and $X$ are further assumed to be independent then Theorem~3 of \cite{han2017sharp} shows that the LSE converges at a rate of $n^{-1/(2+\alpha)}$ when $q\ge 1+ 2/\alpha$. 
% Corollary~\ref{cor:L2_classical} allows for heteroscedastic and heavy-tailed errors, but {\clr shows that this relaxation \emph{might} come at a cost.} {\clg this relaxation comes at a cost.} 
The rate of convergence obtained in~\eqref{eq:rn_Classical} is strictly slower than the minimax rate  for this setup.  A similar sub-$n^{-1/(2+\alpha)}$ rate was found in \cite[Section~3.4.3.1]{VdVW96} for  fixed design regression ($X_1,\ldots, X_n$ are fixed and non-random) with heavy-tailed errors. There are two possible explanations for the rate bound in~\eqref{eq:rn_Classical}: (1) the LSE is not minimax rate optimal under the assumptions of Corollary~\ref{cor:L2_classical} and  there exists some dependence structure between $\epsilon$ and $X$ and a choice of $\F$ such that the convergence rate of LSE is $n^{-{1}/{(\alpha+ 2 q/(q-1)})}$; or (2) the LSE actually converges at an $n^{-1/(2+\alpha)}$ rate and the obtained rate is an artifact of the proof. The optimality of Corollary~\ref{cor:L2_classical} is still an open problem.

\begin{remark}\label{rem:Comp_Han_MCLT} 
 \cite[Proposition 3 and Remark 10]{han2017sharp} argue that under~\eqref{eq:Brack},  the rate of convergence of the LSE can be arbitrarily slow when $\epsilon$ is heteroscedastic and has heavy tails. On surface, this might seem to be at odds with Corollary~\ref{cor:L2_classical}, but in their examples,  both $\F$ and $\E(\epsilon^2|X)$ are unbounded.  This is important because~\eqref{eq:cvar}, the boundedness of $\E(\epsilon^2|X)$, is a crucial assumption in all our results.  We use condition~\eqref{eq:cvar} to provide bracketing entropy bounds for $\{\epsilon(f - f_0):f\in\mathcal{F}\}$ based on the bracketing entropy bounds for $\F-f_0$. Note that if $[\ell, u]$ is the bracket for $f-f_0$, i.e., $\ell \le f - f_0 \le u$ then 
 $\epsilon_+\ell - \epsilon_-u \le \epsilon(f - f_0) \le \epsilon_+u - \epsilon_-\ell,$
  where $\epsilon_+$ and $\epsilon_-$ are the positive and negative parts of $\epsilon$, respectively.  The width of this bracket is $|\epsilon|(u - \ell)$. Under assumption~\eqref{eq:cvar}, we have
$\|\epsilon(u - \ell)\| \le \sigma\|u - \ell\|.$
Therefore under~\eqref{eq:cvar}, $N_{[\,]}(\eta, \{\epsilon(f - f_0):f\in\mathcal{F}\}, \|\cdot\|) \le N_{[\,]}(\eta/\sigma, \F-f_0, \|\cdot\|)$. This crucial conclusion might not hold if~\eqref{eq:cvar} is not satisfied. 

\end{remark}

The rates of convergence  in Corollary~\ref{cor:L2_classical} does not take into account any  structure of $\F$ other than the complexity {\clr (entropy)} of the function class.  Theorem~\ref{thm:SecondMomentRate} improves upon Corollary~\ref{cor:L2_classical} by using the {\clr envelope growth condition on} $\F$ (around $f_0$); see~\cite{audibert2011robust,chen1998sieve,han2018robustness,shen1994convergence,VanDeGeer90} for results that use a similar {\clr envelope growth condition} implicitly or explicitly. Theorem~\ref{thm:SecondMomentRate} shows that the LSE will converge at an $n^{-1/(2+\alpha)}$ rate if $\epsilon$ has  enough moments. 
To better understand the rate in Theorem~\ref{thm:SecondMomentRate}, let us assume that both $A$ and $\Phi$ are constants (do not change with $n$). In this case, $r_n$ (in Theorem~\ref{thm:SecondMomentRate}) can be simplified to
\[ r_n \asymp  \min\left\{ {n^{1/(2+\alpha)}},{n^{(q-1)/(q(2-s))}},  {n^{(q-1)/(q(2-s)+\alpha(q-1))}} \right\}.\]
 Furthermore, observe that
% \vspace{-.2in}
\begin{align}\label{eq:rate}
\begin{split}
\frac{1}{2+\alpha}  \le \frac{q-1}{q(2-s)+\alpha(q-1)} \le \frac{q-1}{q(2-s)}  &\quad \Leftrightarrow\quad  q \ge \frac{2}{ s}.
\end{split}
\end{align}
Thus 
\begin{equation}\label{eq:rn_final_L2}
r_n \asymp  \min\left\{ {n^{1/(2+\alpha)}},  {n^{(q-1)/(q(2-s)+\alpha(q-1))}} \right\} \text{ for all }q\ge 2
\end{equation}
and if $q\ge 2/s$ then $r_n= n^{1/(2+\alpha)}$. 
The above calculations suggest an interesting interplay between $\alpha$, $q$, and $s$. They show that if  $\E(|\epsilon|^{2/s}|X) \le C <\infty$ and $\|F_\delta\|_{2/s} \le C \delta ^{s}$, then the rate of convergence of the LSE under the heavy-tailed heteroscedastic  errors  is $n^{-1/(2+\alpha)}$ and this rate coincides with the rate under sub-Gaussian errors. 
  % This justifies the usage of least squares estimators under heavy-tailed errors in a wide variety of examples.
 However, if $\epsilon$ has less than $2/s$ moments then Theorem~\ref{thm:SecondMomentRate} suggests that there might exist ``hard'' settings where the ``noise'' is too strong and the guaranteed rate of convergence for the  LSE is slower than $n^{-1/(2+\alpha)}$. 
% The optimality of~Theorem~\ref{cor:L2_classical} is ~\cite[Proposition 3]{han2017sharp} argue that the independence assumption is crucial in 

% The deficiency of the LSE is well-known and has motivated the study of robust estimators such as the least absolute deviation estimator \cite{MR1128411} or \cite[Page 336]{VdVW96} and medians-of-mean estimators~\cite{lugosi2016risk}. The optimality of the above result is an open problem ??. 

% The dependence between $\epsilon$ and $X$ is indeed important. \cite[Proposition~3]{han2017sharp} shows that 
% It turns out that the, LSE can have a better rate of convergence when $\F$ has additional local smoothness (around $f_0$); . The following theorem (proved in Section~\ref{sec:proofSecondMomentRate}) improves the rate of convergence of the LSE when $F_\delta(\cdot)$ ``shrinks'' as $\delta \downarrow 0. $
% subsection example_convex_regression (end)
Table~\ref{tab:Example_comp} shows some interesting applications of Theorem~\ref{thm:SecondMomentRate} and compares the results with Theorem~3 of~\cite{han2017sharp}. Both of these theorems consider function classes $\F$ that satisfy $\eqref{eq:Brack}.$ However, Theorem~\ref{thm:SecondMomentRate} uses the {\clr envelope growth condition} of the $\F$ when deriving the rates, while~\cite{han2017sharp} does not assume any structure on $\F$. Table~\ref{tab:Example_comp} shows that when $\F$ is class of H\"{o}lder or Sobolev functions, then the inherent smoothness of the functions involved can help significantly reduce the requirement on $\epsilon$ for the optimal $n^{-1/(2+\alpha)}$ rate of convergence when $\alpha<1$. To see this, observe that  for H\"{o}lder classes  $s = 2/(2+ \alpha)$. Thus when $\alpha<1$, we have  $2/s< 1+ 2/\alpha$, i.e., the moment requirements for Theorem~\ref{thm:SecondMomentRate}  is smaller than that in~\cite[Theorem~3]{han2017sharp}. This is significant, as in contrast to the results of~\cite{han2017sharp}, Theorem~\ref{thm:SecondMomentRate} allows for errors $\epsilon$ to depend on $X.$

 \begin{table}[!ht]
 % \vspace{-.2in}
 \center
 
      \caption[Different choices of $\F$ and their corresponding parameters] {\label{tab:Example_comp}Different choices of $\F$,  their corresponding values of $\alpha$ and $s$, and the number of moments of $\epsilon$ required for the LSE to converge at an  $n^{-1/(2+\alpha)}$ rate. The moment requirements due to~\cite{han2017sharp} are under stronger assumptions, they assume independence between $\epsilon$ and $X$.}
      \centering
      \resizebox{\textwidth}{!}{

  \begin{tabular}{lccccl}
  \toprule
  \multicolumn{4}{c}{Choices of $\F$ and $\rchi$}&\multicolumn{2}{c}{ Moments needed for an $n^{-1/(2+\alpha)}$ rate}\\
  \cmidrule(r){1-4} \cmidrule(l){5-6}
  \multicolumn{1}{c}{$\rchi$} &$\F$  & $\alpha$ & $s$ &  \cite[Theorem 3]{han2017sharp} &Theorem~\ref{thm:SecondMomentRate}\\
  \midrule
  $[0,1]^d$ &$\gamma$-H\"older class & $d/\gamma$ & $\frac{2\gamma}{2\gamma +d}$& $1+2\frac{\gamma}{d}$& $2+ \frac{d}{\gamma}$\\ 
  $[0,1]^d$ & $\gamma$-Sobolev class & $d/\gamma$ & $\frac{2\gamma-1}{2\gamma+d-1}$& $1+2\frac{\gamma}{d}$ & $2+ \frac{2d}{2\gamma-1}$\\ 
  $[0,1]$ & Uniformly Lipschitz   & {\clr 1} & 2/3& 3& 3\\
  $[0,1]$ &   $\gamma$-H\"older class $\cup\, \{\mathbf{1}_{[a,b]}: 0\le a\le b\le 1\}$ & $1/\gamma$& $0$& $1+ 2\gamma$&$\infty$\footnote{Here $s=0$, thus the upper bound on the rate of convergence of the LSE  is $n^{-1/( 2+ \alpha (q-1)/q)}.$}\\
 % $[0,1]$ &{\clr Uniformly bounded Convex}& 1/2 & 1/5 & 5 & 10\\ 
  \bottomrule
  \end{tabular}}
  \end{table}
The proof of Theorem~\ref{thm:SecondMomentRate} (in Section~\ref{sec:proofSecondMomentRate}) is an application of our peeling result, Theorem~\ref{thm:GeneralPeeling}, in conjunction with a classical maximal inequality \cite[Lemma 3.4.2]{VdVW96} for bounded empirical processes.  The maximal inequality in \cite[Lemma 3.4.2]{VdVW96}  applies only to bounded empirical process and cannot be used to control the unbounded empirical process in~\eqref{eq:NecessaryExpectation}. In contrast to the standard peeling argument~\cite[Theorem~3.2.5]{VdVW96}, Theorem~\ref{thm:GeneralPeeling} incorporates a truncation step directly into the peeling argument (see Step $1$ in the proof of Theorem~\ref{thm:GeneralPeeling}) and thus allowing us to use the classical maximal inequality in this setting. To control the unbounded remainder, we observe that it has $q$ moments and use a Markov inequality of $q$th order.  The above two steps will show that $\p(r_n \|\widehat{f}-f_0\|\ge D) = O(D^{-1})$.  To show that the probability of the tail is in fact of a smaller order, we use Talagrand's inequality~\cite[Proposition 3.1]{Gine00}.  

\begin{remark}\label{rem:Smoothness_Examples}
 In each of the examples in Table~\ref{tab:Example_comp}, the complexity $\alpha$ is well known. For  H\"{o}lder, Sobolev, and Lipschitz functions we use standard interpolation inequalities to find $s$; see e.g.,~\cite{Agmon10,chen1998sieve,nirenberg2011elliptic,shen1994convergence,VANG,VdVW96,VANG}. These works also contain other examples for which $\F$ satisfies  the assumptions of Theorem~\ref{thm:SecondMomentRate}. Also see Appendix~\ref{sec:auxiliary_results} for three new interpolation inequalities. 
\end{remark}

\subsection{Example 1: Univariate convex regression} % (fold)
\label{sub:univConvex}

% To illustrate the use of Theorem~\ref{thm:SecondMomentRate}, we consider the case of univariate convex regression. 
We now find the rate of convergence of the convex LSE under heteroscedastic and heavy-tailed  errors.
 Let $\F$ be the class of convex functions on $[0,1]$ and $P_X$ be the uniform distribution on $[0,1]$. Recall that $\widehat{f}$  is only well-defined at the data points $\{X_i\}_{i=1}^n$; see~\eqref{eq:L_2Loss}. In this paper, we consider the canonical extension of $\widehat{f}$, and define $\widehat{f}$ to be the unique left-continuous piecewise linear function on $[0,1]$ with potential kinks at the data points. We are interested in finding the rate of convergence of $\widehat{f}$ when $f_0\in \F$. The class of convex functions in $[0,1]$ is unbounded. However, a simple modification\footnote{\label{foot:footnote1}\cite{han2018robustness} assume that $\epsilon$ is independent of $X$. However, their proof (\cite[Section 5.3.1]{han2018robustness}) goes through if we use the Etemadi's maximal inequality~\cite[Proposition 1.1.2]{MR1666908} and the fact that $\epsilon_i$'s satisfy~\eqref{eq:cvar} instead of L\'{e}vy's inequality for sums of i.i.d random variables~\cite[Theorem 1.1.5]{MR1666908}.} of  \cite[Lemma~5]{han2018robustness} shows that $\|\widehat{f}\|_{\infty} =O_p(1)$ if $\epsilon$ satisfies~\eqref{eq:cvar}. Let  $\F_{n} := \{f\in \F: \|f\|_{\infty}\le C \sqrt{\log n}\}$, where $C$ is a constant. Because $\|\widehat{f}\|_{\infty} =O_p(1)$, we have that  $\p(\widehat{f}\notin \F_{n}) =o(1).$ Now define $\tilde{f} := \argmin_{f \in \F_{n}} \sum_{i=1}^{n} (Y_i - f(X_i))^2,$  then $\p(\tilde{f} = \widehat{f})= 1-o(1).$ Thus the rate of convergence of $\tilde{f}$ coincides with the rate of convergence of $\widehat{f}$, because for every $D>0$
\begin{equation}\label{eq:tilde_hat}
\p(r_n \|\widehat{f}-f_0\|> D)\le \p(r_n \|\tilde{f}-f_0\|> D)+ \p(\tilde{f} \neq \widehat{f}).
\end{equation}
% where the $o(1)$ term does not depend on $D$.
 If $\epsilon$ is uniformly sub-Gaussian or bounded, then classical results \cite[Section 3.4.3.2]{VdVW96} show that $\tilde{f}$ converges at a rate $n^{-2/5}$ up to a $\log n$ factor. In this example,  we will show that the light tail assumption is unnecessary and that Theorem~\ref{thm:SecondMomentRate} implies that $\tilde{f}$ converges at a rate $n^{-2/5}$ (up to a polynomial in $\log n$ factors) if $\epsilon$ satisfies~\eqref{eq:cvar} and $\E(|\epsilon|^3|X)$ is uniformly bounded. 

\begin{figure}[h!]
\centering
\includegraphics[width=.9\textwidth]{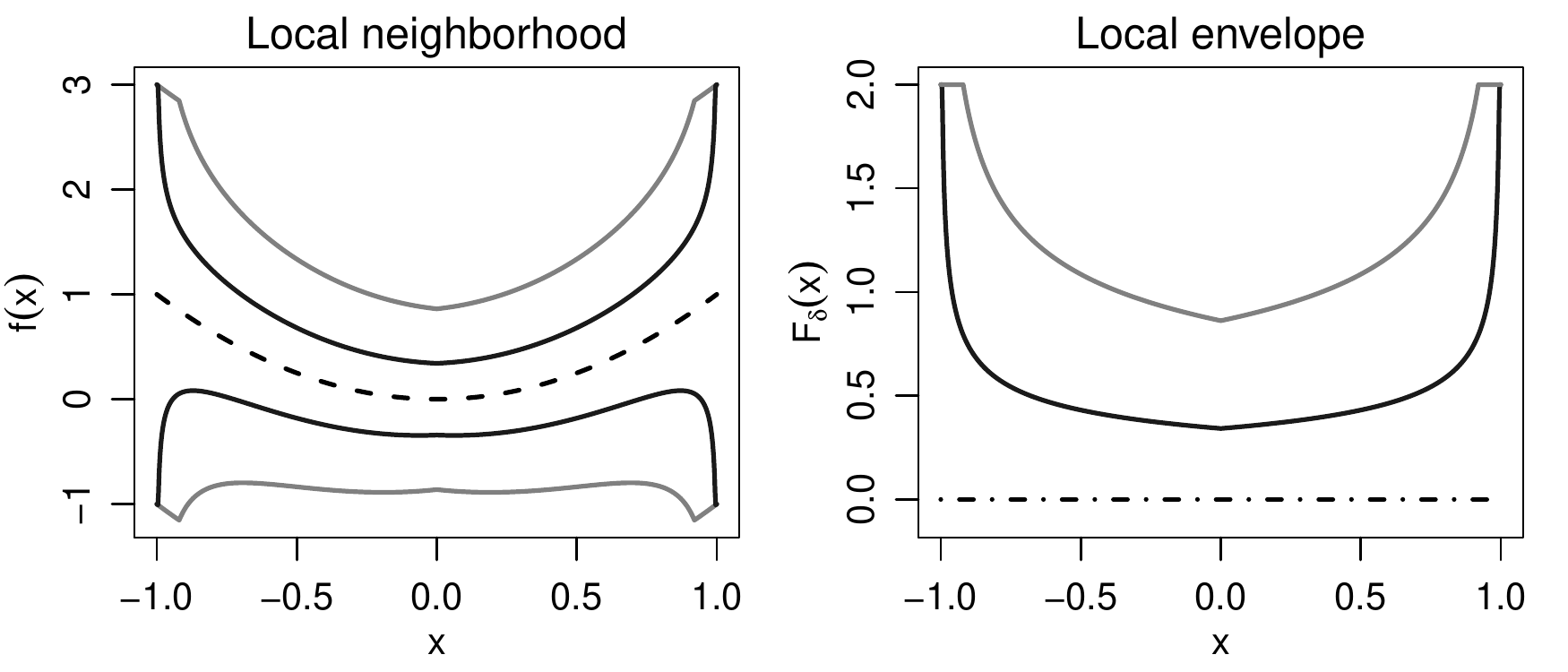}
%% scale=.8
  \caption[]{Illustration of $\F_\delta$ (left panel) and $F_\delta$ (right panel) when $f_0(x)= x^2$  and $\F:= \{f :[0,1] \to \R \,|\, \|f\|_{\infty} \le 2 \text{ and } f \text{ is convex}\}$ for  $\delta=.2$ (solid black) and $\delta=.05$ (solid gray).  Any convex function $f$ that is uniformly bounded by 2 and satisfies $\|f-f_0\| \le .2$ lies in the band created by the solid  gray  lines. The dashed line in the left panel is $f_0$.}
  \label{fig:convexEnve}
\end{figure}

 Theorem~3.1~of~\cite{2015arXiv150600034D} shows that $\F_n$ satisfies~\eqref{eq:Brack} with $A= C (\log n)^{1/4}$ and $\alpha =1/2$. Further, if $\F_{n, \delta}:= \{ f-f_0 : f\in\F_n, \|f-f_0\| \le \delta\}$ and $F_{n, \delta}(\cdot):= \sup_{g \in \F_{n, \delta}} |g(\cdot)|$, in Proposition~\ref{lem:Convex_F_delta}, we show that 
\begin{equation}\label{eq:gciahguiqj}
 \|F_{n,\delta}\|_{\infty} =C\sqrt{\log n}\quad \text{ and }\quad \|F_{n,\delta}\|_3 \le 4\delta^{2/3} \left[\log(1/\delta)\right]^{1/3}\sqrt{ \log{n}};
 \end{equation} 
 see Fig.~\ref{fig:convexEnve} for a plot of the local neighborhood and the local envelope. Suppose there exists a constant $C$ such that $\E(|\epsilon|^3|X) \le C$ for a.e. $P_X$. Because $\Phi=C\sqrt{ \log{n}}$, we have 
\begin{equation}\label{eq:Envelope_Example_convex}
\big\|(|\epsilon| +\Phi) F_{n, \delta}(X)\big\|_3 \le C  \delta^{2/3} \log(1/\delta)^{1/3}\log{n}. 
\end{equation}
% Thus $\epsilon$ and $F_{n,\delta}$ satisfy~\eqref{eq:GeneralEnvelopeCondition_L2} up to a polynomial factor of $\log(1/\delta)$ with $q=2$ and $s=2/3$.
Thus $\tilde{f}$, $\epsilon$, and $\F_n$ satisfy the assumptions of Theorem~\ref{thm:SecondMomentRate} and item~\ref{item:extraLogfactor} of Remark~\ref{rem:assumptions} (also see~Section~\ref{sub:additional_log_factors}) with $\alpha=1/2$, $s=2/3$, $\nu= 1/3$, $\Phi=C \sqrt{ \log{n}}$, and  $A= (\log n)^{1/4}$. Hence by~\eqref{eq:ep_n_extra_log}, a modification of~\eqref{eq:rn_SecondMoment}, we have that both $\tilde{f}$ and $\widehat f$ converge at a rate $n^{-2/5}\log n$ when  $\E(|\epsilon|^3|X) \le C.$  This result seems to be new. 
\section[1]{Rates of convergence of the LSE using the $L_\infty$-entropy}\label{sec:MajorGeneral}
Although~\eqref{eq:Brack} is the most widely used notion of complexity, often function classes also satisfy the stronger entropy condition~\eqref{eq:inf}, especially when $\rchi$ is bounded.
% \footnote{A counter example is the class univariate convex functions on $[0,1]$. They satisfy~\eqref{eq:Brack} with $\alpha=1/2$ but do not satisfy~\eqref{eq:inf}.}
 Moreover, they often satisfy both~\eqref{eq:Brack} and~\eqref{eq:inf} for the same value of the complexity parameter, e.g., H\"{o}lder and Sobolev functions on $[0,1]^d$. The following result (proved in Section~\ref{sec:proof:L_inf_smooth}) shows that the rate of convergence of the LSE  in Theorem~\ref{thm:SecondMomentRate} can be improved if $\F$ satisfies~\eqref{eq:inf}.

\begin{thm}\label{thm:L_inf_smooth}
Suppose $\F$ satisfies~\eqref{eq:inf}, $\epsilon$ satisfies~\eqref{eq:cvar} and~\eqref{eq:MomentCondition}, and $f_0\in \F$.  Let $\Phi:= \sup_{f\in \F} \|f\|_{\infty}$.  Moreover, suppose there exists a constant $C>0$ such that 
\begin{equation}\label{eq:GeneralEnvelopeCondition_L_infty_linf}
\| F_{\delta}\|_{\infty} \le C \Phi^{1-s} \delta^{s},
\end{equation}
 for some $s\in[0,1]$, and let \small{\begin{equation}\label{eq:rate_L_inf_L_2}
r_n := \min\left\{\frac{(nA^{-1})^{1/(2+\alpha)}}{(\sigma + \Phi)^{2/(2+\alpha)}},
% \frac{n^{1/2}}{(\sigma + \Phi)},
\frac{n^{(q-1)/(q(2-s))}}{\Phi^{2/(2-s)}},
\frac{(nA^{-1})^{(q-1)/(q(2-s)+\alpha s(q-1))}}{\Phi^{(q(2-s)+\alpha (s-1)(q-1))/(q(2-s)+\alpha s(q-1))}}\right\}.
 \end{equation}}
   Then, there exists a constant $C > 0$ depending only on $q,s,$ and $\alpha$, such that
\begin{equation}\label{eq:Tailbpound_L_nf_L_inf}
\mathbb{P}\left( r_n \|\widehat{f} - f_0\| \ge D\right)\le C D^{-q + \mathbf{1}\{s=1\}/10}
\end{equation}
  for any $n \ge 1$ and $D> 0$.

 \end{thm}
% The assumptions of the above theorem deserve some comments.
{\clr Assumption~\eqref{eq:GeneralEnvelopeCondition_L_infty_linf} of Theorem~\ref{thm:L_inf_smooth} is an $L_\infty$-envelope growth condition, cf. the $L_2$ version in~\eqref{eq:GeneralEnvelopeCondition_L2}.}  
% Assumption~\eqref{eq:GeneralEnvelopeCondition_L_infty_linf} of Theorem~\ref{thm:L_inf_smooth} on $F_\delta$ is stronger than~\eqref{eq:GeneralEnvelopeCondition_L2} {\clr and is a version of envelope growth condition}.  
Just as in Theorem~\ref{thm:SecondMomentRate}, the tail bound in~\eqref{eq:Tailbpound_L_nf_L_inf} holds for all $n\ge 1$ and the discussion in Remark~\ref{rem:Tailbound} applies to~\eqref{eq:Tailbpound_L_nf_L_inf} as well.  
 % When $\Phi:= \sup_{f\in \F} \|f\|_{\infty} < \infty$ and $\epsilon$ satisfies~\eqref{eq:MomentCondition},  then observe that every uniformly bounded function class $\F$ satisfies~\eqref{eq:GeneralEnvelopeCondition_L_infty_linf} with $s=0.$
The following corollary finds the rate of the LSE if $\F$ does not satisfy any {\clr envelope growth} assumption of~Section~\ref{sub:local_structure_of_} (i.e., $s=0$) and $A$ and $\Phi$ are constants; cf. Corollary~\ref{cor:L2_classical}.

\begin{cor}\label{cor:MajorGeneral} 
Suppose $f_0\in \F$ and $\Phi:= \sup_{f\in \F} \|f\|_{\infty}$ is a constant.
% Moreover, suppose  $\epsilon$ satisfies~\eqref{eq:cvar} and~\eqref{eq:MomentCondition} and $\F$ satisfies~\eqref{eq:Brack}.  Then for any $n \ge 1$ and $D > 0$, we have
Moreover, suppose $\epsilon$ satisfies~\eqref{eq:cvar} and~\eqref{eq:MomentCondition} and $\F$ satisfies~\eqref{eq:inf}, and let 
\begin{equation}\label{eq:rn_Major}
r_n:= \min \left\{  {n^{1/(2+\alpha)}},  {n^{1/2-1/2q}} \right\}.
\end{equation}
  Then, there exists a constant $C>0$  depending only on $q, \alpha, \sigma, A, \Phi,$ and $K_q$, such that $\mathbb{P}\big(r_n\|\widehat{f} - f_0\| \ge D\big) \le {C}{D^{-q}}$  for any $n \ge 1$ and $D > 0$.
  % , we have that $\mathbb{P}\big(r_n\|\widehat{f} - f_0\| \ge D\big) \le {C}{D^{-q}},$
% \begin{equation}\label{eq:Tail_boundLSE}
% \mathbb{P}\left(r_n \|\widehat{f} - f_0\| \ge D\right)\le  \frac{C}{D^{q}},
% \end{equation}
% where $C$ is a constant depending only on $q, \alpha, \sigma, A, \Phi,$ and $K_q$ and
% \todo{Change this?}
\end{cor}

To prove Corollary~\ref{cor:MajorGeneral}, apply Theorem~\ref{thm:L_inf_smooth} with $F_\delta\equiv 2\Phi$, i.e., $s=0$. If $\F$ is such that $s=0$  and satisfies~\eqref{eq:inf}   then Corollary~\ref{cor:MajorGeneral} uses the stronger entropy condition ($ N_{[\,]}(\zeta,\F, \|\cdot\|) \le  N(\zeta,\F, \|\cdot\|_{\infty})$) to show that the LSE converges at an $n^{-1/(2+\alpha)}$ rate under heteroscedastic errors if $q \ge 1+2/\alpha$; compare this to the rate of the LSE obtained in~Corollary~\ref{cor:L2_classical}. It is well known, that the worst case rate for the LSE under only the entropy assumption~\eqref{eq:inf} is~$n^{-1/(2+\alpha)}$ when $\epsilon$ is uniformly sub-Gaussian. Corollary~\ref{cor:MajorGeneral} shows that the heavy-tailed (and heteroscedastic) nature of the $\epsilon$ does not affect this rate as long as~$\epsilon$ has at least $1+2/\alpha$ moments and satisfies~\eqref{eq:cvar}.

Theorem~\ref{thm:L_inf_smooth} shows that the  upper bounds on the rate of convergence of the LSE in~\eqref{eq:rn_Major} can be reduced if $\F$ satisfies the {\clr envelope growth assumption}~\eqref{eq:GeneralEnvelopeCondition_L_infty_linf}.
If $A$ and $\Phi$ are constants and $\F$ satisfies~\eqref{eq:GeneralEnvelopeCondition_L_infty_linf}, we can ignore the middle term in~\eqref{eq:rate_L_inf_L_2} because ${q(2-s)} \le {q(2-s)+ s(q-1)}$ for all $q\ge 2$. Hence,  Theorem~\ref{thm:L_inf_smooth} implies that the  LSE converges at the rate
\begin{equation}\label{eq:e_nfinal_L_inf_L_2}
\max \left\{{n^{-1/(2+\alpha)}} ,  {n^{-(q-1)/(q(2-s)+\alpha s(q-1))}} \right\}.\end{equation}
 This implies if $ q\ge{(2 + \alpha(1-s))}/{(s+ \alpha(1-s))}$, then the LSE converges at an $n^{-1/(2+\alpha)}$ rate. Furthermore
\[\frac{2 + \alpha(1-s)}{s+ \alpha(1-s)} \le 1+ \frac{2}{\alpha}\text{ for all } s\le 1\text{ and }0 \le \alpha<2.\]
{ Thus if $s>0$, then Theorem~\ref{thm:L_inf_smooth} shows that the LSE converges at an $n^{-1/(2+\alpha)}$ rate under weaker assumptions on $\epsilon$ than in Corollary~\ref{cor:MajorGeneral}.}

The proofs of Theorems~\ref{thm:SecondMomentRate} and~\ref{thm:L_inf_smooth} are similar. 
But in the case of Theorem~\ref{thm:SecondMomentRate}, we could apply the readily available maximal inequality~\cite[Lemma 3.4.2]{VdVW96}. Existing maximal inequalities, however,  cannot take the $L_\infty$-covering number into account. For this purpose, we use  a chaining argument~\cite{Dirksen} in conjunction with a new maximal inequality for the maximum over a finite set; see Proposition~\ref{lem:bernouli_maximal}.
% {\clr In the same context Theorem 7.1 of~\cite{2017arXiv170800145K}, we have applied (2.5.5) of~\cite{VdVW96} which is obtained from Bernstein inequality and the chaining based on this inequality can lead to rates only in case $\alpha < 1$. Proposition~\ref{lem:bernouli_maximal}, however, leads to rates over all $\alpha\in[0, 2)$.} 

Proposition~\ref{lem:bernouli_maximal} is also of independent interest and compares favorably to Lemma 8 of~\cite{MR3350040}. 
Our result shows that the maximum of $N$ centered averages converges at the rate of $\sqrt{n^{-1}\log N}$, if $\log N =o(n)$ and the envelope has finite $q\ge 2$ moments; see~\eqref{eq:q-moment}. On the other hand,~\cite[Lemma 8]{MR3350040} requires $\log N = O(n^{1-2/q})$ for a $\sqrt{n^{-1}\log N}$ rate of convergence; see Section~\ref{rem:B1Application} in Appendix~\ref{sec:a_new_maximal_inequality_for_finite_maximums} for details.  

% \begin{remark}\label{rem:Comp_3.1}
% If $f_0,\, \F,$ and $\epsilon$ satisfy the assumptions of Theorem~\ref{thm:L_inf_smooth}, then they also satisfy the assumptions of Theorem~\ref{thm:SecondMomentRate}. But the rate of convergence of the LSE in Theorem~\ref{thm:L_inf_smooth} (\eqref{eq:rate_L_inf_L_2} and~\eqref{eq:e_nfinal_L_inf_L_2}) is faster than that in Theorem~\ref{thm:SecondMomentRate} (\eqref{eq:rn_SecondMoment} and~\eqref{eq:rn_final_L2}, respectively). Also note that $(2+ \alpha(1-s))/(s+ \alpha(1-s))$ is strictly smaller than $2/s$. 
% \end{remark}
\subsection{Example 2: Multivariate and multiple index smooth regression models} % (fold)
\label{sub:multi_index} 
In this section, we consider the example of multivariate regression when the unknown function is known to be smooth. Let $\rchi :=[0,1]^d$ and $P_X$ be the uniform distribution on $\rchi$; the unit hypercube can be replaced a bounded and convex subset of $\R^d$.
 For any vector $k = (k_1, \ldots, k_d)$ of $d$ positive integers, define the differential operator
$ D^k := \frac{\partial^{\|k\|_1}}{\partial x_1^{k_1}\cdots\partial x_d^{k_d}},$  where $\|k\|_1 := \sum k_i$. Define the class of real valued functions on $\rchi$
  \[
  \mathcal{F}_{\gamma,d} := \left\{f:\,\max_{\|k\|_1 \le [\gamma]}\sup_{x}|D^kf(x)| + \max_{\|k\|_1 = [\gamma]}\sup_{x\neq y}\frac{|D^kf(x) - D^kf(y)|}{|x - y|^{\gamma - [\gamma]}} \le 1\right\}.
  \]
  % (The constant $1$ in $\mathcal{F}_{\gamma, d}$ represents the bound on the smoothness norm.) 
  Theorem 2.7.1 of~\cite{VdVW96} implies that there exists a constant $C$ depending only on $\gamma$ and $d$ such that 
  \[
  \log N(\nu, \mathcal{F}_{\gamma,d}, \norm{\cdot}_{\infty}) \le C\nu^{-d/\gamma}\quad\mbox{for all}\quad \nu > 0.
  \]
  
   % and $\Lambda(\mathcal{X}^{1})$ denotes the Lebesgue measure of the set $\{y:\,|x - y| \le 1\mbox{ for some }x\in \mathcal{X}\}$. 

  \paragraph{Multivariate smooth regression} Suppose $f_0 \in \F_{\gamma,d}$ for some $\gamma >d/2$ and $\widehat{f}$ is  defined as in~\eqref{eq:L_2Loss} with $\F=\F_{\gamma,d}$.   
  % We will apply Theorem~\ref{thm:L_inf_smooth} to show that $\widehat{f}$ is minimax rate optimal. 
  By~\cite[Lemma~2]{chen1998sieve}, we have that $\F_{\gamma, d}$ satisfies~\eqref{eq:GeneralEnvelopeCondition_L_infty_linf} with $s={2\gamma}/(2\gamma +d)$.  Theorem~\ref{thm:L_inf_smooth} and~\eqref{eq:e_nfinal_L_inf_L_2} imply that
\begin{equation}\label{eq:Tgif}
n^{\gamma/( 2\gamma+d)}\|\widehat{f} - f_0\| = O_p(1), \quad \text{if} \quad q\ge 2+ \frac{2 d \gamma-d^2}{2\gamma^2 +d^2}.
\end{equation}
Because~\cite[Theorem 3.2]{Gyorfi02} shows that $n^{-\gamma/( 2\gamma+d)}$ is the global minimax rate of convergence,~\eqref{eq:Tgif} show that LSE is minimax optimal under~\eqref{eq:cvar} and~\eqref{eq:MomentCondition} if $q \ge 2+ {(2 d \gamma-d^2)}/{(2\gamma^2 +d^2)}$.   Note that $\F_{\gamma,d}$ also satisfies the assumptions of~Theorem~\ref{thm:SecondMomentRate} but it would lead to a sub-optimal result; see Table~\ref{tab:Example_comp}.

  \paragraph{Multiple index smooth regression} In the above setup the rate of convergence of the LSE is strongly affected by the dimension. When $d$ is large, a widely used semiparametric alternative that ameliorates the curse of dimensionality is the multiple index model~\cite{MR1891738,2017arXiv170800145K}. In multiple index model, the true regression function is assumed to belong to 
  \[
\M_{\gamma, d, d_1}:= \left\{x\mapsto f(Bx):\,f\in\F_{\gamma,d_1}\mbox{ and } B\in\mathbb{R}^{d_1\times d}\mbox{ satisfying }\norm{B}_{2} \le 1\right\},
  \]
  for some $d_1 \le d$. Note that {\clr the global} minimax optimal rate in the multiple index model is $n^{-\gamma/(2 \gamma+ d_1)}$~\cite[Page 129]{MR1272079}, a faster rate than that in~\eqref{eq:Tgif}. We now find a sufficient condition under which the LSE achieves the minimax rate.  

  Because $\norm{B}_{2} \le 1,$ it can be easily shown that  there exists a constant $C$ (depending only on $d$) such that 
  \[  \log N(\nu, \M_{\gamma, d, d_1}, \norm{\cdot}_{\infty}) \le C\nu^{-d_1/\gamma}\quad\mbox{for all}\quad \nu > 0.
\]
 Thus $\M_{\gamma, d, d_1}$ is much less ``complex'' than $\F_{\gamma, d}$, when $d_1$ is smaller than $d$.  By Proposition~\ref{prop:Interpolation_for_multiple_index_models}, we have that $\M_{\gamma, d, d_1}$ satisfies~\eqref{eq:GeneralEnvelopeCondition_L_infty_linf} with $s = 2\gamma/(2\gamma + d_1).$ Now, Theorem~\ref{thm:L_inf_smooth} shows that if 
  \[
  \widehat{f}:= \argmin_{f  \in \M_{\gamma, d, d_1}}\sum_{i=1}^{n} (Y_i - f(X_i))^2,
  \]
  then
  \[  n^{\gamma/(2\gamma+d_1)}\|\widehat{f} - f_0\| = O_p(1), \quad \text{when} \quad q\ge 2 + \frac{2 d_1\gamma - d_1^2}{2\gamma^2 + d_1^2}.\]

  \paragraph{Additive model regression} An even simpler function class than $\M_{\gamma, d, d_1}$ is given by functions that are separable in their coordinates~\cite{MR994249,FriedmanStuetzle81}. Formally, define
  \[
  \mathcal{A}_{\gamma} := \Big\{x\in\mathbb{R}^d\mapsto f(x) = \sum_{j=1}^d f_j(x_j):\,f_j\in\mathcal{F}_{\gamma,1}\Big\}.
  \]
  In this case it can be shown that there exists a constant $C$ such that 
  \[
  \log N(\nu, \mathcal{A}_{\gamma}, \|\cdot\|_{\infty}) \le C\nu^{-1/\gamma}\quad\mbox{for all}\quad \nu > 0.
  \]
 By Proposition~\ref{prop:Interpolation_for_Additive_Models}, it follows that $\mathcal{A}_{\gamma}$ satisfies~\eqref{eq:GeneralEnvelopeCondition_L_infty_linf} with $s = 2\gamma/(2\gamma+1)$. Thus Theorem~\ref{thm:L_inf_smooth} shows that if
  \[
  \widehat{f} := \argmin_{f\in\mathcal{A}_{\gamma}}\sum_{i=1}^n (Y_i - f(X_i))^2,
  \]
  then {\clr for $\gamma > 1/2$,}
  \[ n^{\gamma/(2\gamma + 1)}\|\widehat{f} - f_0\| = O_p(1)\quad\mbox{ when }\quad q \ge 2 + \frac{2\gamma - 1}{2\gamma^2 + 1}.\]

  The function classes above can also be replaced by other smoothness classes such as Sobolev, Nikolskii, or Besov spaces \cite{nickl2007bracketing,2020arXiv200810979G}. Furthermore, using the proofs of Propositions~\ref{prop:Interpolation_for_Additive_Models} and~\ref{prop:Interpolation_for_multiple_index_models}, one can also consider combination of function spaces $\M_{\gamma, d, d_1}$ and $\mathcal{A}_{\gamma},$ wherein some coordinates are modeled through linear combinations and the remaining coordinates are modeled through additive model.
% subsection subsection_name (end)

% section back_to_ (end)

\section{Rate of convergence for VC-type classes} % (fold)
\label{sec:rate_of_convergence_for_vc_type_classes}

In Sections~\ref{sec:L2all} and~\ref{sec:MajorGeneral}, we showed that the local envelope $F_\delta$ can affect the rate of convergence of the LSE. We showed that if $\F$ satisfies~\eqref{eq:Brack} or~\eqref{eq:inf} with $\alpha$ and the envelope growth parameter  is non-zero, then the LSE converges at a rate $n^{-1/(2+\alpha)}$ even when $\epsilon$ has only few moments. If $\F$ is the class of {\clr totally bounded} smooth functions (e.g., Sobolev, H\"{o}lder, or Besov spaces),  then $s$  depends only on the smoothness of the functions in the class and not on the choice of $f_0$; recall that $\F_\delta$ is the local neighborhood of $f_0$ in $\F$. However, it turns out that for certain function classes $\F$, the {\clr envelope growth parameter $s$} can depend on $f_0.$ For example, in Proposition~\ref{lem:Convex_F_delta}, we show that if $\F$ is the class uniformly bounded convex functions on $[0,1]$, then   $F_\delta(x) \le C\Phi^{1/3} \delta^{2/3}\max\{x^{-1/3}, (1 - x)^{-1/3}\}$ for any $f_0\in \F$. But if $f_0$ is a linear function (or piecewise linear) then~\cite[Lemma A.3]{guntuboyina2015global} shows that $F_\delta(x) \le C \delta \max\{x^{-1/2}, (1 - x)^{-1/2}\}$, i.e., $F_\delta$ has a smaller $L_2$-norm when $f_0$ is linear. This change in local behavior of $\F_\delta$, when $f_0$ belongs to a particular subclass of functions, drives the adaptive behavior of the LSE in shape-constrained regression; see e.g.,~\cite{bellec2018sharp,chatterjee2015risk,chatterjee2015adaptive,guntuboyina2018nonparametric} and references therein. Furthermore, in these examples it turns out that  $\F$ satisfies~\eqref{eq:Unif_entr}, when $f_0$ belongs to these special subclasses of $\F$. In the following theorem (proved in Section~\ref{sec:Vc-type-proof}) we use the envelope growth condition to find the worst-case rate of convergence of the LSE when $\F$ satisfies~\eqref{eq:Unif_entr} around $f_0$ and $\epsilon$ satisfies~\eqref{eq:cvar}. {\clr Theorem \ref{thm:Vc-type} can be used to prove that the LSE can attain the local minimax rate of convergence in the sense of \eqref{eq:local-minimax-rate}.} In this section, we do not make any assumptions on the higher order moments of $\epsilon$. This is done with the goal of keeping the result simple. Furthermore, it turns out that LSE is rate optimal in certain scenarios with just two finite moments.
\begin{thm}\label{thm:Vc-type}
Suppose $\F$ satisfies~\eqref{eq:Unif_entr}, $\epsilon$ satisfies~\eqref{eq:cvar}, and  $f_0 \in \F$. Let $\Phi:= \sup_{f\in \F} \|f\|_{\infty}$.  Assume that $\sigma,$ $\Phi$, and $A$ (in~\eqref{eq:Unif_entr}) are constants. Moreover, suppose  there exists a constant $C>0$  such that 
\begin{equation}\label{eq:GeneralEnvelopeCondition_VC}
\big\|  F_{\delta}\big\| \le C  \Phi^{1-s} \delta^{s},
\end{equation}
for some $s\in[0,1]$, and for $\alpha$, $\beta$ as in~\eqref{eq:Unif_entr}, let
\begin{equation}\label{eq:rn_vc_combined}
r_n := 
\begin{cases}
n^{1/(2(2-s))} &\mbox{if }\alpha\in[0, 2)\mbox{ and }\beta \ge 0,\\
(n^{1/2}/\log n)^{1/(2 - s)} &\mbox{if }\alpha = 2\mbox{ and }\beta = 0,\\
n^{1/(\alpha(2 - s))} &\mbox{if }\alpha > 2\mbox{ and }\beta = 0,
\end{cases}
\end{equation}
 Then, there exists a constant $C > 0$ depending only on $\sigma, \Phi, A,$ and $\alpha$ such that
% for some constant $C$, then  there is a constant  $C$ depending only on $ \sigma$  and $\alpha$ such that 
\[
\p\left( r_n \|\widehat{f} - f_0\|\ge D \right) \le {C}{D^{-4 (2-s)/3}}
\] for any $n \ge 1$ and $D> 0$.
\end{thm}

\begin{remark}\label{rem:assumptions_VC}
We make some observations about the assumptions and  conclusions of Theorem~\ref{thm:Vc-type}.
\begin{enumerate}

  \item The assumption~\eqref{eq:GeneralEnvelopeCondition_VC} has a different structure than those in~Theorems~\ref{thm:SecondMomentRate} and~\ref{thm:L_inf_smooth}. In contrast to assumption~\eqref{eq:GeneralEnvelopeCondition_L_infty_linf}, we only require a bound on $\|F_\delta\|$. No control is required for higher moments of $\epsilon$.
  \item The above theorem provides the rates of convergence of the LSE even when $\F$ is non-Donsker, i.e., $\alpha\ge 2$.

  \item  The assumption that $A$ and $\Phi$ do not change with $n$ is made  to keep the presentation simple. In the proof of the result (Section~\ref{sec:Vc-type-proof}), we provide explicit finite sample tail bounds that allow $A$ and $\Phi$ to depend on $n$. See~\eqref{eq:en_vC_case1},~\eqref{eq:en_vC_case2}, and~\eqref{eq:en_vC_case_3} to find the exact relationship between $r_n$, $\Phi$, and $A$ for the three situations considered in~\eqref{eq:rn_vc_combined}.

\item If $s=1$, then it is clear that the obtained rate of convergence, $n^{-1/2}$, of the LSE cannot be improved when $\alpha \in [0,2)$ and $\beta \ge 0$.

\item  If $s<1$ and $\epsilon$ satisfies higher moment assumptions, then
the rates obtained in~\eqref{eq:rn_vc_combined} can be improved  using the tools developed in this paper along with Theorem 2.1 of~\cite{van2011local}.\label{item:Moremoments}

\item 
Just as in Theorems~\ref{thm:SecondMomentRate} and~\ref{thm:L_inf_smooth}, we can relax the bound on $\|F_\delta\|$ in~\eqref{eq:GeneralEnvelopeCondition_VC} to be of the form $C \Phi \delta^s (\log(1/\delta))^\nu$. This relaxation will reduce $r_n$ in~\eqref{eq:rn_vc_combined} by a polynomial  in $\log n$ factor; the order  of the polynomial will depend only on $s, \nu$, and $\alpha$. In particular, if $\alpha\in [0,2)$ and $\beta\ge 0$, then the rate of convergence of the LSE is $n^{-1/(2(2-s))}(\log n)^{\nu}$ if  $\|F_\delta\| \le C \Phi \delta^s (\log(1/\delta))^\nu.$ \label{item:extraLogfactor_VC}

\item To prove Theorem~\ref{thm:Vc-type}, we use the refined Dudley's chaining inequality (Lemma~\ref{lem:UniformCoveringPhi}) in conjunction with Theorem~\ref{thm:GeneralPeeling}.
\end{enumerate}
\end{remark}

% section rate_of_convergence_for_vc_type_classes (end)
\begin{remark}[Comparison with Theorem~1 of~\cite{han2018robustness}]\label{rem:Rboustness}
Theorem~\ref{thm:Vc-type} is an improvement over Theorem~2 of~\cite{han2018robustness}. Our result allows $\epsilon$ to depend on $X$ and $\epsilon$ to have only $2$-moments, whereas~\cite[Theorem~1]{han2018robustness} requires the error to have $L_{2,1}$ moments (i.e., $\int_{0}^{\infty} \sqrt{\p(|\epsilon|> t)}\, dt <\infty$) and be independent of $X.$ Finite $L_{2,1}$ moments  imply finite second moments. But the converse is not true. Our proof of Theorem~\ref{thm:Vc-type} uses our new peeling result (Theorem~\ref{thm:GeneralPeeling}) and is different from the proof of~\cite[Theorem~1]{han2018robustness}. {\clr Furthermore, Theorem~\ref{thm:Vc-type} and discussion in Section~\ref{sub:univConvex} (see footnote~\ref{foot:footnote1}) can be used to establish oracle inequalities for both convex and isotonic LSEs under heteroscedastic errors with only two moments. The proves will be almost identical to (but improve upon) Theorems~3 and 5 of ~\cite{han2018robustness}. }
% and establish the rate adaptive behavior of the LSE~\cite{bellec2018sharp,chatterjee2015risk,han2018robustness} when $\E(\epsilon^2|X)\le \sigma^2$ and $\epsilon$ is  allowed to depend on $X.$ 
% {\clr Using the tools developed in Theorem~\ref{thm:Vc-type}, we can establish oracle inequalities for both convex and isotonic LSEs under heteroscedastic errors with only two moments. The proofs are almost identical to the proof of Theorems 3 and 5 of~\cite{han2018robustness} and will not be repeated here; we can use techniques in the proof of Theorem~\ref{thm:Vc-type} to bound the multiplier process in Eq. (4.10) of~\cite{han2018robustness} (without assuming independence between $\epsilon$ and $X$).}

\end{remark}

\subsection{Example 3: Univariate isotonic regression} % (fold)
\label{sub:example_3_univariate_isotonic_regression}
{\clr Let $\F$ be the set of nondecreasing functions on $[0,1]$ and $P_X$ be any nonatomic probability measure on $[0,1]$. In a fixed design setting,~\cite{zhang2002risk} shows that  $\|\widehat{f}-f_0\|_n= O_p( n^{-1/3})$ when $\epsilon$ has finite variance. Moreover, when $f_0\equiv 0$ (or any other constant), the LSE satisfies $\|\widehat{f}-f_0\|_n= O_p(\sqrt{\log n/n})$. In this case, the LSE is near minimax rate optimal; when $f_0\equiv 0$ the local minimax rate of convergence is $n^{-1/2}$~\citep{gao2020estimation}.~\cite{han2018robustness}  establish the above rates in the random design setting when $\epsilon$ has finite $L_{2,1}$ moment and is independent of $X.$ The number of finite moments required of $\epsilon$ was unknown when $\epsilon$ is allowed to depend on $X.$ In this example,  We will show that both the independence and the finite $L_{2,1}$ moment assumption in~\cite{han2018robustness} can be removed when $f_0\equiv 0$ and~$\epsilon$ satisfies~\eqref{eq:cvar}.\footnote{If $f_0$ is non-constant, then $\F$ does not satisfy~\eqref{eq:Unif_entr} and Theorem~\ref{thm:Vc-type} does not apply.  Furthermore, in this case one can show that $s=0$. Thus an application of~Theorem~\ref{thm:SecondMomentRate} yields an $n^{-1/3}$ rate only when $q=\infty$; see Remark~\ref{rem:Rboustness} for a discussion. } 
}

Note that $\F$ is unbounded, but the discussion in Section~\ref{sub:univConvex} (see footnote~\ref{foot:footnote1} in page~\pageref{foot:footnote1}) and \cite[Lemma 5]{han2018robustness}  show that if $\epsilon$ satisfies~\eqref{eq:cvar}, then $\|\widehat{f}\|_{\infty} =O_p(1).$  Thus, following the arguments of Section~\ref{sub:univConvex},  it is easy to see that the rate of convergence of the isotonic LSE matches (up to a polynomial in $\log n$ factor) the rate of convergence of LSE  when  $\F$ is the set of nondecreasing functions uniformly bounded by $1$.~\cite[Example 3.8]{MR2243881} show that if $ f_0\equiv 0$ and $\F:=\{f: [0,1] \to [-1,1]| f \text{ is nondecreasing}\}$, then 
\begin{equation}\label{eq:F_del_monoton}
F_\delta(x) = \min\Big\{1, \delta\, \max\big(P_X[0,x], P_X[x,1]\big)^{-1/2}\Big\} \text{ for all } x\in [0,1], 
\end{equation}
where for every $0\le a\le b\le 1$, $P_X[a,b] := \p(X\in [a,b]).$ Furthermore, \cite[Example 3.8]{MR2243881} show that $\|F_\delta\| \le C \delta \sqrt{\log(1/\delta)}.$  Therefore, Theorem~\ref{thm:Vc-type} in conjunction with the arguments in item~\ref{item:extraLogfactor_VC} of Remark~\ref{rem:assumptions_VC} (also see item~\ref{item:extraLogfactor} in Remark~\ref{rem:assumptions})  implies that $\widehat{f}$ converges at an $n^{-1/2}$ rate (up to a polynomial in $\log n$ factor) when $\epsilon$ satisfies~\eqref{eq:cvar}.  Similar results exist when $\epsilon$ is independent of $X$~\cite{han2018robustness} or under the fixed design setting~\cite{chatterjee2015risk,2017arXiv170506386G,zhang2002risk}. However, to the best of our knowledge our result here is new and reduces the assumptions on $\epsilon$ for optimal convergence of the LSE. 

\begin{remark}[Extension to piecewise constant functions and adaptive rates]\label{rem:Extension}
In the above example, we showed that the isotonic LSE converges at a parametric rate (up to $\log n$ factors) when $f_0$ is a constant function. The above result can be generalized to case when $f_0$ is piecewise constant functions with $K$-pieces to show that that the 
$\|\widehat{f}-f_0\| \le C \log n\sqrt{K/n} $ with high probability; {see \cite{han2018robustness} and~\cite{chatterjee2015risk} for similar results under independence and fixed design settings, respectively. It should be noted that under the fixed design setting,~\cite{gao2020estimation} shows that the minimax rate with respect to $L_2(\mathbb{P}_n)$ loss is $\sqrt{K/n} \sqrt{\log \log (n/K)}$ when $K\ge 2$.}  Furthermore, if $f_0$ can be ``approximated well'' by a piecewise constant function then Theorem~\ref{thm:Vc-type} can be used to find sharp rate upper bounds on $\|\widehat f -f_0\|$; see \cite[Section 3 and Theorem 3]{han2018robustness} for an excellent and elaborate discussion on this. 

Also, as discussed in the beginning of Section~\ref{sec:rate_of_convergence_for_vc_type_classes},~\cite[Lemma A3]{guntuboyina2015global} shows that $F_\delta$ satisfies~\eqref{eq:F_del_monoton} when $f_0$ is a linear function and $\F$ is the class of uniformly bounded convex functions. Thus a similar almost parametric rate can be proved if $f_0$ is a linear (piecewise linear or well approximated by piecewise linear function) function on $[0,1]$ and  $\F$ is the set of convex functions on $[0,1]$.
% \todo[inline]{Add a discussion on oracle inequalities and connect to Example~\ref{sub:example_3_univariate_isotonic_regression}. }

% {\clr Using the tools developed in Theorem~\ref{thm:Vc-type}, we can establish oracle inequalities for both convex and isotonic LSEs under heteroscedastic errors with only two moments. The proofs are almost identical to the proof of Theorems 3 and 5 of~\cite{han2018robustness} and will not be repeated here; we can use techniques in the proof of Theorem~\ref{thm:Vc-type} to bound the multiplier process in Eq. (4.10) of~\cite{han2018robustness} (without assuming independence between $\epsilon$ and $X$).}
% Using the tools developed in using the tools used in our Theorem~\ref{thm:Vc-type}  Finally, we mention that following the proof of Theorems 3 and 5 of~\cite{han2018robustness} provides an oracle inequality for both convex and isotonic LSEs under heteroscedastic errors with only two moments. The proof is almost identical and will not be repeated here. The reason for this is that Eq. (4.10) in~\cite{han2018robustness} is a crucial step and the right hand side of this inequality with Rademacher variables $\varepsilon_i$ replaced by the noise $\epsilon_i$ can be bounded using the tools used in our Theorem~\ref{thm:Vc-type} under heteroscedastic noise (i.e., without assuming independence of $\epsilon$ and $X$).

\end{remark}
% subsection example_3_univariate_isotonic_regression (end)

\subsection{A lower bound for rate of convergence of the LSE} % (fold)
\label{sub:a_lower_bound_for_rate_of_convergence_of_the_lse}
In this and the previous sections, we have only presented upper bounds for the rate of convergence of the LSE under various entropy bounds and {\clr the envelope growth parameter}. One of the main messages from the previous sections is that {\clr the envelope growth parameter} dampens the effect of the number of moments $q$ of the response/errors on the rate upper bound of the LSE. We will now prove {\clr that, under heavy-tailed error, the envelope growth parameter plays a role in the rate of convergence the LSE in the worst case. We do this by constructing a function class for which the rate of convergence proved in Theorem \ref{thm:Vc-type} under \eqref{eq:GeneralEnvelopeCondition_VC} is tight.}
 % Theorem~\ref{thm:Lowerbnd} partly answers the question: \emph{Should the rate of convergence of the LSE depend on the {\clr envelope growth parameter}?} 
 This lower bound heavily borrows from the machinery developed in Proposition 2 of~\cite{han2018robustness}, but with some crucial changes.   
\begin{thm}\label{thm:Lowerbnd}
Let $P_X$ be the uniform distribution on $[0, 1]$. Fix any $s\in[0, 1]$. Then, there exists (1) a function  class $\widetilde\F$ satisfying~\eqref{eq:Unif_entr} with $\alpha\in [0,2]$ and $\beta \ge 0$; (2) a function $f_0 \in \widetilde\F$ such that the local envelope around $f_0$ satisfies~\eqref{eq:GeneralEnvelopeCondition_VC}; and (3) an error distribution for $\epsilon$ with $\|\epsilon\|_2 <\infty$, such that, for some constant $C > 0$ that depends on $s$,
\begin{equation}\label{eq:lwrbnd}
\liminf_{n\to\infty}\, \mathbb{P}\Big(\|\widehat{f}-f_0\| \ge Cn^{-1/(2(2-s))} (\log n)^{-2\mathbf{1}\{s\notin \{0, 1\}\}/(s(2-s))}\Big) \ge \frac{1}{10}.
\end{equation}
% for large enough $n.$
\end{thm}
{\clr
 % {\clg Theorem~\ref{thm:Lowerbnd} (proved in Section~\ref{sec:proof_of_theorem_thm:lowerbnd}) shows that the rate of convergence of the LSE under the assumptions of Theorem~\ref{thm:Vc-type} cannot be much faster than $n^{-1/(2(2-s))}$.}
 The above result can thought as a worst case pointwise asymptotic lower rate bound, i.e., for the function class $\widetilde{\F}$ there exists a function $f_0$ for  which LSE does not converge at rate ``much'' faster than $n^{-1/(2(2-s))}$.
}% \todo{Blowing rate, above. Cite~\cite{rakhlin2017empirical}?}
% subsection a_lower_bound_for_rate_of_convergence_of_the_lse (end)
\section{Misspecification} % (fold)
% \section{Misspecification and oracle-type inequalities} % (fold)
\label{sec:oracle_ineq}
The results in the previous sections find the rate of convergence of the LSE $\widehat{f}$ when $f_0\in\mathcal{F}$ and errors have finite number of moments. Recall that $f_0(x) = \mathbb{E}[Y|X = x]$. A crucial step in finding upper bounds for $\|\widehat{f}-f_0\|$ is proving
\begin{equation}\label{eq:Quadratic-lower-bound}
\mathbb{E}[(Y - f(X))^2] - \mathbb{E}[(Y - f_0(X))^2] = \|f - f_0\|^2,
\end{equation}
for any function $f\in\mathcal{F}$; see \eqref{eq:lower-bound-main} in the proof of Theorem~\ref{thm:GeneralPeeling}. A natural next step is the study of the LSE when $f_0\notin\F$. LSEs under misspecification have received a lot of attention but most works assume restrictive conditions on $\epsilon$; see e.g.,~\cite{birge1998minimum,han2018robustness,koltchinskii2006local,MR2829871,VdVW96}. The techniques developed in this paper can be used to relax the assumptions on $\epsilon$ and allow for heteroscedastic and heavy-tailed errors. 

 If the true conditional expectation ($f_0$) does not belong to the class $\mathcal{F}$ but $\mathcal{F}$ is a convex set, then defining
\[
\bar{f} := \argmin_{f\in\mathcal{F}}\,\mathbb{E}[(Y - f(X))^2],
\]
we have that for any $f\in\mathcal{F}$,
\[
\mathbb{E}[(Y - \bar{f}(X))(f - \bar{f})(X)] ~{\clr\le}~ 0.
\]
This implies that
$\mathbb{E}[(Y - f(X))^2] - \mathbb{E}[(Y - \bar{f}(X))^2]\, {\clr\ge}\, \|f - \bar{f}\|^2.$
Using this fact instead of~\eqref{eq:Quadratic-lower-bound}, the results proved in previous sections imply the same rate bounds for $\|\widehat{f} - \bar{f}\|${\clr; also, see Theorem 3.2.5 of~\cite{VdVW96}}. Thus, when $\F$ is a convex set, the proofs of Theorems~\ref{thm:SecondMomentRate} and~\ref{thm:L_inf_smooth} go through by replacing $\epsilon$ with $\xi = Y - \bar{f}(X)$; note that $\E(\xi|X)\neq0$. Most of the function classes considered in this paper are convex sets and hence our results do not require the well-specification assumption. This discussion concludes that an analogue of Theorem 5.1 of~\cite{MR2829871} holds even when the response has finite number of moments. In our case, the tail probabilities will decay polynomially and not exponentially as in~\cite{MR2829871}.
% Note that finite number of moments of response or error is equivalent under uniform boundedness of function spaces. 
The analysis, however, is different if $\mathcal{F}$ is a non-convex set. Examples of non-convex function spaces include multivariate quasiconvex functions, single or multiple index models, and sparse linear or non-linear models; e.g., see Section~\ref{sub:multi_index}. When $\F$ is a non-convex set,~\eqref{eq:Quadratic-lower-bound} (or its inequality version) may not hold and finding a rate upper  bound for $\|\widehat{f}-\bar{f}\|$ requires different proof techniques. However, even in this case, the tools developed in this paper can be used, because the proofs for rate bounds under misspecification hinge on the control of an empirical process analogous to~\eqref{eq:Characterization}; see~\cite[Eqns. (1) and (2)]{MR2490397} and~\cite[Theorem 5.2]{MR2829871}. {Following the proof of Theorem 2 of~\cite{massart2006risk}, one can prove the generalizations of Theorems~\ref{thm:SecondMomentRate} and~\ref{thm:L_inf_smooth} for heavy-tailed error $\epsilon$. The choices of $\ell(\cdot, \cdot)$ and $w(\cdot)$ are given in Section 2.3 of~\cite{massart2006risk}. The major change comes in using Proposition 3.1 of~\cite{Gine00} instead of Bousquet's inequality on page 2348 of~\cite{massart2006risk}. We leave the details to the reader.}   

 \section{Concluding remarks} % (fold)
 \label{sec:conclusion}
 Least squares estimators in nonparametric regression models are known to be minimax rate optimal when $\epsilon$ is sub-Gaussian and when $\F$ satisfies appropriate entropy assumptions.  We show that in a wide variety of cases, the LSE attains the same rate of convergence even when $\epsilon$  is neither sub-Gaussian nor independent of $X$. We find sufficient moment conditions on $\epsilon$ under which the rate of convergence of the LSE  under heavy-tailed errors matches the rate of the LSE under sub-Gaussian errors.
 % , i.e., the LSE is ``robust'' to heavy-tailed errors. 
 Our sufficient conditions depend on the complexity ($\alpha$) and the local structure ($s$) of the function class $\F$.  
 % The results justify the use of LSE even under heteroscedastic and heavy-tailed errors. 
   In this paper, all our results focus on the squared error loss but our results can be easily generalized to other smooth loss functions by modifying the proof of Theorem~\ref{thm:GeneralPeeling}; see~\cite[Theorem~3]{chen1998sieve}.

The necessity of our conditions is under investigation. Interestingly, the local structure of $\F$ seems to play a role in the rate of convergence of the LSE only when the errors have only finitely many (conditional) moments, i.e., $q < \infty$. In particular, if the errors are sub-Gaussian, then the local structure can be completely ignored to derive the $n^{-1/(2+\alpha)}$ rate of convergence. {\clr The dependence of the rate of convergence of the LSE on the envelope growth parameter ($s$) is an open problem which is left for future investigation. In particular, for each function class $\mathcal{F}$ and $f_0\in\mathcal{F}$ with the envelope growth parameter $s\in[0, 1]$, the dependence of $s$ on the rate of convergence of the LSE remains unsolved. Recall that Theorem~\ref{thm:Lowerbnd} only establishes the existence of a function class for which the rate of convergence of the LSE is driven by the envelope growth parameter. In this sense, Theorem~\ref{thm:Lowerbnd} is only a worst case result.}

% When the errors are sub-Gaussian then one doesn't need to consider the local structure ($s$) of $\F$ when deriving the rate of convergence of the LSE. However, it turns out that the careful analysis of the local structure is when the errors have only finitely many (conditional) moments.      

In this paper, we have exclusively focused on conditions under which the LSE is ``rate optimal.'' Even though the LSE can attain the global minimax rate of convergence under these conditions, it can be lacking in other aspects such as the tail behavior~\cite{brownlees2015empirical} and accuracy/confidence trade-off~\cite{lugosi2016risk}. Also, one might want a single procedure for all $q\ge2$, instead of changing the procedure depending on some conditions. Several authors~\cite{lugosi2016risk,brownlees2015empirical,lecue2020robust} have taken such concerns into consideration and developed alternative estimators such as median-of-means, Catoni's loss estimators, etc. However, these estimators are still not satisfactory because they do not apply to some of the function classes we can accommodate. 
{Furthermore, the LSE is often favorable because of its intuitive nature, and adaptivity properties for shape-constrained classes.}

% Furthermore, the LSE is not only favorable because of its intuitive nature but also because of its adaptivity properties for shape-constrained classes (which are important in practice). 

% Indeed in certain scenarios, estimators such as  Median-of-means, PAC Bayesian, or $\rho$-estimators are advantageous when the errors are heavy-tailed. However, in settings considered here (shape constrained or smoothness classes), LSE's are often used because of other nice properties, e.g., adaptivity, ease of computation of the estimator, interpretability, among others. 
\bigskip
\begin{appendix}
\begin{center}
{\large\bf APPENDIX}
\end{center}
\section{Interpolation inequalities} % (fold)
\label{sec:auxiliary_results}
In this section, we state  three interpolation inequalities that find the local envelope and the envelope growth parameter for the examples considered in the paper.  The proofs are in Section~\ref{sec:proof_of_propositions_in_appendix_} of the supplementary file.
% section auxiliary_results (end)
\begin{prop}[Local envelope for bounded convex function]\label{lem:Convex_F_delta}
Let 
\begin{equation}\label{eq:ConvexClass}
\F := \left\{f:[0,1]\to[-\Phi,\Phi] \; | \;f \text{ is convex}\right\}
\end{equation}
and $P_X$ be the uniform distribution on $[0,1].$ Fix any $f_0\in \F$, then 
% \[ \mathcal{F}_{\delta} = \{f \in \F:\,\int_0^1 (f - f_0)^2(x)dx \le \delta^2\}\qquad \text{ and } \qquad 
% F_{\delta}(x) = \sup_{f\in\mathcal{F}_\delta}\,|(f - f_0)(x)|.\]
for any $x\in [0,1]$, $F_{\delta}(x) \le  \min\{2(2\Phi)^{1/3}\delta^{2/3}\max\{x^{-1/3}, (1 - x)^{-1/3}\}, 2\Phi\}.$ 
Thus $\|F_\delta\|_{\infty} =2\Phi$ and $\|F_{\delta}\|_3 \le 4 \Phi^{1/3}\delta^{2/3} [\log(0.5\Phi^2/\delta^2)]^{1/3}.$
\end{prop}
% For any function $f:[0,1]^d\to\mathbb{R}$, define 
% \[
% \|f\|_2 = \left(\int_{[0,1]^d} f^2(x)dx\right)^{1/2}.
% \]
\begin{prop}[Local envelope for additive Models]\label{prop:Interpolation_for_Additive_Models}
Suppose $f:[0,1]^d \to \mathbb{R}$ can be written as $f(x) = f_1(x_1) + f_2(x_2)$ for some functions $f_j:[0,1]^{d_j}\to\mathbb{R}, j = 1,2$, for every $x = (x_1^{\top}, x_2^{\top})^{\top}\in[0, 1]^d$ with $x_1\in[0,1]^{d_1}, x_2\in[0,1]^{d_2}$. If $P_X$ is the uniform distribution on $[0,1]^d$ and $f_j\in\mathcal{F}_{\gamma_j,d_j}(L),$ $j = 1, 2$ where for any $\gamma>0$ and dimension $d$ 
\begin{equation}\label{eq:F_gama_d_def}
\mathcal{F}_{\gamma,d}(L) := \left\{f :\rchi\to\mathbb{R}: f/L \in \F_{\gamma, d}\right\}.
\end{equation} 
 Then $\|f\|_{\infty} \le 5\left(\|f\|_2^{c_1} + \|f\|_2^{c_2}\right)\left(L^{1-c_1} + L^{1-c_2}\right),$
where $c_j := 2\gamma_j/(2\gamma_j + d_j),\, j = 1, 2.$ In particular, if $f(x) = \sum_{j=1}^d f_j(x_j)$ for functions $f_j:[0, 1] \to \mathbb{R}$ such that $f_j\in\mathcal{F}_{\gamma,1}(L)$, then $\|f\|_{\infty} \le 5d\|f\|_2^{c}L^{1-c},$ where $ c = 2\gamma/(2\gamma + 1).
$
\end{prop}
\begin{prop}[Local envelope for multiple index models]\label{prop:Interpolation_for_multiple_index_models}
Suppose $f(x) = m(Bx) - m_0(B_0x)$ for functions $m,m_0:\mathbb{R}^p\to\mathbb{R}$ satisfying $m,m_0\in\mathcal{F}_{\gamma,p}(L)$ (defined in~\eqref{eq:F_gama_d_def}) and $B, B_0\in\mathbb{R}^{p\times d}$ with $p < d$. If $X\in\mathbb{R}^d$ is a random vector such that $((BX)^{\top}, (B_0X)^{\top})$ has a density with respect to the Lebesgue measure that is lower bounded by $\underline{C} > 0$, then
\[
\|m\circ B - m_0\circ B_0\|_{\infty} \le 10 \underline{C}^{-c/2}\|m\circ B - m_0\circ B_0\|^{c}L^{1-c},
\]
where $\|m\circ B - m_0\circ B_0\| := (\mathbb{E}[|m(BX) - m_0(B_0X)|^2])^{1/2}$ and $c = 2\gamma/(2\gamma + p)$.
\end{prop}

\section{A new maximal inequality for finite maximums} % (fold)
\label{sec:a_new_maximal_inequality_for_finite_maximums}
The following maximal inequality (proved in Section~\ref{sub:proof_of_proposition_lem:bernouli_maximal} in the supplement) will be used in the proof of Theorem~\ref{thm:L_inf_smooth} but is also  of independent interest. 
{\clr Proposition~\ref{lem:bernouli_maximal} is an analogue of Nemirovski's inequality; see~\cite{dumbgen2010nemirovski},~\cite{massart2013around}, and~\cite[Chapter 11.2 and 11.3]{boucheron2013concentration}. }
% The proof of Proposition~\ref{lem:bernouli_maximal} is in Section~\ref{sub:proof_of_proposition_lem:bernouli_maximal} of the supplementary file. 
\begin{prop}\label{lem:bernouli_maximal}
Let $X_1, \ldots, X_n$ be mean zero independent random variables in $\R^p,\, p\ge 1$. Suppose for all $i \in \{1, \ldots, n\}$
\begin{equation}\label{eq:q-moment}
 \mathbb{E}[\xi_i^q] < \infty \quad \text{where}\quad \xi_i := \max_{1\le j\le p}|X_{i,j}|
 \end{equation}
 and  $X_i := (X_{i,1}, \ldots, X_{i,p})$. If $V_{n,p} := \max_{1\le j\le p}\sum_{i=1}^n \mathbb{E}[X_{i,j}^2]$, then
\begin{equation}\label{eq:maximal_bernouli}
\mathbb{E}\left[\max_{1\le j\le p}\left|\sum_{i=1}^n X_{i,j}\right|\right] \le \sqrt{6V_{n,p}\log(1 + p)} + \sqrt{2}(3\log(1 + p))^{1 - 1/q}\left(2\sum_{i=1}^n \mathbb{E}[\xi_i^q]\right)^{1/q}.
\end{equation}
\end{prop}
\subsection{Example 4: An application of Proposition~\ref{lem:bernouli_maximal}} % (fold)
\label{rem:B1Application}
Suppose we have $n$ i.i.d.~pairs $(X_i, \epsilon_i), 1\le i\le n$ such that $\mathbb{E}[\epsilon_i|X_i] = 0$ and $\mathbb{E}[\epsilon_i^2|X_i] \le \sigma^2$ a.e. $P_X$. Let $\{f_1, \ldots, f_N\}$ be a collection of functions from $\rchi$ to $\mathbb{R}$ and let $F(\cdot) := \max_{1\le j\le N}|f_j(\cdot)|$ denote their envelope. Then Proposition~\ref{lem:bernouli_maximal} (with $q = 2$) yields
\begin{align*}
\mathbb{E}\left[\max_{1 \le j\le N}|\mathbb{G}_n[\epsilon f_j(X)]|\right] &\lesssim \sqrt{{\log N}}\max_{1\le j\le N}\mathbb{E}^{1/2}\left[\frac{1}{n}\sum_{i=1}^n \epsilon_i^2f_j^2(X_{i})\right] \\
&\qquad+ \frac{(\log N)^{1-1/q}}{n^{1/2 - 1/q}}\left(\frac{1}{n}\sum_{i=1}^n \mathbb{E}[\epsilon_i^qF^q(X_i)]\right)^{1/q}\\
&\lesssim \sigma\sqrt{{\log N}}\max_{1\le j\le N}\|f_j\| + \sigma \|F\|\left({\log N}\right)^{1/2}\quad\mbox{(with $q = 2$)}\\
&\lesssim \sigma\sqrt{{\log N}}\left[\max_{1\le j\le N}\|f_j\| + \|F\|\right],
\end{align*}
which implies that 
\begin{equation}\label{eq:applictaionb1}
 \max_{1\le j\le N}|\mathbb{G}_n[\epsilon f_j(X)]| = O_p(\sqrt{\log N})
 \end{equation} whenever $\|F\| = O(1)$. In contrast, Lemma 8 of~\cite{MR3350040} implies
\begin{align*}
\mathbb{E}\left[\max_{1\le j\le N}|\mathbb{G}_n[\epsilon f_j(X)]|\right] &\lesssim \sqrt{{\log N}}\max_j\mathbb{E}^{1/2}\left[\frac{1}{n}\sum_{i=1}^n \epsilon_i^2f^2_j(X_{i})\right] \\
&\qquad + \frac{\log N}{\sqrt{n}}\sqrt{\mathbb{E}\left[\max_{1\le i\le n}|\epsilon_i|^2F^2(X_i)\right]}.
\end{align*}
Even when $\|F\|_{\infty} \le C < \infty$, under only the second moment assumption, $\mathbb{E}[\max_{1\le i\le n}|\epsilon_i|^2] = O(n)$ and hence the second term on the right hand side will be of the order $\log N$. Thus Lemma 8 of~\cite{MR3350040} will imply that $\max_{1\le j\le N}|\mathbb{G}_n[\epsilon f_j(X)]| = O_p(\log N)$. Thus~\eqref{eq:applictaionb1} is a significant improvement, as the above calculation now implies that the lasso estimator is minimax rate optimal under just the conditional second moment assumption~\eqref{eq:cvar}, when the covariates are coordinate-wise bounded; see Theorem 11.1 of~\cite{hastie2015statistical}. Proposition~\ref{lem:bernouli_maximal} can also be used in proving consistency of the multiplier bootstrap under finite moment assumptions~\cite[Remark 5.2]{kuchibhotla2018moving}.

% \end{proof}
% section a_new_maximal_inequality_for_finite_maximums (end)

\section{A refined peeling result} % (fold)
\label{sec:PeelingResult}
In this section, we state a new peeling result. The proof is in Section~\ref{sec:proof_of_theorem_thm:generalpeeling} of the supplementary file. The result is a key component in the proofs of the rate results (Theorems~\ref{thm:SecondMomentRate},~\ref{thm:L_inf_smooth}, and~\ref{thm:Vc-type}) in the paper.  It is this refinement that helps us prove fast rates of convergence of the LSE in previously inaccessible cases. Before stating the result, we will introduce some notations.
Let \begin{equation}\label{eq:hatf}
\widehat{f} := \argmin_{f\in \F} \mathbb{M}_n(f) \qquad \text{where} \qquad \mathbb{M}_n(f) := \frac{1}{n} \sum_{i=1}^{n} (Y_i- f(X_i))^2
\end{equation}
and $T(f; \epsilon, X) := 2\epsilon(f - f_0)(X) - (f - f_0)^2(X).$ Furthermore, let  $U: \R\times\rchi\times \R \to \R$ be such that 
\begin{equation}\label{eq:U_def}
 \sup_{f: \|f-f_0\| \le \delta}|T(f; \epsilon, X)|\le U(\epsilon, X; \delta),
\end{equation}
for all values of $\epsilon$, $X$, and $\delta$. If $\Phi:= \sup_{f\in\F} \|f\|_{\infty}$, then a trivial choice is $U(\epsilon, X; \delta) = 4 (|\epsilon| + \Phi) \Phi$ (we use this choice in the proof of Theorem~\ref{thm:SecondMomentRate}). Now for  any $B > 0$, let 
\begin{equation}\label{eq:TB_def}
T_B(f; \epsilon, X, \delta) := T(f; \epsilon, X)\mathbf{1}\{U(\epsilon, X; \delta) \le B\}.
\end{equation}
Theorem~\ref{thm:GeneralPeeling} below is useful because it provides tail bounds for $\|\widehat{f}-f_0\|$ in terms of upper bounds on a bounded (note that $T(\cdot;\cdot,\cdot,\cdot)$ is unbounded while $T_B(\cdot;\cdot,\cdot,\cdot)$ is bounded) empirical process and most existing maximal inequalities provide upper bounds for only bounded empirical processes.

% Note that this function is different from the one in the previous section.
% section a_more_refined_set_of_maximal_inequalities (end)
\begin{thm}[Peeling with truncation] \label{thm:GeneralPeeling} Suppose $f_0\in \F$,  $\Phi:= \sup_{f\in\F} \|f\|_{\infty}$,  and there exists a real valued function $\phi_n(\cdot; \cdot)$ such that
\begin{equation}\label{eq:Max_ineq_requirement11}
\sup_{\{f: \|f-f_0\|\le \delta\}} \E T^2_B(f; \varepsilon, X, \delta) \le  4(\sigma+\Phi)^2 \delta^2 
\end{equation}
and
\begin{equation}\label{eq:Max_ineq_requirement}
\mathbb{E}\left[\sup_{\delta/2 \le \|f - f_0\| \le \delta}\mathbb{G}_n(T_B(f; \epsilon, X, \delta))\right] \le \phi_n(\delta; B),
\end{equation}
for  every $n$ and any $\delta, B > 0$.
Further, if there exists $\gamma \ge 2$ and $s_{\gamma}(\cdot)$ such that $\mathbb{E}\left[U^\gamma(\epsilon, X; \delta)\right] \le s_{\gamma}(\delta),$ then there exists a universal constant $C$ such that for any positive $D, \varepsilon_n,$  and $\{B_k\}_{k = 1}^{\infty}$,  $\tau \ge 1$, and $n\ge1$ we have
\begin{align}
\mathbb{P}\left(\|\widehat{f} - f_0\| \ge D\varepsilon_n\right) &\le 
\left(\frac{C}{\sqrt{n}(D\varepsilon_n)^2}\right)^\tau\sum_{k = 1}^{\infty}\frac{\phi_n^\tau(2^kD\varepsilon_n, B_k)}{2^{2k\tau}} + \left(\frac{ 8 C\sqrt{\tau}(\sigma + \Phi)}{\sqrt{n}(D\varepsilon_n)}\right)^\tau\nonumber\\
 &\qquad + \left(\frac{C\tau}{n(D\varepsilon_n)^2}\right)^\tau\sum_{k = 1}^{\infty} \frac{B_k^\tau}{2^{2k\tau}} + \frac{16}{(D\varepsilon_n)^2}\sum_{k = 0}^{\infty} \frac{s_{\gamma}(2^kD\varepsilon_n)}{2^{2k}B_k^{\gamma - 1}}\label{eq:Error_of_fhat}
\end{align}
\end{thm}

\begin{remark}\label{rem:Dydic Peeling}
\cite[Theorem 3]{shen1994convergence} propose an iterative (non-dyadic) version of the dyadic peeling argument presented here. After a thorough investigation, we have found that their approach does not lead to better rates. It should be also noted that a version of Theorem~\ref{thm:GeneralPeeling} for general loss function can be derived using Theorem~3 of~\cite{shen1994convergence}. We refrain from this to keep the paper focused on the LSE.
\end{remark}

\begin{remark}\label{rem:Unbounded_empirical_process}
The motivation behind the truncation argument in Theorem~\ref{thm:GeneralPeeling} was to use existing maximal inequalities  to  get $\phi_n$ in~\eqref{eq:Max_ineq_requirement}. There are however maximal inequalities that allow for unbounded stochastic processes; see e.g.,~\cite[Lemma 6.12]{MR1915446},~\cite{van2011local}, and~\cite[Theorem 1.9]{mendelson2016upper}. We did not use them in this paper, because~\cite[Lemma 6.12]{MR1915446} and~\cite[Theorem~3.1]{van2011local} lead to slower rates in Theorems~\ref{thm:SecondMomentRate} and~\ref{thm:Vc-type} and \cite[Theorem 1.9]{mendelson2016upper} is not applicable for the examples in this paper. 
\end{remark}

% \begin{remark}\label{rem:Hoffmann-J}
% Is Hoffman-J{\o}rgensen the right truncation process? Apparently not, the procedure of~\cite{shen1994convergence} is better. More...
% \end{remark}
% \todo[inline]{How is this different from the existing peeling?}
% \todo[inline]{Change $q$. It is not the same as the number of moments of $\epsilon$. }
\noeqref{eq:Error_of_fhat}
\noeqref{eq:U_def}
\noeqref{eq:Max_ineq_requirement}

\end{appendix}

\bibliographystyle{apalike}
\bibliography{SigNoise}

\newpage
\setcounter{section}{0}
\setcounter{equation}{0}
\setcounter{figure}{0}

% \cftsetindents{section}{1em}{2.5em}
% \cftsetindents{subsection}{1.5em}{3em}

\renewcommand{\thesection}{S.\arabic{section}}
\renewcommand{\theequation}{S.\arabic{equation}}
\begin{center}
  \Large {\bf Supplement to ``On Least Squares Estimation Under Heteroscedastic and Heavy-Tailed Errors''}
  \end{center}
% {
%   \hypersetup{linkcolor=blue}
%   \input{2ndRevSuppDraft0.toc}
% }

\section{Discussion on the local envelope}\label{sec:discussion_on_the_local_envelope}
{\clr In this section, we will provide a heuristic argument that suggests that the local envelope is an important quantity to consider in the study of the rate of convergence of the LSE under heavy-tailed noise.

Assuming equivalence of the left hand side of~\eqref{eq:Characterization} and~\eqref{eq:NecessaryExpectation}, the rate of convergence of the LSE $\widehat{f}$ defined on the function space $\mathcal{F}$ is characterized by
\[
\delta~~\mapsto~~\overline{\phi}_n(\delta) ~:=~ \mathbb{E}\left[\sup_{f\in\mathcal{F}:\,\|f - f_0\| \le \delta}\left|\mathbb{G}_n[\epsilon(f - f_0)(X)]\right|\right].
\]
Theorem 1.4.4 (and Remark 1.4.6) of~\cite{MR1666908}, with $p = q = 1$ implies that
\begin{equation}\label{eq:main-conclusion}
\begin{split}
\overline{\phi}_n(\delta) ~&\asymp~ \mathbb{E}\left[\sup_{f\in\mathcal{F}:\|f - f_0\| \le \delta}\left|\mathbb{G}_n[S(f, \epsilon)\mathbbm{1}\{\left|\epsilon\right|F_{\delta}(X) \le t_0\}]\right|\right]\\ 
~&\quad+~ \frac{1}{\sqrt{n}}\mathbb{E}\left[\max_{1\le i\le n}\left|\epsilon_i\right|F_{\delta}(X_i)\right],
\end{split}
\end{equation}
where $S(f, \epsilon) ~:=~ \epsilon(f - f_0)(x),$ for $F_{\delta}(x) := \sup_{f\in\mathcal{F}_{\delta}}|(f - f_0)|(x)$
and
\[
t_0 ~:=~ \inf\left\{t:\,\sum_{i=1}^n \mathbb{P}\left(\left|\epsilon_i\right|F_{\delta}(X_i) \ge t\right) \le \frac{1}{8}\right\}.
\]
If $\epsilon F_\delta$ satisfies
\[ \mathbb{P}\left(\frac{|\epsilon F_\delta(X)|}{C_\delta} \ge t\right) = \frac{\log^2 2}{t^q \log^2 (1+t)} \text{ for } t\ge 1,\]
 then $\|\epsilon F_\delta(X)\|_q\asymp C_\delta$ and~\cite[Proposition 1.4.1]{MR1666908} implies
\begin{align*}
n^{1/q} \frac{\|\epsilon F_\delta(X)\|_q}{\log^{2/q} (n/\|\epsilon F_\delta(X)\|_q)} ~\lesssim~ \mathbb{E}\left[\max_{1\le i\le n}\left|\epsilon_i\right|F_{\delta}(X_i)\right] ~\lesssim~ n^{1/q} {\|\epsilon F_\delta(X)\|_q},
\end{align*}
where $ \|\epsilon F_\delta(X)\|_q= \left(\mathbb{E}\left[\left|\epsilon_i\right|^qF_{\delta}^q(X_i)\right]\right)^{1/q}.$
Further
\[ t_0\asymp n^{1/q} \frac{\|\epsilon F_\delta(X)\|_q}{\log^{2/q} (n/\|\epsilon F_\delta(X)\|_q)}. \]

Both of the terms on the right of~\eqref{eq:main-conclusion} depend on the $ \|\epsilon F_\delta(X)\|_q$; the first summand of~\eqref{eq:main-conclusion} via the truncation $t_0$ and the second summand directly. 
 There also exist cases where the second summand in the right hand side~\eqref{eq:main-conclusion} is of the same order as the left hand side of~\eqref{eq:main-conclusion}; see \eqref{eq:tightlwrbound} in the proof of Theorem~\ref{thm:Lowerbnd} for an example.
% Take $\mathcal{F} = \{x\mapsto\langle x, \beta\rangle:\,\|\beta\|_2 \le 1\}$, $X_i, 1\le i\le n$ are isotropic random vectors in $\mathbb{R}^d$ satisfying, for some $K\ge1$,
% \[
% \mathbb{P}\left(\|X_i\|_2/K \ge t\right) \le \frac{1}{1 + t^p}\quad\mbox{for} t\ge0.
% \]
% Finally, take $\epsilon$ to be a Rademacher variable independent of $X$.\footnote{We thank the reviewer for providing this example.} In this case,
% \[
% \mathcal{F}_{\delta} = \{x\mapto\langle x, \beta\rangle:\,\|\beta - \beta_0\|_2 \le \delta\}\quad\mbox{and}\quad F_{\delta}(x) = \sup_{\beta:\,\|\beta - \beta_0\|_2 \le \delta}|\langle X, \beta - \beta_0\rangle| = \delta\|X\|_2.
% \]
% In this case, the left hand side of~\eqref{eq:main-conclusion} is
% \[
% \mathbb{E}\left[\sup_{\beta:\,\|\beta - \beta_0\|_2 \le \delta}\left|\frac{1}{\sqrt{n}}\sum_{i=1}^n \epsilon_i\langle X_i, \beta - \beta_0\rangle\right|\right] = {\delta}\mathbb{E}\left\|\frac{1}{\sqrt{n}}\sum_{i=1}^n \epsilon_iX_i\right\|_2 \le \delta\sqrt{\mathbb{E}[\|X_i\|_2^2]}.
% \]
% The second summand on the right hand side of~\eqref{eq:main-conclusion} is given by
% \[
% \frac{1}{\sqrt{n}}\mathbb{E}\left[\max_{1\le i\le n}|\epsilon_i|F_{\delta}(X_i)\right] = \frac{\delta}{\sqrt{n}}\mathbb{E}\left[\max_{1\le i\le n}\right]
% \]

Thus, in general, the lower bound on $\overline{\phi}_n(\delta)$ depends on $ \|\epsilon F_\delta(X)\|_q$. Assuming

\[ \overline{\phi}_n(\delta) ~\gg~ \E\left[\sup_{f\in \F: \|f-f_0\|\le \delta} \mathbb{G}_n(f-f_0)^2(X)\right],\]
the lower bound on \[\mathbb{E}\left[\sup_{f\in\mathcal{F}:\,\|f - f_0\| \le \delta}\big|\mathbb{G}_n[2\epsilon(f - f_0)(X) - (f - f_0)^2(X)]\big|\right]\] will depend on the ``size'' of the local envelope $F_{\delta}$. Based on this lower bound, one can find a lower bound on $\|\widehat{f}_n- f_0\|$ using Proposition~6 of~\cite{han2017sharp}. 
% Hence, one can not claim that rate of the LSE will never depend on the local envelope $F_\delta.$ 
The argument above is only a heuristic and formalizing this is beyond the scope of the current paper.  }
\section{Proof of propositions in Appendix~\ref{sec:auxiliary_results}} % (fold)
\label{sec:proof_of_propositions_in_appendix_}
In this section, we prove the propositions in Appendix~\ref{sec:auxiliary_results} of the main paper.
% section auxiliary_results (end)

\subsection{Proof of Proposition~\ref{lem:Convex_F_delta}} % (fold)
\label{sub:proof_of_proposition_lem:convex_f_delta}

% subsection proof_of_proposition_lem:convex_f_delta (end)

A convex function $f$ bounded by $\Phi$ on $[0, 1]$ is Lipschitz on any sub-interval $[a,b]$ with Lipschitz constant $2\Phi/\min\{a, 1-b\}$. Fix any $x\in(0, 1/2]$. On any interval $[a,b]\subseteq[0,1]$ containing $x$, $f$ and $f_0$ are both Lipschitz with Lipschitz constant $2\Phi/a$ (which implies that $f - f_0$ is Lipschitz with Lipschitz constant $4\Phi/a$). Using the interpolation for Lipschitz functions~\cite[Lemma 2]{chen1998sieve}, we have
\begin{equation}\label{eq:interpoconvex}
|(f - f_0)(x)| \le 2\left(\int_{a}^{b} |f - f_0|^2(t)dt\right)^{1/3}\left(\frac{2\Phi}{a}\right)^{1/3}.
\end{equation}
Hence if $\mathcal{F}_{\delta} := \{f:[0,1]\to[-\Phi,\Phi]\, |\,\|f - f_0\| \le \delta\text{ and } f\in \F\}$, then for every $0< x\le 1/2,$ we have
\[
\sup_{f\in\mathcal{F}_{\delta}}|(f - f_0)(x)| \le 2(2\Phi)^{1/3}\delta^{2/3}x^{-1/3},
\]
where we replaced $a$ in~\eqref{eq:interpoconvex} by $x$ by taking limit $a\downarrow x$. 
Thus, by symmetry
\[
F_{\delta}(x) \le 2(2\Phi)^{1/3}\delta^{2/3}\max\{x^{-1/3}, (1 - x)^{-1/3}\}.
\]
However, $\|f\|_\infty \le \Phi$ for every $f\in \F$ thus 
\[
F_{\delta}(x) \le \min \left\{2(2\Phi)^{1/3}\delta^{2/3}\max\{x^{-1/3}, (1 - x)^{-1/3}\}, 2\Phi\right\}
\]
To find the upper bound on $\|F_\delta\|_3$, observe that  
\[
\int_0^{1/2} |F_{\delta}^3(x)|dx =\int_0^{\eta} |F_{\delta}^3(x)|dx + \int_{\eta}^{1/2} |F_{\delta}^3(x)|dx \le  8\Phi^3\eta + (2(2\Phi)^{1/3}\delta^{2/3})^{3}\log(1/(2\eta)).
\]
Taking $\eta = \Phi^{-2}\delta^2$ implies
\[
\|F_{\delta}\|_3^3 \le 16(2\Phi)\delta^2\log(0.5\Phi^2/\delta^2). \qedhere
\]
% For any function $f:[0,1]^d\to\mathbb{R}$, define 
% \[
% \|f\|_2 = \left(\int_{[0,1]^d} f^2(x)dx\right)^{1/2}.
% \]
\subsection{Proof of Proposition~\ref{prop:Interpolation_for_Additive_Models}} % (fold)
\label{sub:proof_of_proposition_prop:interpolation_for_additive_models}

% subsection proof_of_proposition_prop:interpolation_for_additive_models (end)
Consider $f(x_1, x_2) = f_1(x_1) + f_2(x_2)$ for two functions $f_1:\mathbb{R}^{d_1}\to\mathbb{R}, f_2:\mathbb{R}^{d_2}\to\mathbb{R}$. Define
\[
\bar{f}_1(x_1) = f_1(x_1) - \int f_1(t)dt\quad\mbox{and}\quad \bar{f}_2(x_2) = f_2(x_2) - \int f_2(t)dt.
\]
Because $f_j\in\mathcal{F}_{\gamma_j,d_j}(L)$ are $\gamma_j$-smooth, $\bar{f}_j, j=1,2$ are also $\gamma_j$-smooth. Hence by~\cite[Lemma 2]{chen1998sieve}, we have 
\[
\|\bar{f}_1\|_{\infty} \le 2\|\bar{f}_1\|_2^{c_1}L^{1-c_1}\quad\mbox{and}\quad \|\bar{f}_2\|_{\infty} \le 2\|\bar{f}_2\|_2^{c_1}L^{1-c_1},
\]
where $c_j = 2\gamma_j/[2\gamma_j + d_j]$, $j = 1,2$. Observe now that
\[
\mathbb{E}[(\bar{f}_1(X_1) + \bar{f}_2(X_2))^2] = \|\bar{f}_1\|_2^2 + \|\bar{f}_2\|_2^{2}.
\]
Therefore,
\begin{align*}
\sup_{(x_1,x_2)\in[0,1]^{d_1+d_2}}\,|\bar{f}_1(x_1) + \bar{f}_2(x_2)| &\le  \|\bar{f}_1\|_{\infty} + \|\bar{f}_2\|_{\infty}\\ 
&\le 2\left[(\|\bar{f}_1\|_2 + \|\bar{f}_2\|_2)^{c_1} + (\|\bar{f}_1\|_2 + \|\bar{f}_2\|_2)^{c_2}\right](L^{1-c_1} + L^{1-c_2}).
\end{align*}
Because 
\begin{equation}\label{eq:teoi9}
\|f\|_1 \le  \|f_1\|_1 + \|f_2\|_1 \text{ and }\|f\|^2 = \|\bar{f}_1\|^2+ \|\bar{f}_2\|^2 +\|\bar{f}_1\|_1^2+ \|\bar{f}_2\|_1^2,
\end{equation} we get
\begin{align*}
\|f\|_{\infty} &\le \int |f(x_1, x_2)|dx_1dx_2 + 4 \left[\|f\|_2^{c_1} + \|f\|_2^{c_2}\right](L^{1-c_1} + L^{1-c_2})\\
&\le \|f\|_2 + 4[\|f\|_2^{c_1} + \|f\|_2^{c_2}](L^{1-c_1} + L^{1-c_2}).
\end{align*}
Furthermore, since $\|f\|_{\infty} \le L$, we get $\|f\|_2^{1-c_1} \le L^{1-c_1}$ and $\|f\|_2^{1-c_2} \le L^{1-c_2}$. Thus, we have
\[
\|f\|_{\infty} \le 5(\|f\|_2^{c_1} + \|f\|_2^{c_2})(L^{1-c_1} + L^{1-c_2}).
\]
If $f(x_1,\ldots,x_d) = \sum_{j=1}^d f_j(x_j)$ and $f_j(\cdot)$ is $\gamma$-smooth for all $1\le j\le d$, then
\[
\|f\|_{\infty} \le 5d\|f\|_{2}^{c}L^{1-c}, \qquad \text{where}\quad c = 2\gamma/[2\gamma + 1].\qedhere
\]

\subsection{Proof of Proposition~\ref{prop:Interpolation_for_multiple_index_models}} % (fold)
\label{sub:proof_of_proposition_prop:interpolation_for_multiple_index_models}

Define
\[
y := \begin{pmatrix}y_1\\y_2\end{pmatrix} := \begin{pmatrix}Bx\\B_0x\end{pmatrix}\in\mathbb{R}^{2p}.
\]
Then we can write $f(y) = f_1(y_1) + f_2(y_2)$ where $f_1(y_1) = m(y_1)$ and $f_2(y_2) = -m_0(y_2)$. Observe that if $Y = ((BX)^{\top}, (B_0X)^{\top})^{\top}$ has a density with respect to the Lebesgue measure bounded away from zero, that is, $p_Y(y) \ge \underline{C}$ with $p_Y(\cdot)$ representing the pdf of $Y$, then
\begin{equation}\label{eq:L2-Norms-bounding}
\|f\|_{2,Y} = (\mathbb{E}[f^2(Y)])^{1/2} \ge \underline{C}^{1/2}\left(\int f^2(y)dy\right)^{1/2}.
\end{equation}
Applying Proposition~\ref{prop:Interpolation_for_Additive_Models} with $\gamma_1 = \gamma_2 = \gamma$ and $d_1 = d_2 = p$ we get
\[
\|f\|_{\infty} \le 10\Big[\int f^2(y)dy\Big]^{c/2}L^{1-c},
\]
where $c = 2\gamma/(2\gamma+ p)$. 
% (The factor $4$ instead of $8$ follows from the proof of Proposition~\ref{prop:Interpolation_for_Additive_Models}.)
Hence from~\eqref{eq:L2-Norms-bounding}, we get
\[
\|f\|_{\infty} \le  10\underline{C}^{-c/2}\|f\|_{2,Y}^cL^{1-c},
\] 
which implies the result.

% subsubsection proof_of_proposition_prop:interpolation_for_multiple_index_models (end)
% section proof_of_propositions_in_appendix_ (end)

\section{Proof of Proposition~\ref{lem:bernouli_maximal}} % (fold)
\label{sub:proof_of_proposition_lem:bernouli_maximal}
  For any $B > 0$, define
  \[
  T_B(x) = \begin{cases}-B, &\mbox{if }x \le -B,\\x, &\mbox{if }|x| \le B,\\B, &\mbox{if }x \ge B\end{cases}\quad\mbox{and}\quad T_B^c(x) = x - T_B(x). 
  \]
  It is clear that $x = T_B(x) + T_B^c(x)$ for all $B > 0$. Hence,
  \begin{align}\label{eq:main-decomposition-heavy-tailed}
  \begin{split}
  \mathbb{E}\left[\max_{1\le j\le p}\left|\sum_{i=1}^n X_{i,j}\right|\right] &= \mathbb{E}\left[\max_{1\le j\le p}\left|\sum_{i=1}^n \{T_B(X_{i,j}) - \mathbb{E}[T_B(X_{i,j})]\}\right|\right]\\
  &\qquad+ \mathbb{E}\left[\max_{1\le j\le p}\left|\sum_{i=1}^n \{T_B^c(X_{i,j}) - \mathbb{E}[T_B^c(X_{i,j})]\}\right|\right]\\
  &=\mathbf{I} + \mathbf{II}.
  \end{split}
  \end{align}
Because $|T_B(X_{i,j})| \le B$, Eq. (2) of~\cite{Geer13} implies
  \[
  \mathbf{I} \le \sqrt{6V_{n,p}\log(1 + p)} + 3B\log(1 + p).
  \]
  Further, 
  \begin{align*}
  \mathbf{II} \le 2\sum_{i=1}^n \mathbb{E}\left[\max_{1\le j\le p}\left|T_B^c(X_{i,j})\right|\right] &\le 2\sum_{i=1}^n \mathbb{E}\left[(\xi_i - B)_+\right]\\
  &= 2\sum_{i=1}^n \int_{0}^{\infty} \mathbb{P}(\xi_i \ge B + t)dt\\
  &= 2\sum_{i=1}^n \int_B^{\infty} \frac{qt^{q-1}}{qt^{q-1}}\mathbb{P}(\xi_i \ge t)dt\\
  &\le \frac{2}{qB^{q-1}}\sum_{i=1}^n \int_B^{\infty} qt^{q-1}\mathbb{P}(\xi_i \ge t)dt \le \frac{2}{qB^{q-1}}\sum_{i=1}^n \mathbb{E}[\xi_i^q].
  \end{align*}
  Combining the bounds on $\mathbf{I}$ and $\mathbf{II}$, we conclude
  \[
  \mathbb{E}\left[\max_{1\le j\le p}\left|\sum_{i=1}^n X_{i,j}\right|\right] \le \sqrt{6V_{n,p}\log(1 + p)} + 3B\log(1 + p) + \frac{2}{qB^{q-1}}\sum_{i=1}^n \mathbb{E}[\xi_i^q].
  \]
  Minimizing over $B > 0$ yields
  \begin{align*}
  \mathbb{E}\left[\max_{1\le j\le p}\left|\sum_{i=1}^n X_{i,j}\right|\right] &\le \sqrt{6V_{n,p}\log(1 + p)} + \frac{q(3\log(1 + p))^{1 - 1/q}}{q-1}\left(\frac{2(q-1)}{q}\sum_{i=1}^n \mathbb{E}[\xi_i^q]\right)^{1/q}\\%\frac{q}{q-1}\\
  &\le \sqrt{6V_{n,p}\log(1 + p)} + \sqrt{2}(3\log(1 + p))^{1 - 1/q}\left(2\sum_{i=1}^n \mathbb{E}[\xi_i^q]\right)^{1/q}. 
  \end{align*}

\section{Proof of Theorem~\ref{thm:GeneralPeeling}} % (fold)
\label{sec:proof_of_theorem_thm:generalpeeling}
The proof follows along the standard peeling argument~\cite[Theorem~3.2.5]{VdVW96}. But the crucial observation here is that the empirical processes involved here are not bounded and thus to be able to apply the rich literature of maximal inequalities we truncate the empirical process involved in the peeling step. The proof is split into 3 main steps. \\

\noindent\textbf{Step 1: Peeling and truncation.} 
From the definition of $\widehat{f}$, it follows that
\begin{align*}
\mathbb{P}\left(\|\widehat{f} - f_0\| \ge D\varepsilon_n\right) &\le \mathbb{P}\left(\sup_{\|\widehat{f} - f_0\| \ge D\varepsilon_n}\,\mathbb{M}_n(f_0) - \mathbb{M}_n(f) \ge 0\right)\\
&\le \mathbb{P}\left(\bigcup_{k = 0}^{\infty}\left\{\sup_{f\in\mathcal{A}_{k}}\mathbb{M}_n(f_0) - \mathbb{M}_n(f) \ge 0\right\}\right),
\end{align*}
where $$\mathcal{A}_{k} = \{f:\,2^{k-1}D\varepsilon_n \le \|f - f_0\| \le 2^kD\varepsilon_n\}.$$ From the definition of $f_0$, we obtain
\begin{equation}\label{eq:lower-bound-main}
\mathbb{E}\left[\mathbb{M}_n(f) - \mathbb{M}_n(f_0)\right] = \mathbb{E}\left[ f^2(X) - 2Y(f(X)-f_0(X)) -f_0^2(X)\right] = \|f - f_0\|^2,
\end{equation}
as $\E(Y|X) = f_0(X)$. Thus
\begin{align}
\mathbb{P}\left(\|\widehat{f} - f_0\| \ge D\varepsilon_n\right) &\le \mathbb{P}\left(\bigcup_{k = 1}^{\infty} \left\{\sup_{f\in\mathcal{A}_{k}}\mathbb{G}_n(T(f;\epsilon, X)) \ge \sqrt{n}(2^{k-1}D\varepsilon_n)^2\right\}\right)\nonumber.
\end{align}
Observe that $T$ is unbounded. To control the tail probabilities, we will truncate $T$ at a sequence $\{B_k\}_{k=1}^{\infty}$. Thus
\begin{align}
\mathbb{P}\left(\|\widehat{f} - f_0\| \ge D\varepsilon_n\right)
 &\le \mathbb{P}\left(\bigcup_{k = 1}^{\infty}\left\{\sup_{f\in\mathcal{A}_k}\mathbb{G}_n(T_{B_k}(f; \epsilon, X,2^kD\varepsilon_n)) \ge \sqrt{n}(2^{k-1}D\varepsilon_n)^2/2\right\}\right)\label{eq:p1p2split}\\
&\quad + \mathbb{P}\left(\bigcup_{k = 1}^{\infty}\left\{\sup_{f\in\mathcal{A}_k}\mathbb{G}_n(T(f; \epsilon, X) - T_{B_k}(f; \epsilon, X,2^kD\varepsilon_n)) \ge \sqrt{n}(2^{k-1}D\varepsilon_n)^2/2\right\}\right)\nonumber\\
&=: \mathbf{P}_1 + \mathbf{P}_2.\nonumber
\end{align}
Observe that $\mathbf{P}_1$ corresponds to the bounded part and $\mathbf{P}_2$ corresponds to the unbounded part.  \\

\noindent\textbf{Step 2: The unbounded part.}
To bound $\mathbf{P}_2$, observe that by definition of $U(\cdot, \cdot; \cdot)$ 
\[
\sup_{f\in\mathcal{A}_{k}}|T(f; \epsilon, X)| \le U(\epsilon, X; 2^{k}D\varepsilon_n).
\]
Therefore, for any $f\in\mathcal{A}_{k}$, we have  
\begin{align*}
|\mathbb{G}_n(T(f; \epsilon, X) - T_{B_k}(f; \epsilon, X))| &\le \sqrt{n}(\mathbb{P}_n + P)|T(f; \epsilon, X)|\mathbf{1}\{U(\epsilon, X; 2^kD\varepsilon_n) \ge B_k\}\\ 
&\le \sqrt{n}(\mathbb{P}_n + P)U(\epsilon, X; 2^kD\varepsilon_n)\mathbf{1}\{U(\epsilon, X; 2^kD\varepsilon_n) \ge B_k\},
\end{align*}
and hence,
\begin{align}\label{eq:P2_bound}
\begin{split}
\mathbf{P}_2 &= \mathbb{P}\left(\bigcup_{k = 1}^{\infty}\left\{\sup_{f\in\mathcal{A}_k}\mathbb{G}_n(T(f; \epsilon, X) - T_{B_k}(f; \epsilon, X,2^kD\varepsilon_n)) \ge \sqrt{n}(2^{k-1}D\varepsilon_n)^2/2\right\}\right)\\
&\le\sum_{k = 1}^{\infty} \mathbb{P}\left(\sup_{f\in\mathcal{A}_k}\mathbb{G}_n(T(f; \epsilon, X) - T_{B_k}(f; \epsilon, X,2^kD\varepsilon_n)) \ge \sqrt{n}(2^{k-1}D\varepsilon_n)^2/2\right)\\
&\le\sum_{k = 1}^{\infty} \frac{1}{ \sqrt{n}(2^{k-1}D\varepsilon_n)^2/2}\mathbb{E}\left(\sup_{f\in\mathcal{A}_k}\mathbb{G}_n(T(f; \epsilon, X) - T_{B_k}(f; \epsilon, X,2^kD\varepsilon_n))\right)\\
&\le\sum_{k = 1}^{\infty} \frac{1}{ \sqrt{n}(2^{k-1}D\varepsilon_n)^2/2}\mathbb{E}\Big(\sqrt{n}(\mathbb{P}_n + P)U(\epsilon, X; 2^kD\varepsilon_n)\mathbf{1}\{U(\epsilon, X; 2^kD\varepsilon_n) \ge B_k\}\Big)\\
&\le \sum_{k = 1}^{\infty} \frac{4\mathbb{E}\left[U(\epsilon, X; 2^kD\varepsilon_n)\mathbf{1}\{U(\epsilon, X; 2^kD\varepsilon_n) \ge B_k\}\right]}{(2^{k-1}D\varepsilon_n)^2}.
\end{split}
\end{align}
Because $\mathbb{E}\left[U^\gamma(\epsilon, X; \delta)\right] \le s_{\gamma}(\delta),$ we have
\begin{equation}
\mathbb{E}[U(\epsilon, X; 2^kD\varepsilon_n)\mathbf{1}\{U(\epsilon, X; 2^kD\varepsilon_n) \ge B_k\}] \le \frac{\mathbb{E}[U^{\gamma}(\epsilon, X; 2^kD\varepsilon_n)]}{B_k^{\gamma - 1}} \le \frac{s_{\gamma}(2^kD\varepsilon_n)}{B_k^{\gamma - 1}}.
\end{equation}
Thus 
\begin{equation}\label{eq:Gamma_U_moment}
\mathbf{P}_2 \le \frac{16}{(D\varepsilon_n)^2}\sum_{k = 1}^{\infty} \frac{s_{\gamma}(2^kD\varepsilon_n)}{2^{2k}B_k^{\gamma - 1}}.
\end{equation}

\noindent\textbf{Step 3: The bounded part.} Now to bound $\mathbf{P}_1$, observe that
\begin{align}\label{eq:f_split1}
\begin{split}
\mathbf{P}_1
&\le \sum_{k = 1}^{\infty} \mathbb{P}\left(\sup_{f\in\mathcal{A}_k}\mathbb{G}_nT_{B_k}(f; \epsilon, X,2^kD\varepsilon_n) \ge \sqrt{n}(2^{k-1}D\varepsilon_n)^2/2\right).
\end{split}
\end{align}
 We will use Markov's inequality to bound each of the term on the right in the above display. For any  $\tau \ge 1$, we have 
\[
\mathbb{P}\left(\sup_{f\in\mathcal{A}_k}\mathbb{G}_nT_{B_k}(f; \epsilon, X,2^kD\varepsilon_n) \ge \frac{\sqrt{n}(2^{k-1}D\varepsilon_n)^2}{2}\right) \le \frac{8^\tau \mathbb{E}\left[\sup_{f\in\mathcal{A}_k}\left|\mathbb{G}_nT_{B_k}(f; \epsilon, X,2^kD\varepsilon_n)\right|^\tau \right]}{\left(\sqrt{n}(2^kD\varepsilon_n)^2\right)^\tau}.
\]
% \todo[inline]{Define $\phi_n$.}
Using Proposition 3.1 of~\cite{Gine00},  for every $\tau\ge 1$,\footnote{We get non trivial bound   for every $\tau>0$ because $T_B(\cdot;\cdot,\cdot,\cdot)$ uniformly bounded.} we get
\begin{align}
\mathbb{E}\left[\sup_{f\in\mathcal{A}_k}\left|\mathbb{G}_nT_{B_k}(f; \epsilon, X,2^kD\varepsilon_n)\right|^\tau\right] &\le C^\tau\phi_n^\tau(2^kD\varepsilon_n, B_k) + C^\tau\tau^{\tau/2}\sup_{f\in\mathcal{A}_k}\left(\mathbb{E}\left[T_{B_k}^2(f; \epsilon, X,2^kD\varepsilon_n)\right]\right)^{\tau/2}\nonumber\\ &\qquad+ C^\tau\frac{\tau^\tau}{n^{\tau/2}}\mathbb{E}\left[\max_{1\le i\le n}\sup_{f\in\mathcal{A}_k}|T_{B_k}(f; \epsilon_i, X_i,2^kD\varepsilon_n)|^\tau\right]\label{eq:GineLatala}\\
&\le C^\tau \left[\phi_n^\tau(2^kD\varepsilon_n, B_k) + \tau^{\tau/2}(\sigma + M)^\tau(2^kD\varepsilon_n)^\tau +  \frac{\tau^\tau B_k^\tau}{n^{\tau/2}}\right]\nonumber,
\end{align}
where $C$ is a universal constant.  This is because by assumption of the theorem, we have 
% \[\sup_{\{f: \|f-f_0\|\le \delta\}} \E T^2_B(f; \varepsilon, X, \delta) \le  4(\sigma+\Phi)^2 \delta^2\]
 { \begin{align*}
   \sup_{f\in\mathcal{A}_k}\mathbb{E}\left[T_{B_k}^2(f; \epsilon, X,2^kD\varepsilon_n)\right] & \le 4 (\sigma + \Phi)^2 (2^kD\varepsilon_n)^2
   \end{align*}
   and $|T_{B_k}(f; \epsilon_i, X_i,2^kD\varepsilon_n)| \le B_k.$ }
% \todo[inline]{Should have a a $2^q$ right? Don't we have $T_B -P T_B$?}
Thus the above three displays combined with~\eqref{eq:Gamma_U_moment}, we get
\begin{align*}
\mathbb{P}\left(\|\widehat{f} - f_0\| \ge D\varepsilon_n\right) &\le \left(\frac{C}{\sqrt{n}(D\varepsilon_n)^2}\right)^\tau\sum_{k = 1}^{\infty}\frac{\phi_n^\tau(2^kD\varepsilon_n, B_k)}{2^{2k\tau}} + \left(\frac{ 8 C\sqrt{\tau}(\sigma + \Phi)}{\sqrt{n}(D\varepsilon_n)}\right)^\tau \sum_{k = 1}^{\infty}\frac{1}{2^{k\tau}}\\ &\qquad + \left(\frac{C\tau}{n(D\varepsilon_n)^2}\right)^\tau\sum_{k = 1}^{\infty} \frac{B_k^\tau}{2^{2k\tau}} + \frac{16}{(D\varepsilon_n)^2}\sum_{k = 1}^{\infty} \frac{s_{\gamma}(2^kD\varepsilon_n)}{2^{2k}B_k^{\gamma - 1}}.\qedhere
\end{align*}

% section proof_of_theorem_thm:generalpeeling (end)

\section{Proof of Theorem~\ref{thm:SecondMomentRate}} % (fold)
\label{sec:proofSecondMomentRate}

The proof proof of Theorem~\ref{thm:SecondMomentRate} will be split into two steps. In the first step, we find $\phi_n(\delta,B)$ satisfying the assumptions of Theorem~\ref{thm:GeneralPeeling}. In the second step we apply Theorem~\ref{thm:GeneralPeeling} with an appropriate choice of $\beta, U$, and  $s_{\gamma}$.\\

\noindent\textbf{Finding $\phi_n(\delta,B)$:} We will find a choice for $\phi_n(\delta, B)$ that satisfies~\eqref{eq:Max_ineq_requirement}. To find $\phi_n(\cdot, \cdot)$, we will apply Lemma 3.4.2 of~\cite{VdVW96}. Recall that 
\begin{equation}
T(f; \epsilon, X) = 2\epsilon(f - f_0)(X) - (f - f_0)^2(X),
\end{equation}
Observe from the definition, we can choose the local envelope to be
\begin{equation}\label{eq:U_SecondRate}
U(\epsilon, X, \delta) := 2( |\epsilon| + 2\Phi) F_\delta(X)
\end{equation}
since it satisfies~\eqref{eq:U_def}. Recall that $T_B(f; \epsilon, X, \delta) = T(f; \epsilon, X)\mathbf{1}\{U(\epsilon, X; \delta) \le B\},$ where $T$ is defined in Appendix~\ref{sec:PeelingResult}.  Observe that by definition
\begin{equation}\label{eq:L_inf_bound_classical}
\|T_B(f; \cdot, \cdot ,\delta)\|_{\infty} \le B,
\end{equation}
and for $f$ satisfying $\|f - f_0\| \le \delta$,
\begin{equation}\label{eq:L_2bound_classical}
P\left(T^2_B(f; \cdot, \cdot ,\delta)\right)\le P(T^2(f; \cdot, \cdot, \delta)) \le 4\sigma^2\delta^2 + 4\Phi^2\delta^2 = (2\sigma + 2 \Phi)^2 \delta^2.
\end{equation}
Finally, from the calculations of Section 3.4.3.2 of~\cite{VdVW96} it follows that for any $\delta>0$, we have 
\begin{equation}\label{eq:Entropy_secondRate}
  N_{[\,]}(\eta, \{T_{B}(f; \epsilon, X,\delta): f\in\mathcal{F}\}, \|\cdot\|) \le N_{[\,]}(\eta/(2\sigma + 2\Phi), \mathcal{F}, \|\cdot\|).
\end{equation}
From Lemma 3.4.2 of~\cite{VdVW96}, we get
\begin{equation}\label{eq:vdvMaximalSecond}
  \mathbb{E}\left[\sup_{\delta/2 \le \|f - f_0\| \le \delta}\mathbb{G}_n(T_{B}(f; \epsilon, X,\delta))\right] \le CJ_{[\,]}(\sqrt{2}(2\sigma + 2\Phi)\delta)\left(1 + \frac{J_{[\,]}(\sqrt{2}(2\sigma + 2\Phi)\delta)B}{\sqrt{n}(2\sigma + 2 \Phi)^2 2\delta^2}\right),
\end{equation}
where by~\eqref{eq:Brack},
\begin{equation}\label{eq:J_def}
\begin{split}
  J_{[\,]}(\beta) &= \int_0^{\beta} \sqrt{\log N_{[\,]}(\eta, \{T_B(f; \epsilon, X):\,f\in\mathcal{F}\}, \|\cdot\|)}d\eta\\
  &\le \int_0^{\beta}\sqrt{A(\eta/(2\sigma+2\Phi))^{-\alpha}}d\eta\\
  &= A^{1/2}(2\sigma+2\Phi)^{\alpha/2}\beta^{1-\alpha/2}/(1-\alpha/2).
  % \\
  % % &\le \frac{A^{1/2}\delta^{1-\alpha/2}}{2-\alpha}.
  \end{split}
\end{equation}
Therefore $J_{[\,]}(\sqrt{2}(2\sigma + 2\Phi)\delta) \le A^{1/2}2^{3/2 -\alpha/4}(2\sigma+2\Phi)\delta^{1-\alpha/2}/(2-\alpha)$ and
  \begin{align*}
  \mathbb{E}\left[\sup_{\delta/2 \le \|f - f_0\| \le \delta}\mathbb{G}_n(T_{B}(f; \epsilon, X,\delta))\right] &\le C\frac{A^{1/2}(\sigma + \Phi)\delta^{1-\alpha/2}}{2 - \alpha} + C\frac{A\delta^{2 - \alpha}B}{(2 - \alpha)^2\sqrt{n}\delta^2}\\
  &\le CA^{1/2}\frac{(\sigma+\Phi)\delta^{1 - \alpha/2}}{2 - \alpha} + CAB\frac{ \delta^{- \alpha}}{(2 - \alpha)^2\sqrt{n}}\\
  &:= \phi_n(\delta, B).
  \end{align*}
With $\phi_n(\cdot, \cdot)$ defined, we will apply Theorem~\ref{thm:GeneralPeeling} to find the tail probability bound.\\ %The above $\phi_n(\cdot, \cdot)$ is identical to the one used in the proof of Theorem~\ref{thm:SecondMomentRate} (see~\eqref{eq:phi_secondmomentrate}).\\

  \noindent\textbf{Application of Theorem~\ref{thm:GeneralPeeling}:}
  % \end{enumerate}
 To apply Theorem~\ref{thm:GeneralPeeling}, we need $\beta,\, \gamma$, and $s_\gamma(\delta)$. From assumption~\eqref{eq:GeneralEnvelopeCondition_L2}, we can choose
\begin{equation}\label{eq:gamma_S_31}
\gamma = q \qquad \text{ and }\qquad s_q(\delta) = \mathbb{E}\big(|U(\epsilon, X, \delta)|^q\big) \le 4^q\|(|\epsilon| + \Phi)F_{\delta}\|_q^q \le 4^qC^q\Phi^{2q}\delta^{qs}.
\end{equation}
% \todo[inline]{Possibility of using $K_q^q +\delta^{qs}$ instead of $K_q^q +\Phi^q$ but need to check if it will be helpful!}
then $s_q(\delta)$ satisfies $\mathbb{E}\left[U^\gamma(\epsilon, X; \delta)\right] \le s_{\gamma}(\delta)$. We will chose $\tau$ later. Theorem~\ref{thm:GeneralPeeling} now implies that
\begin{align*}
\mathbb{P}\left(\|\widehat{f} - f_0\| \ge D\varepsilon_n\right) &\le 
\left(\frac{C}{\sqrt{n}(D\varepsilon_n)^2}\right)^\tau\sum_{k = 1}^{\infty}\frac{\phi_n^\tau(2^kD\varepsilon_n, B_k)}{2^{2k\tau}} + \left(\frac{ 8 C\sqrt{\tau}(\sigma + \Phi)}{\sqrt{n}(D\varepsilon_n)}\right)^\tau\\ &\qquad + \left(\frac{C\tau}{n(D\varepsilon_n)^2}\right)^\tau\sum_{k = 1}^{\infty} \frac{B_k^\tau}{2^{2k\tau}} + \frac{16}{(D\varepsilon_n)^2}\sum_{k = 0}^{\infty} \frac{s_{q}(2^kD\varepsilon_n)}{2^{2k}B_k^{q - 1}},\\
&\le \left(\frac{C}{\sqrt{n}(D\varepsilon_n)^2}\right)^\tau\sum_{k = 1}^{\infty}\frac{\phi_n^\tau(2^kD\varepsilon_n, B_k)}{2^{2k\tau}} + \left(\frac{ 8 C\sqrt{\tau}(\sigma + \Phi)}{\sqrt{n}(D\varepsilon_n)}\right)^\tau\\ &\qquad + \left(\frac{C\tau}{n(D\varepsilon_n)^2}\right)^\tau\sum_{k = 1}^{\infty} \frac{B_k^\tau}{2^{2k\tau}} + \frac{4^{q+2}C^q\Phi^{2q}}{(D\varepsilon_n)^2}\sum_{k = 1}^{\infty} \frac{(2^kD\varepsilon_n)^{qs}}{2^{2k}B_k^{q - 1}},
\end{align*}
From the definition of $\phi_n(\delta; B)$, write $\phi_n(\delta; B) = A_1(\delta) + BA_2(\delta)$, where
\begin{equation}\label{eq:Intro_A_1_A_2}
A_1(\delta) := \frac{KA^{1/2}(\sigma+\Phi)\delta^{1-\alpha/2}}{(2-\alpha)},\quad\mbox{and}\quad A_2(\delta) := \frac{KA\delta^{-\alpha}}{n^{1/2}(2-\alpha)^2}.
\end{equation}
This implies
\[
\phi_n^{\tau}(2^kD\varepsilon_n; B_k) \le 2^{\tau}A_1^{\tau}(2^kD\varepsilon_n) + 2^{\tau}B_k^{\tau}A_2^{\tau}(2^kD\varepsilon_n).
\]
Hence
\begin{align}\label{eq:Finding_B_before}
\begin{split}
\mathbb{P}\left(\|\widehat{f} - f_0\| \ge D\varepsilon_n\right) &\le 
\sum_{k = 1}^{\infty}\left(\frac{CA_1(2^kD\varepsilon_n)}{2^{2k}\sqrt{n}(D\varepsilon_n)^2}\right)^{\tau} + \left(\frac{ 8 C\sqrt{\tau}(\sigma + \Phi)}{\sqrt{n}(D\varepsilon_n)}\right)^\tau\\ &\qquad + \sum_{k = 1}^{\infty} \frac{C^{\tau}B_k^\tau}{2^{2k\tau}}\left(\frac{\tau}{n(D\varepsilon_n)^2} + \frac{A_2(2^kD\varepsilon_n)}{\sqrt{n}(D\varepsilon_n)^2}\right)^\tau\\ 
&\qquad+ \frac{4^{q+2}C^q\Phi^{2q}}{(D\varepsilon_n)^2}\sum_{k = 1}^{\infty} \frac{(2^kD\varepsilon_n)^{qs}}{2^{2k}B_k^{q - 1}},
\end{split}
\end{align}
for a  universal constant $C.$ We now choose $B_k$ to balance the summands of the last two terms
\[
{B_k}\left(\frac{C\tau}{n(2^kD\varepsilon_n)^2} + \frac{CA_2(2^kD\varepsilon_n)}{\sqrt{n}(2^kD\varepsilon_n)^2}\right) = \left(\frac{4^{q+2}C^q\Phi^{2q}}{(2^kD\varepsilon_n)^{2-qs}B_k^{q-1}}\right)^{1/\tau}.
\]
Equivalently,
\begin{equation}\label{eq:B_kChoice_SecondMomentRate_Proof}
B_k = \left(\frac{C\tau}{n(2^kD\varepsilon_n)^2} + \frac{CA_2(2^kD\varepsilon_n)}{\sqrt{n}(2^kD\varepsilon_n)^2}\right)^{-\tau/(\tau + q-1)}\left(\frac{4^{q+2}C^q\Phi^{2q}}{(2^kD\varepsilon_n)^{2-qs}}\right)^{1/(\tau + q-1)}.
\end{equation}
Hence the last two terms in~\eqref{eq:Finding_B_before} become
\begin{align}\label{eq:th625}
\begin{split}
&\sum_{k=1}^{\infty} \frac{4^{q+2}C^q\Phi^{2q}}{(2^kD\varepsilon_n)^{2-qs}}\left(\frac{4^{q+2}C^q\Phi^{2q}}{(2^kD\varepsilon_n)^{2-qs}}\right)^{-(q-1)/(\tau + q-1)}\left(\frac{C\tau}{n(2^kD\varepsilon_n)^2} + \frac{CA_2(2^kD\varepsilon_n)}{\sqrt{n}(2^kD\varepsilon_n)^2}\right)^{\tau(q-1)/(\tau + q-1)}\\
&= \sum_{k=1}^{\infty} \left(\frac{4^{q+2}C^q\Phi^{2q}}{(2^kD\varepsilon_n)^{2-qs}}\right)^{\tau/(\tau + q-1)}\left(\frac{C\tau}{n(2^kD\varepsilon_n)^2} + \frac{CA_2(2^kD\varepsilon_n)}{\sqrt{n}(2^kD\varepsilon_n)^2}\right)^{\tau(q-1)/(\tau + q-1)}\\
&= \sum_{k=1}^{\infty} \left(\frac{4^{q+2}C^q\Phi^{2q}}{(2^kD\varepsilon_n)^{2-qs}}\right)^{\tau/(\tau + q-1)}\left(\frac{C\tau}{n(2^kD\varepsilon_n)^2} + \frac{CKA}{{n}(2-\alpha)^2(2^kD\varepsilon_n)^{2+\alpha}}\right)^{\tau(q-1)/(\tau + q-1)},
\end{split}
\end{align}
where the last equality follows from the definition of $A_2(\cdot)$ in~\eqref{eq:Intro_A_1_A_2}. Substituting this in~\eqref{eq:Finding_B_before} and using the definition of $A_1(\cdot)$, we get
\begin{align}\label{eq:forextra_logterm}
\begin{split}
&\mathbb{P}\left(\|\widehat{f} - f_0\| \ge D\varepsilon_n\right)\\ &\le 
\sum_{k = 1}^{\infty}\left(\frac{CKA^{1/2}(\sigma+\Phi)}{\sqrt{n}(2^kD\varepsilon_n)^{1 + \alpha/2}(2-\alpha)}\right)^{\tau} + \left(\frac{ 8 C\sqrt{\tau}(\sigma + \Phi)}{\sqrt{n}(D\varepsilon_n)}\right)^\tau\\ &\quad+ \sum_{k=1}^{\infty} \left[\left(\frac{4^{q+2}C^q\Phi^{2q}}{(2^kD\varepsilon_n)^{2-qs}}\right)\left(\frac{C\tau}{n(2^kD\varepsilon_n)^2} + \frac{CKA}{{n}(2-\alpha)^2(2^kD\varepsilon_n)^{2+\alpha}}\right)^{(q-1)}\right]^{{\tau/(\tau + q-1)}}\\
% &\le \left(\frac{CKA^{1/2}(\sigma+\Phi)}{\sqrt{n}(D\varepsilon_n)^{1+\alpha/2}(2-\alpha)}\right)^{\tau}\sum_{k=1}^{\infty}\frac{1}{2^{k\tau(1+\alpha/2)}} + \left(\frac{8C\sqrt{\tau}(\sigma + \Phi)}{\sqrt{n}(D\varepsilon_n)}\right)^{\tau}\\
% &\qquad+ \left(\frac{C^{q-1}\tau^{q-1}4^{q+2}C^q\Phi^{2q}}{n^{q-1}(D\varepsilon_n)^{2-qs+2(q-1)}}\right)^{\tau/(\tau + q-1)}\sum_{k=1}^{\infty} \frac{1}{2^{\tau kq(2-s)/(\tau+q-1)}}\\
% &\qquad+ \left(\frac{(CKA)^{q-1}4^{q+2}C^q\Phi^{2q}}{n^{q-1}(2-\alpha)^2(D\varepsilon_n)^{2-qs+(2+\alpha)(q-1)}}\right)^{\tau/(\tau+q-1)}\sum_{k=1}^{\infty}\frac{1}{2^{\tau kq(2-s)/(\tau+q-1)}}\\
&\le \left(\frac{CKA^{1/2}(\sigma+\Phi)}{\sqrt{n}(D\varepsilon_n)^{1+\alpha/2}(2-\alpha)}\right)^{\tau} + \left(\frac{8C\sqrt{\tau}(\sigma + \Phi)}{\sqrt{n}(D\varepsilon_n)}\right)^{\tau}\\
&\qquad+ \left(\frac{C^{q-1}\tau^{q-1}4^{q+2}C^q\Phi^{2q}}{n^{q-1}(D\varepsilon_n)^{2-qs+2(q-1)}}\right)^{\tau/(\tau + q-1)} + \left(\frac{(CKA)^{q-1}4^{q+2}C^q\Phi^{2q}}{n^{q-1}(2-\alpha)^2(D\varepsilon_n)^{2-qs+(2+\alpha)(q-1)}}\right)^{\tau/(\tau+q-1)}.
\end{split}
\end{align}In the inequalities above the constant $C$ could be different in different lines. Now choose $\varepsilon_n$ so that the following inequalities are satisfied:
\begin{align*}
\frac{A^{1/2}(\sigma+\Phi)}{\sqrt{n}\varepsilon_n^{1+\alpha/2}} \le 1\quad&\Leftrightarrow\quad \varepsilon_n \ge A^{1/(2+\alpha)}(\sigma+\Phi)^{2/(2+\alpha)}n^{-1/(2+\alpha)},\\
\frac{(\sigma + \Phi)}{\sqrt{n}\varepsilon_n} \le 1\quad&\Leftrightarrow\quad \varepsilon_n \ge {(\sigma + \Phi)}{n^{-1/2}},\\
\frac{C^q\Phi^{2q}}{n^{q-1}\varepsilon_n^{q(2-s)}} \le 1\quad&\Leftrightarrow\quad \varepsilon_n \ge {C^{1/(2-s)}\Phi^{2/(2-s)}}{n^{-(q-1)/(q(2-s))}},\\
\frac{A^{q-1}C^q\Phi^{2q}}{n^{q-1}\varepsilon_n^{2-qs + (2+\alpha)(q-1)}} \le 1\quad&\Leftrightarrow\quad \varepsilon_n \ge (A^{q-1}C^q\Phi^{2q})^{1/(2-qs+(2+\alpha)(q-1))}n^{-1/(2+\alpha+(2-qs)/(q-1))}.
\end{align*}
Take
\[
\varepsilon_n := \max\left\{\frac{(\sigma+\Phi)^{2/(2+\alpha)}}{(nA^{-1})^{1/(2+\alpha)}},
\frac{(\sigma + \Phi)}{n^{1/2}},
\frac{C^{1/(2-s)}\Phi^{2/(2-s)}}{n^{(q-1)/(q(2-s))}},
\frac{(A^{q-1}C^q\Phi^{2q})^{1/(2-qs+(2+\alpha)(q-1))}}{n^{1/(2+\alpha+(2-qs)/(q-1))}}\right\},
\]
for which, the tail bound becomes
\begin{align*}
\mathbb{P}(\|\widehat{f} - f_0\| \ge D\varepsilon_n) &\le \left(\frac{CK}{D^{1+\alpha/2}(2-\alpha)}\right)^{\tau} + \left(\frac{8C\sqrt{\tau}}{D}\right)^{\tau} + \left(\frac{C^{q-1}\tau^{q-1}4^{q+2}}{D^{q(2-s)}}\right)^{\tau/(\tau+q-1)}\\
&\qquad+ \left(\frac{(CK2^{2-\alpha})^{q-1}4^{q+2}}{(2-\alpha)^2D^{2-qs+(2+\alpha)(q-1)}}\right)^{\tau/(\tau+q-1)}.
\end{align*}
Take $\tau$ such that $\tau \ge q$ and $\tau q(2-s)/(\tau+q-1) \ge q$ or equivalently $\tau \ge \max\{q, (q-1)/(1-s)\}$. In case $s = 1$, fix any $\eta > 0$ and take $\tau$ such that $\tau \ge q$ and $\tau q/(\tau + q-1) = q - \eta$ or equivalently $\tau \ge \max\{q, (q-1)(q-\eta)/\eta\}$. This would imply that for all $D > 0$,
\begin{equation}\label{eq:TrueTailBound_L2}
\mathbb{P}(\|\widehat{f} - f_0\| \ge D\varepsilon_n) \le {C} D^{-q+\eta\mathbf{1}\{s=1\}},
\end{equation}for a constant ${C} > 0$ depending only on $q, s, \alpha,$ and $\eta$. Because $A\ge1$, $\varepsilon_n^{-1}$ is equal to $r_n$ in~\eqref{eq:rn_SecondMoment}.

\subsection{Additional log factors in~(\ref{eq:GeneralEnvelopeCondition_L2})} % (fold)
\label{sub:additional_log_factors}
Suppose $F_\delta$ and $\epsilon$ do not satisfy~\eqref{eq:GeneralEnvelopeCondition_L2} but satisfy \begin{equation}\label{eq:GeneralEnvelopeCondition_L2_extra}
\big\|(|\epsilon| +\Phi) F_{\delta}(X)\big\|_q \le C \Phi^2 \delta^{s} \log^{\nu}(1/\delta). 
\end{equation}
Then, by modifying $s_q$ in~\eqref{eq:gamma_S_31} and incorporating the changes in~\eqref{eq:Finding_B_before},~\eqref{eq:B_kChoice_SecondMomentRate_Proof}, and~\eqref{eq:th625}, it can be shown that $\widehat{f}$ will satisfy the following modification of~\eqref{eq:forextra_logterm},
\begin{align*}\label{eq:test567}
&\mathbb{P}\left(\|\widehat{f} - f_0\| \ge D\varepsilon_n\right)\\ &\le \left(\frac{CKA^{1/2}(\sigma+\Phi)}{\sqrt{n}(D\varepsilon_n)^{1+\alpha/2}(2-\alpha)}\right)^{\tau} + \left(\frac{8C\sqrt{\tau}(\sigma + \Phi)}{\sqrt{n}(D\varepsilon_n)}\right)^{\tau}\\
&\qquad+ \left(\frac{C^{q-1}\tau^{q-1}4^{q+2}C^q\Phi^{2q} ( \log(1/D \epsilon_n))^{\nu}}{n^{q-1}(D\varepsilon_n)^{2-qs+2(q-1)}}\right)^{\tau/(\tau + q-1)} + \left(\frac{(CKA)^{q-1}4^{q+2}C^q\Phi^{2q}( \log(1/D \epsilon_n))^{\nu}}{n^{q-1}(2-\alpha)^2(D\varepsilon_n)^{2-qs+(2+\alpha)(q-1)}}\right)^{\tau/(\tau+q-1)}.
\end{align*}
We will now chose $\varepsilon_n$ that satisfies the following inequalities
\begin{align*}
\frac{A^{1/2}(\sigma+\Phi)}{\sqrt{n}\varepsilon_n^{1+\alpha/2}} \le 1\quad&\Leftrightarrow\quad \varepsilon_n \ge A^{1/(2+\alpha)}(\sigma+\Phi)^{2/(2+\alpha)}n^{-1/(2+\alpha)},\\
\frac{(\sigma + \Phi)}{\sqrt{n}\varepsilon_n} \le 1\quad&\Leftrightarrow\quad \varepsilon_n \ge {(\sigma + \Phi)}{n^{-1/2}},\\
\frac{C^q\Phi^{2q} ( \log(1/\epsilon_n))^{\nu} }{n^{q-1}\varepsilon_n^{q(2-s)}} \le 1\quad&\Leftrightarrow\quad \varepsilon_n \ge {C^{1/(2-s)}\Phi^{2/(2-s)}}{n^{-(q-1)/(q(2-s))}} (\log n)^{\nu/(q(2-s))},\\
\frac{A^{q-1}C^q\Phi^{2q} ( \log(1/ \epsilon_n))^{\nu}}{n^{q-1}\varepsilon_n^{2-qs + (2+\alpha)(q-1)}} \le 1\quad&\Leftrightarrow\quad \varepsilon_n \ge (A^{q-1}C^q\Phi^{2q} (\log n)^{\nu})^{1/(2-qs+(2+\alpha)(q-1))}n^{-1/(2+\alpha+(2-qs)/(q-1))}.
\end{align*}
Thus by choosing the $\tau$ as before, we have that $\widehat{f}$ will satisfy~\eqref{eq:TrueTailBound_L2} with 
{\small \begin{equation}\label{eq:ep_n_extra_log}
\varepsilon_n := \max\left\{\frac{(\sigma+\Phi)^{2/(2+\alpha)}}{(nA^{-1})^{1/(2+\alpha)}},
\frac{C^{1/(2-s)}\Phi^{2/(2-s)} }{[n^{(q-1)} (\log n)^{-\nu}]^{1/(q(2-s))}},
\left(\frac{A^{q-1}\Phi^{2q} (\log n)^{\nu}}{n^{q-1}}\right)^{1/(2-qs+(2+\alpha)(q-1))}\right\}.
\end{equation}}
% subsection additional_log_factors (end)

\section{Proof of Theorem~\ref{thm:L_inf_smooth}} % (fold)
\label{sec:proof:L_inf_smooth}
As in the proof of Theorem~\ref{thm:SecondMomentRate}, we will first find an appropriate $\phi_n(\delta, B)$. We will use the maximal inequality derived in~Proposition~\ref{lem:bernouli_maximal}  in conjunction with techniques borrowed from~\cite{MR1411488} (also see~\cite[Theorem 3.5]{Dirksen})  to find this.
Recall that $T(f; \epsilon, X) := 2\epsilon(f - f_0)(X) - (f - f_0)^2(X).$
By symmetrization and contraction (Theorem 3.1.21 and Corollary 3.2.2 of~\cite{Gine16}, respectively), we get that 
\begin{equation}\label{eq:contrac_sym_xf}
\phi_n(\delta, B) \le \E \left(\sup_{f\in \F_\delta} \mathbb{G}_nT(f; \epsilon, X, \delta)\right) \le \E \left(\sup_{f\in \F_\delta} \mathbb{G}_n (\epsilon+ 16 \Phi) (f-f_0)\right).
\end{equation}
 Define a stochastic process $X(\cdot)$ on $\F$ as 
\[
X(f) := \mathbb{G}_n (\epsilon+ 16 \Phi) (f-f_0).
\]
We will now bound $\E (\sup_{f\in \F_\delta} X(f))$. Let $\{f_0\} = S_0 \subset S_1 \subset \cdots \subset S_m \subset \cdots$ be a sequence of incremental subsets of $\mathcal{F}_{\delta}$. We take these sets $S_i$ so that $\log|S_{i+1}| \le 2^{i+1}$. Let $\varepsilon_{2,i}$ denote the smallest $\varepsilon$ so that $\log N(\varepsilon_{2,i}, \mathcal{F}_{\delta}, \|\cdot\|) \le 2^i$ and $\varepsilon_{\infty, i}$ denote the smallest $\varepsilon$ so that $\log N(\varepsilon_{\infty,i}, \mathcal{F}_{\delta}, \|\cdot\|_{\infty}) \le 2^{i}$. Now set $A_i = \{f_1, \ldots, f_{2^{2^i}}\}$ as the $\varepsilon_{2,i}$-net of $\mathcal{F}_{\delta}$ with respect to $\|\cdot\|$ and $B_i = \{g_1,\ldots,g_{2^{2^i}}\}$ as the $\varepsilon_{\infty,i}$-net of $\mathcal{F}_{\delta}$ with respect to $\|\cdot\|_{\infty}$. Define the partition of $\mathcal{F}_{\delta}$ by
\[
\{B_2(f_j,\varepsilon_{2,i})\cap B_{\infty}(g_k, \varepsilon_{\infty,i}):\,1\le j,k\le 2^{2^i}\}.
\]
The number of sets in this partition is bounded above by $(2^{2^i})^2 = 2^{2^{i+1}}$. Take $S_{i+1}\subset \mathcal{F}_{\delta}$ of cardinality at most $2^{2^{i+1}}$ by taking one element in each of $B_{2}(f_j, \varepsilon_{2,i})\cap B_{\infty}(g_k, \varepsilon_{\infty,i})$ for all $1\le j,k\le 2^{2^i}$. For any $f\in\mathcal{F}_{\delta}$, we let $\pi_if$ be the element in $S_i$ that is closest to $f$ so that
\[
\max_{f\in\mathcal{F}_{\delta}}\|f - \pi_if\| \le \varepsilon_{2,i-1}\quad\mbox{and}\quad \max_{f\in\mathcal{F}_{\delta}}\|f - \pi_if\|_{\infty} \le \varepsilon_{\infty,i-1}.
\]
Using these $\pi_if$, we can write
% \[
% X(f) - X(f_0) = X(f) - X(\pi_Tf) + \sum_{i=1}^{T} \{X(\pi_if) - X(\pi_{i-1}f)\}.
% \]
\[
X(f) - X(f_0) =  \sum_{i\ge 1} \{X(\pi_if) - X(\pi_{i-1}f)\}.
\]
Thus 
% \begin{align}
% \E \left(\sup_{f\in \F_\delta} X(f)\right) &= \E \left(\sup_{f\in \F_\delta}[X(f) - X(f_0)]\right) \label{eq:chaining_Xf} \\
% &\le  \E \left(\sup_{f\in \F_\delta} [X(f) - X(\pi_Tf)]\right) + \sum_{t=1}^{T} \E \left(\sup_{f\in \F_\delta} \{X(\pi_tf) - X(\pi_{t-1}f)\}\right).
% \end{align}
\begin{align}
\E \left[\sup_{f\in \F_\delta} X(f)\right] = \E \left[\sup_{f\in \F_\delta}[X(f) - X(f_0)]\right] 
&\le  \sum_{t\ge 1} \E \left[\sup_{f\in \F_\delta} \{X(\pi_tf) - X(\pi_{t-1}f)\}\right].\label{eq:chaining_Xf} 
\end{align}
% Observe now that by symmetrization and Cauchy-Schwarz inequality
% \begin{align}\label{eq:tail_part}
% \begin{split}
%   \mathbb{E}\left[\sup_{f\in\mathcal{F}_{\delta}}|X(f) - X(\pi_Tf)|\right]
% ={}& \mathbb{E}\left[\sup_{f\in\mathcal{F}_{\delta}}\Big|\mathbb{G}_n (\epsilon + 16 \Phi) [f- \pi_Tf]\Big|\right]\\
% \le{}& 2 \mathbb{E}\left[\sup_{f\in\mathcal{F}_{\delta}}\Big|\frac{1}{\sqrt{n}}\sum_{i=1}^n R_i(\epsilon_i + 16 \Phi) [f(X_i)- \pi_Tf(X_i)]\Big|\right]\\
% % ={}&2\mathbb{E}\left[\sup_{f\in\mathcal{F}_{\delta}}\left|\frac{1}{\sqrt{n}}\sum_{i=1}^n R_i[(2\epsilon_i + f + \pi_Tf - 2f_0)\mathbf{1}\{U(\epsilon, X; \delta) \le B\}](f - \pi_Tf)(X_i)\Big|\right]\\
% \le{}& 2\mathbb{E}\left[\left(\sum_{i=1}^n (|\epsilon_i| + 16\Phi)^2\right)^{1/2}\sup_{f\in\mathcal{F}_{\delta}}\|f - \pi_Tf\|_{n}\right]\\
% \le{}& 2n^{1/2}(\sigma + 16\Phi)\varepsilon_{\infty,T-1}.
% \end{split}
% \end{align}
% We will now bound each of the expectations in the second term. 
 By symmetrization,  for every $t \ge 1$ we have 
\begin{align}
\mathbb{E}\left[\sup_{f\in\mathcal{F}_{\delta}}|X(\pi_{t}f) - X(\pi_{t-1}f)|\right]\label{eq:sym_chain_part}
={}& \mathbb{E}\left[\sup_{f\in\mathcal{F}_{\delta}}\Big|\mathbb{G}_n (\epsilon + 16 \Phi) [\pi_{t}f- \pi_{t-1}f]\Big|\right]\\
\le{}& 2 \frac{1}{\sqrt{n}}\mathbb{E}\Big[\sup_{f\in\mathcal{F}_{\delta}}\Big|\sum_{i=1}^n R_i(\epsilon_i + 16 \Phi) [\pi_{t}f(X_i)- \pi_{t-1}f(X_i)]\Big|\Big].
\end{align}
Observe that the number of possible pairs $(\pi_t f, \pi_{t-1}f)$ is bounded by $|S_t| |S_{t-1}|$. Let $N_{t+1}= 2^{2^{t+1}}$ and  $\{(g_j, h_j)\}_{j=1}^{N_t}$ denote all such pairs.  Thus in~\eqref{eq:sym_chain_part}, the supremum is over a finite set.  Thus 
\begin{align}
&\frac{1}{\sqrt{n}}\mathbb{E}\left[\sup_{f\in\mathcal{F}_{\delta}}\Big|\sum_{i=1}^n R_i(\epsilon_i + 16 \Phi) [\pi_{t}f(X_i)- \pi_{t-1}f(X_i)]\Big|\right] \\
\le{}&\frac{1}{\sqrt{n}}\mathbb{E}\left[\max_{1\le j \le N_{t+1}}\Big|\sum_{i=1}^n R_i(\epsilon_i + 16 \Phi) [g_j(X_i)- h_j(X_i)]\Big|\right]. \label{eq:sup_to_max}
\end{align}
To bound the above expectation, we will use the maximal inequality in Proposition~\ref{lem:bernouli_maximal} with 
\begin{equation}\label{eq:xi_i_def}
x_{i,j}=(\epsilon_i + 16 \Phi) [g_j(X_i)- h_j(X_i)]\qquad\text{and}\qquad \xi_i:= \max_{1\le j\le N_{t+1}} |(\epsilon_i + 16 \Phi) [g_j(X_i)- h_j(X_i)]|.
\end{equation}
Observe that 
\begin{equation}\label{eq:sup_exp_bound}
\mathbb{E}\left[\sum_{i=1}^n x_{i,j}^2\right] \le \mathbb{E}\left[\sum_{i=1}^n (\epsilon_i + 16 \Phi) [g_j(X_i)- h_j(X_i)]^2\right] \le (\sigma + 16 \Phi)^2 n \|g_j-h_j\|
\end{equation}
and 
\begin{equation}\label{eq:xi_pmom}
\mathbb{E}[\xi_i^q] =\mathbb{E}\left[\Big|(\epsilon_i + 16 \Phi)^q \max_{1\le j\le N_{t+1}}  [g_j(X_i)- h_j(X_i)]^q\Big|\right] \le (\|\epsilon\|_q + 16 \Phi)^q\max_{1\le j\le N_{t+1}} \|g_j- h_j\|_{\infty}^q.
\end{equation}
Thus by~Proposition~\ref{lem:bernouli_maximal}, we have that 
\begin{align}\label{eq:apply_bernouli}
\begin{split}
&\frac{1}{\sqrt{n}}\mathbb{E}\left[\max_{1\le j \le N_{t+1}}\Big|\sum_{i=1}^n R_i(\epsilon_i + 16 \Phi) [g_j(X_i)- h_j(X_i)]\Big|\right] \\
\le{}&(\sigma + 16 \Phi)\sqrt{\log N_{t+1}}  \, \max_{1\le j \le N_{t+1}}  \|g_j- h_j\| \\
&\qquad+ \left(\log N_{t+1}\right)^{1 - 1/q} n^{1/q-1/2}(\|\epsilon\|_q + 16 \Phi)\max_{1\le j\le N_{t+1}} \|g_j- h_j\|_{\infty}.
\end{split}
\end{align}
We will now bound $\max_{1\le j \le N_{t+1}}  \|g_j- h_j\|$ and $\max_{1\le j \le N_{t+1}}  \|g_j- h_j\|_{\infty}.$ By definition of  $\{(g_j, h_j)\}_{j=1}^{N_t}$, we have that
\[ \max_{1\le j \le N_{t+1}}  \|g_j- h_j\| = \sup_{f\in \F_\delta} \|\pi_t f- \pi_{t-1}f\| \le \sup_{f\in \F_\delta} \| f- \pi_{t}f\|+ \sup_{f\in \F_\delta} \|f- \pi_{t-1}f\| \le 2 \varepsilon_{2, t-2}.\]
Similarly, we also have 
\[\max_{1\le j \le N_{t+1}}  \|g_j- h_j\|_{\infty} \le 2 \varepsilon_{\infty, t-2}.\]
Thus combining~\eqref{eq:sym_chain_part},~\eqref{eq:sup_to_max}, and~\eqref{eq:apply_bernouli}, we have that 
\begin{align*}
 \sum_{i\ge 1} \mathbb{E}\left[\sup_{f\in\mathcal{F}_{\delta}}|X(\pi_if) - X(\pi_{i-1}f)|\right]
&\lesssim (\sigma+\Phi)\sum_{t\ge 1} \varepsilon_{2,t-2}\sqrt{\log N_{t+1}}\\ &\qquad+ \frac{\|\epsilon\|_q+ \Phi}{n^{1/2 - 1/q}}\sum_{t\ge 1} \varepsilon_{\infty,t-2}(\log N_{t+1})^{1-1/q}\\
&\lesssim (\sigma+\Phi)\sum_{t\ge 1} 2^{t/2}\varepsilon_{2,t-2} + \frac{\|\epsilon\|_q+ \Phi}{n^{1/2 - 1/q}}\sum_{t\ge 1} 2^{t(1-1/q)}\varepsilon_{\infty,t-2},
\end{align*}
where $\varepsilon_{2,-1} = \delta$ and $\varepsilon_{2,-1} = \|F_\delta\|_{\infty}.$ Thus 
\begin{align}\label{eq:chaining_final}
\begin{split}
\E \left(\sup_{f\in \F_\delta} \mathbb{G}_nT(f; \epsilon, X, \delta)\right)
\le{}&\sum_{i\ge 1} \mathbb{E}\left[\sup_{f\in\mathcal{F}_{\delta}}|X(\pi_if) - X(\pi_{i-1}f)|\right]\\
\lesssim{}& (\sigma+\Phi)\Big(\delta+ \sum_{t\ge 2} 2^{t/2}\varepsilon_{2,t-2}\Big) + \frac{\|\epsilon\|_q+ \Phi}{n^{1/2 - 1/q}}\Big( \|F_\delta\|_{\infty}+ \sum_{t\ge 2} 2^{t(1-1/q)}\varepsilon_{\infty,t-2}\Big)\\
\lesssim{}& (\sigma+\Phi)\Big(\delta+ \sum_{t\ge 0} 2^{t/2}\varepsilon_{2,t}\Big) + \frac{\|\epsilon\|_q+ \Phi}{n^{1/2 - 1/q}}\Big( \|F_\delta\|_{\infty}+ \sum_{t\ge 0} 2^{t(1-1/q)}\varepsilon_{\infty,t}\Big)
\end{split}
\end{align}
By (2.37) of~\cite[Page 22]{MR3184689}, we have that 
\begin{equation}\label{eq:e2_sum}
\sum_{t\ge 0} 2^{t/2}\varepsilon_{2,t} \le \int_{0}^\delta \sqrt{\log N(\eta, \F_\delta, \|\cdot\|)} \, d\eta \text{  and  } \sum_{t\ge 0} 2^{t/2}\varepsilon_{\infty,t} \le \int_{0}^{\|F_\delta\|_{\infty}} \{\log N(\eta, \F_\delta, \|\cdot\|_{\infty})\}^{1-1/q} \, d\eta. 
\end{equation}
Thus if $\log N(\eta, \F_\delta, \|\cdot\|_{\infty}) \le A \eta^{-\alpha}$ for some $\alpha \in [0, 2)$. Then 
\[ \sum_{t\ge 0} 2^{t/2}\varepsilon_{2,t} \le A^{1/2} \delta^{1-\alpha/2 }/(1-\alpha/2)\]
and if $\alpha(1-1/q) < 1$ then 
\[\sum_{t\ge 0} 2^{t/2}\varepsilon_{\infty,t} \le  A^{1-1/q}  \|F_\delta\|_{\infty}^{1- \alpha(1-1/q)}. \]
Thus if $\alpha(1-1/q) < 1$, then 
\begin{align}
\Phi_n(\delta, B)&\le C (\sigma+\Phi)\Big(\delta+  \frac{A^{1/2} \delta^{1-\alpha/2 }}{1-\alpha/2}\Big) + C \frac{\|\epsilon\|_q+ \Phi}{n^{1/2 - 1/q}}\Big( \|F_\delta\|_{\infty}+ A^{1-1/q}  \|F_\delta\|_{\infty}^{1- \alpha(1-1/q)}\Big)\label{eq:L_infy_phi_n}\\
&\le C (\sigma+\Phi)\Big(\delta+  \frac{A^{1/2} \delta^{1-\alpha/2 }}{1-\alpha/2}\Big) + C \frac{\|\epsilon\|_q+ \Phi}{n^{1/2 - 1/q}}\Big( \Phi^{1-s} \delta^{s}+ A^{1-1/q}  (\Phi^{1-s} \delta^s)^{1- \alpha(1-1/q)}\Big).
\end{align}

\noindent\textbf{Application of Theorem~\ref{thm:GeneralPeeling}:}
  % \end{enumerate}
 To apply Theorem~\ref{thm:GeneralPeeling}, we need $\tau,\, \gamma$, and $s_\gamma(\delta)$. From assumption~\eqref{eq:GeneralEnvelopeCondition_L_infty_linf}, we can choose
\begin{equation}\label{eq:gamma_S_3}
\gamma = q \qquad \text{ and }\qquad s_q(\delta) = \mathbb{E}\big(|U(\epsilon, X, \delta)|^q\big) \le 4^q\|(|\epsilon| + \Phi)F_{\delta}\|_q^q \le C\Phi^{q+q(1-s)}\delta^{qs}.
\end{equation}
% \todo[inline]{Possibility of using $K_q^q +\delta^{qs}$ instead of $K_q^q +\Phi^q$ but need to check if it will be helpful!}
then $s_q(\delta)$ satisfies the conditions of Theorem~\ref{thm:GeneralPeeling}. We will chose $\tau$ later. Thus by~\eqref{eq:Error_of_fhat}, we have that 
\begin{align}\label{eq:probbound_beforeB_L-inf_smooth}
\begin{split}
\mathbb{P}\left(\|\widehat{f} - f_0\| \ge D\varepsilon_n\right) &\le 
\left(\frac{C}{\sqrt{n}(D\varepsilon_n)^2}\right)^{\tau}\sum_{k = 1}^{\infty}\frac{\phi_n^{\tau}(2^kD\varepsilon_n, B_k)}{2^{2k\tau}} + \left(\frac{ 8 C\sqrt{\tau}(\sigma + \Phi)}{\sqrt{n}(D\varepsilon_n)}\right)^{\tau}\\ &\qquad + \left(\frac{C\tau}{n(D\varepsilon_n)^2}\right)^{\tau}\sum_{k = 1}^{\infty} \frac{B_k^\tau}{2^{2k\tau}} + \frac{16}{(D\varepsilon_n)^2}\sum_{k = 0}^{\infty} \frac{s_{q}(2^kD\varepsilon_n)}{2^{2k}B_k^{q - 1}}\\
&\le \left(\frac{C}{\sqrt{n}(D\varepsilon_n)^2}\right)^\tau\sum_{k = 1}^{\infty}\frac{\phi_n^\tau(2^kD\varepsilon_n, B_k)}{2^{2k\tau}} + \left(\frac{ 8 C\sqrt{\tau}(\sigma + \Phi)}{\sqrt{n}(D\varepsilon_n)}\right)^\tau\\ &\qquad + \left(\frac{C\tau}{n(D\varepsilon_n)^2}\right)^\tau\sum_{k = 1}^{\infty} \frac{B_k^\tau}{2^{2k\tau}} + \frac{4^{q+2}C^q\Phi^{2q}}{(D\varepsilon_n)^2}\sum_{k = 1}^{\infty} \frac{(2^kD\varepsilon_n)^{qs}}{2^{2k}B_k^{q - 1}}.
\end{split}
\end{align}  
We will now choose $B_k$ to minimize the upper bound and then choose the smallest $\varepsilon_n$ such that the right hand side is a does not depend on $n$ and goes to zero as $D$ increases to infinity.
Substituting~\eqref{eq:L_infy_phi_n} in the probability bound~\eqref{eq:probbound_beforeB_L-inf_smooth}, we get
\begin{align}\label{eq:Finding_B_before_L_infty_proof}
\begin{split}
\mathbb{P}\left(\|\widehat{f} - f_0\| \ge D\varepsilon_n\right) &\le 
\sum_{k = 1}^{\infty}\left(\frac{CA_1(2^kD\varepsilon_n)}{2^{2k}\sqrt{n}(D\varepsilon_n)^2}\right)^{\tau} + \left(\frac{ 8 C\sqrt{\tau}(\sigma + \Phi)}{\sqrt{n}(D\varepsilon_n)}\right)^\tau\\ &\qquad + \sum_{k = 1}^{\infty} \frac{C^{\tau}B_k^\tau}{2^{2k\tau}}\left(\frac{\tau}{n(D\varepsilon_n)^2} \right)^\tau+ \frac{4^{q+2}C^q\Phi^{2q}}{(D\varepsilon_n)^2}\sum_{k = 1}^{\infty} \frac{(2^kD\varepsilon_n)^{qs}}{2^{2k}B_k^{q - 1}},
\end{split}
\end{align}
for some universal constant $C$ and 
\begin{equation}\label{eq:A1_def_L_inf}
A_1(\delta):= C (\sigma+\Phi)\Big(\delta+  \frac{A^{1/2} \delta^{1-\alpha/2 }}{1-\alpha/2}\Big) + C \frac{\|\epsilon\|_q+ \Phi}{n^{1/2 - 1/q}}\Big( \Phi^{1-s} \delta^{s}+ A^{1-1/q}  (\Phi^{1-s} \delta^s)^{1- \alpha(1-1/q)}\Big).
\end{equation}
We now see that the choice of  $B_k$ that balances the summands of the last two terms is
\begin{equation}\label{eq:B_kChoice_SecondMomentRate_Proof1}
B_k = \left(\frac{C\tau}{n(2^kD\varepsilon_n)^2} \right)^{-\tau/(\tau + q-1)}\left(\frac{4^{q+2}C^q\Phi^{2q}}{(2^kD\varepsilon_n)^{2-qs}}\right)^{1/(\tau + q-1)}.
\end{equation}
Hence the last two terms in~\eqref{eq:Finding_B_before_L_infty_proof} become
\begin{align}\label{eq:th6251}
\begin{split}
&\sum_{k=1}^{\infty} \frac{4^{q+2}C^q\Phi^{2q}}{(2^kD\varepsilon_n)^{2-qs}}\left(\frac{4^{q+2}C^q\Phi^{2q}}{(2^kD\varepsilon_n)^{2-qs}}\right)^{-(q-1)/(\tau + q-1)}\left(\frac{C\tau}{n(2^kD\varepsilon_n)^2} \right)^{\tau(q-1)/(\tau + q-1)}\\
&= \sum_{k=1}^{\infty} \left(\frac{4^{q+2}C^q\Phi^{2q}}{(2^kD\varepsilon_n)^{2-qs}}\right)^{\tau/(\tau + q-1)}\left(\frac{C\tau}{n(2^kD\varepsilon_n)^2} \right)^{\tau(q-1)/(\tau + q-1)}
\end{split}
\end{align}
Substituting this  in the probability bound in~\eqref{eq:Finding_B_before_L_infty_proof} and simplifying, we obtain
\begin{align}
&\mathbb{P}\left(\|\widehat{f} - f_0\| \ge D\varepsilon_n\right)\nonumber\\ &\le 
\sum_{k = 1}^{\infty}\left(\frac{CA_1(2^kD\varepsilon_n)}{2^{2k}\sqrt{n}(D\varepsilon_n)^2}\right)^{\tau} + \left(\frac{ 8 C\sqrt{\tau}(\sigma + \Phi)}{\sqrt{n}(D\varepsilon_n)}\right)^\tau\label{eq:B_kL_inf_L2}\\ &\quad+ \sum_{k=1}^{\infty} \left(\frac{4^{q+2}C^q\Phi^{2q}}{(2^kD\varepsilon_n)^{2-qs}}\right)^{\tau/(\tau + q-1)}\left(\frac{C\tau}{n(2^kD\varepsilon_n)^2} \right)^{\tau(q-1)/(\tau + q-1)}.\nonumber
\end{align}
Now we use the definitions of $A_1(\cdot)$ in~\eqref{eq:A1_def_L_inf}, to get
\begin{align*}
&\mathbb{P}\left(\|\widehat{f} - f_0\| \ge D\varepsilon_n\right)\\
 \le{}& 
\sum_{k=1}^{\infty} \left(\frac{C(\sigma+ \Phi)}{\sqrt{n}(2^kD\varepsilon_n)}\right)^{\tau}  + \sum_{k=1}^{\infty}\left(\frac{C \sqrt{A} (\sigma+ \Phi) }{(1-\alpha/2)\sqrt{n}(2^kD\varepsilon_n)^{1+\alpha/2}}\right)^{\tau} +  \sum_{k=1}^{\infty}C \left(\frac{ \Phi^{1-s} (\|\epsilon\|_q+ \Phi)}{n^{1 - 1/q} (2^k D \epsilon_n)^{2-s}}\right)^\tau \\
 &+  \sum_{k=1}^{\infty}C \left(\frac{   A^{1-1/q}  \Phi^{(1-s) (1- \alpha(1-1/q))}(\|\epsilon\|_q+ \Phi)}{n^{1 - 1/q} (2^k D \epsilon_n)^{2-s+ \alpha s (1-1/q)}}\right)^\tau +\left(\frac{ C\sqrt{\tau}(\sigma + \Phi)}{\sqrt{n}(D\varepsilon_n)}\right)^\tau\\
 &+\sum_{k=1}^{\infty} \left(\frac{\Phi^{2q} \tau^{q-1}}{ n^{q-1}(2^kD\varepsilon_n)^{2q-qs}}\right)^{\tau/(\tau + q-1)}\\
\le{}&  \left(\frac{C \sqrt{A} (\sigma+ \Phi) 2^{1-\alpha/2}}{(1-\alpha/2)\sqrt{n}(D\varepsilon_n)^{1+\alpha/2}}\right)^{\tau} +C \left(\frac{ \Phi^{1-s} (\|\epsilon\|_q+ \Phi)}{n^{1 - 1/q} ( D \epsilon_n)^{2-s}}\right)^\tau  \\
 &+ \left(\frac{   A^{1-1/q}  \Phi^{(1-s) (1- \alpha(1-1/q))}(\|\epsilon\|_q+ \Phi)}{n^{1 - 1/q} ( D \epsilon_n)^{2-s+ \alpha s (1-1/q)}}\right)^\tau + \left(\frac{ 8 C\sqrt{\tau}(\sigma + \Phi)}{\sqrt{n}(D\varepsilon_n)}\right)^\tau+ \left(\frac{\Phi^{2q} \tau^{q-1}}{ n^{q-1}(D\varepsilon_n)^{2q-qs}}\right)^{\tau/(\tau + q-1)}.
\end{align*}
We now choose $\varepsilon_n$ so that the following inequalities are satisfied:
\begin{align*}
\frac{A^{1/2}(\sigma + \Phi)}{n^{1/2}\varepsilon_n^{1 + \alpha/2}} \le 1\quad&\Leftrightarrow\quad\varepsilon_n \ge \frac{A^{1/(2+\alpha)}(\sigma + \Phi)^{2/(2+\alpha)}}{n^{1/(2+\alpha)}},\\
\frac{ \Phi^{1-s} (\|\epsilon\|_q+ \Phi)}{n^{1 - 1/q} ( D \epsilon_n)^{2-s}} \le 1 \quad &\Leftrightarrow\quad \varepsilon_n \ge \frac{ \Phi}{n^{(q - 1)/q(2-s)} },\\
\frac{   A^{1-1/q}  \Phi^{(1-s) (1- \alpha(1-1/q))}(\|\epsilon\|_q+ \Phi)}{n^{1 - 1/q} ( D \epsilon_n)^{2-s+ \alpha s (1-1/q)}}\le 1 \quad &\Leftrightarrow\quad \varepsilon_n \ge \left(\frac{A^{q-1}\Phi^{(2q-qs+ \alpha s(q-1) -\alpha(q-1))}}{n^{q-1}}\right)^{1/(2q-qs+ \alpha s(q-1) )},\\
\frac{\sigma+\Phi}{\sqrt{n}\varepsilon_n} \le 1\quad&\Leftrightarrow\quad \varepsilon_n \ge \frac{(\sigma + \Phi)}{n^{1/2}},\\
\frac{C^q\Phi^{2q}}{n^{q-1}\varepsilon_n^{2q-qs }} \le 1\quad&\Leftrightarrow\quad \varepsilon_n \ge \frac{\Phi^{2/(2-s)}}{n^{(q-1)/(q(2-s))}}.
\end{align*}
Define
\[
\varepsilon_n := \max\left\{\frac{(\sigma + \Phi)^{2/(2+\alpha)}}{(nA^{-1})^{1/(2+\alpha)}},
\frac{(\sigma + \Phi)}{n^{1/2}},
\frac{\Phi^{2/(2-s)}}{n^{(q-1)/(q(2-s))}},
\left(\frac{A^{q-1}\Phi^{(q(2-s)+ \alpha s(q-1) -\alpha(q-1))}}{n^{q-1}}\right)^{1/(q(2-s)+ \alpha s(q-1) )}\right\}.
\]
From the definition of $\varepsilon_n$, the probability bound simplifies to
\begin{align*}
\mathbb{P}\left(\|\widehat{f} - f_0\| \ge D\varepsilon_n\right) &\le D^{-\tau(1+\alpha/2)} +D^{-\tau(2-s)}+ D^{-(2-s +\alpha s (1-1/q)) \tau} +  \tau^{\tau} D^{-\tau} + \frac{\tau^{\tau(q-1)/(\tau+q-1)}}{D^{[(2q-qs)\tau]/(\tau + q-1)}}.
\end{align*}
Choose $\tau$ such that $\tau \ge q$ and $[(2q-qs)\tau ]/(\tau+q-1) \ge q$ or { equivalently $\tau \ge \max\{q, (q-1)/(1-s)\}$. If $s = 1$ choose $\tau$ such that $\tau \ge q-\eta$ and $[(2-q)\tau+2\tau(q-1)]/(\tau+q-1) \ge q-\eta$ or equivalently $\tau = \max\{q-\eta,(q-1)(q-\eta)/\eta\}$. This choice of $\tau$ implies that}
\[
\mathbb{P}(\|\widehat{f} - f_0\| \ge D\varepsilon_n) \le {C}^qD^{-q + \eta\mathbf{1}\{s=1\}},
\]  
for a constant ${C} > 0$ depending only on $q,s,\alpha,$ and $\eta$. Because $A \ge 1$, $r_n = \varepsilon_n^{-1}$.

\section{Proof of Theorem~\ref{thm:Vc-type}} % (fold)
\label{sec:Vc-type-proof}
Just as in the proofs of the earlier theorems, we will first find an appropriate $\phi_n(\cdot,\cdot)$ and then apply Theorem~\ref{thm:GeneralPeeling} with the right choice of $\tau$ and $s_2(\cdot)$; $\gamma=2$ here. Note that by~definition of~$\Phi$ and $F_\delta$, we  have that 
\begin{equation}\label{eq:udef_1}
U(\epsilon, X, \delta):=2  ( |\epsilon| + \Phi) F_\delta
\end{equation}
satisfies~\eqref{eq:U_def}. Recall that for this theorem, we have assumed that $\epsilon$ has only two moments and $\E(\epsilon^2)\le \sigma^2$. The following lemma proved in Section~\ref{sec:proof_UniformCovering} finds the bound $\phi_n(\cdot, \cdot)$ for various values of $\delta$ and $B$.

\begin{lemma}[Bound on $\phi_n(\cdot,\cdot)$]\label{lem:UniformCoveringPhi}
Let~$\Phi(\cdot, \cdot)$ be as defined in~\eqref{eq:Max_ineq_requirement}. If $\epsilon$ and $\mathcal{F}$ satisfy the assumptions of~Theorem~\ref{thm:Vc-type}, then there exists a constant depending only on $\alpha \ge 0$ and $\beta > 0$\footnote{Here $\alpha$ and $\beta$ are as described in~\eqref{eq:Unif_entr}.} such that 
\begin{enumerate}
\item If $\alpha\in[0, 2)$ and $\beta \ge 0$, then  
\begin{equation}\label{eq:Kolt_Pollard}
\phi_n(\delta, B) \le C A^{1/2}(\sigma + \Phi)\Phi \delta^s, \text{ for every  } \delta , B>0.\footnote{ Note that, we can choose $B$ to be $\infty$.}
\end{equation}

\item If $\alpha=2$ and $\beta=0$, then 
\begin{equation}\label{eq:AlphaEqual2a}
\phi_n(\delta, B) \le C A^{1/2}(\sigma + \Phi)\Phi \delta^s\log\left(A^{-1}n\right), \text{ for all  } n \ge A \text{ and every }  \delta, B >0.
\end{equation}

\item If $\alpha> 2$ and $\beta=0$, then 
\begin{equation}\label{eq:NonDonskera}
\phi_n(\delta, B) \le C A^{1/\alpha}n^{1/2 - 1/\alpha}(\sigma + \Phi)\Phi \delta^s,  \text{ for all  } n \ge A \text{ and every }  \delta, B >0.
\end{equation}
\end{enumerate}
\end{lemma}

\noindent\textbf{Application of Theorem~\ref{thm:GeneralPeeling}:} To apply Theorem~\ref{thm:GeneralPeeling}, we need $\tau,\, \gamma$, and $s_\gamma(\delta)$. If we choose
\begin{equation}\label{eq:gamma_S_2}
\tau= \gamma = 2 \qquad \text{ and }\qquad s_2(\delta) = \E\big(|U(\epsilon, X, \delta)|^2\big)\equiv  2^2 \sigma^2  \Phi^2 \delta^{2s},
\end{equation}
then $s_q(\delta)$ satisfies  $\mathbb{E}\left[U^\gamma(\epsilon, X; \delta)\right] \le s_{\gamma}(\delta)$. Thus by~\eqref{eq:Error_of_fhat}, we have that 
  \begin{align}\label{eq:probbound_beforeB2_VC}
  \begin{split}
  \mathbb{P}\left(\|\widehat{f} - f_0\| \ge D\varepsilon_n\right) &\le 
\left(\frac{C}{\sqrt{n}(D\varepsilon_n)^2}\right)^2\sum_{k = 1}^{\infty}\frac{\phi_n^2(2^kD\varepsilon_n, B_k)}{2^{4k}} + \left(\frac{ C\sqrt{2}(\sigma + \Phi)}{\sqrt{n}D\varepsilon_n}\right)^2\sum_{k = 1}^{\infty}\frac{1}{2^{2k}}\\ &\qquad + \left(\frac{2C}{n(D\varepsilon_n)^2}\right)^2\sum_{k = 1}^{\infty} \frac{B_k^2}{2^{4k  }} + \frac{16}{(D\varepsilon_n)^2}\sum_{k = 0}^{\infty} \frac{s_2(2^kD\varepsilon_n)}{2^{2k}B_k}.
\end{split}
\end{align}
We will now choose $\left\{B_k\right\}$ so as to minimize the right hand side. In contrast to the earlier proofs the $\phi_n(\delta, B)$ obtained in Lemma~\ref{lem:UniformCoveringPhi} do not depend on $B$. Thus the choice of $\left\{B_k\right\}$ that minimizes the bound in~\eqref{eq:probbound_beforeB2_VC} is the value for which the last two terms are equal order for every $k \ge 0$, i.e., 
\begin{equation}\label{eq:B_kbalanace1}
\left(\frac{1}{n(D\varepsilon_n)^2}\right)^2 \frac{ B_k^2}{2^{4k}} = \frac{  \sigma^2 \Phi^2  }{(D\varepsilon_n)^2} \frac{(2^k D \varepsilon_n)^{2s}}{2^{2k}B_k}.
\end{equation}
Thus, we choose 
\begin{equation}\label{eq:B_k1_VC}
  B_k=  \left[ n^2 (D\varepsilon_n)^{2(2+s)}  \sigma^2 \Phi^2   2^{2k(1+2s)}\right]^{1/3 }.
\end{equation}
Then 
\begin{align}\label{eq:B_kbackToBound_VC}
\begin{split}
\left(\frac{2C}{n(D\varepsilon_n)^2}\right)^2\sum_{k = 1}^{\infty} \frac{B_k^2}{2^{4k  }}={}&    \left(\frac{2C}{n(D\varepsilon_n)^2}\right)^2 \left[n^2 (D\varepsilon_n)^{2(2+s)}  \sigma^2 \Phi^2  \right]^{2/3 }\sum_{k = 1}^{\infty} \frac{1}{2^{8k(1-s)/3}}\\
 ={}& \frac{4C^2 \sigma^3 \Phi^3 }{n^{2/3}(D\varepsilon_n)^{4(2-s)/3}}.
\end{split}
\end{align}
Substituting this in~\eqref{eq:probbound_beforeB2_VC}, we get
 \begin{align}
 &\mathbb{P}\left(\|\widehat{f} - f_0\| \ge D\varepsilon_n\right) \\
 \le{}& 
\left(\frac{C}{\sqrt{n}(D\varepsilon_n)^2}\right)^2\sum_{k = 1}^{\infty}\frac{\phi_n^2(2^kD\varepsilon_n, B_k)}{2^{4k}} + \left(\frac{ C\sqrt{2}(\sigma + \Phi)}{\sqrt{n}D\varepsilon_n}\right)^2 + C \frac{4C^2 \sigma^3 \Phi^3 }{n^{2/3}(D\varepsilon_n)^{4(2-s)/3}} \nonumber\\
\le{}& C\left[\left(\frac{1}{\sqrt{n}(D\varepsilon_n)^2}\right)^2\sum_{k = 1}^{\infty}\frac{\phi_n^2(2^kD\varepsilon_n, B_k)}{2^{4k}} + \left(\frac{ (\sigma + \Phi)}{\sqrt{n}D\varepsilon_n}\right)^2 + \frac{ \sigma^3 \Phi^3 }{n^{2/3}(D\varepsilon_n)^{4(2-s)/3}}\right], \label{eq:aftera1a2aplit_andchoosebk_VC}
\end{align}
where $C$ is a universal constant. We will now compute the bound in the~\eqref{eq:aftera1a2aplit_andchoosebk_VC} for all the three cases in~Lemma~\ref{lem:UniformCoveringPhi} by substituting the appropriate values of $\phi_n(\delta, \cdot).$ To aid in this calculation observe that each of the bound on $\Phi_n(\delta, B)$ in~Lemma~\ref{lem:UniformCoveringPhi} is of the form $h(A, \Phi, n,\alpha, \beta) \delta^s$, where $h(A, \Phi, n,\alpha, \beta)$ is a function of it inputs.  Thus~\eqref{eq:aftera1a2aplit_andchoosebk_VC} can be rewritten as
% \todo[inline]{Here $h(A, \Phi, n,\alpha, {\clr \beta})$ is the $\beta$ in~\eqref{eq:Unif_entr}, but it seems like the bounds do not actually depend on $\beta$. I guess, we should still keep it.}
 \begin{align}\label{eq:boundInH}
\begin{split}
 &\mathbb{P}\left(\|\widehat{f} - f_0\| \ge D\varepsilon_n\right)\\
  \le {}& C\left[\left(\frac{1}{\sqrt{n}(D\varepsilon_n)^2}\right)^2\sum_{k = 1}^{\infty}\frac{\phi_n^2(2^kD\varepsilon_n, B_k)}{2^{4k}} + \left(\frac{ (\sigma + \Phi)}{\sqrt{n}D\varepsilon_n}\right)^2 + \frac{ \sigma^3 \Phi^3 }{n^{2/3}(D\varepsilon_n)^{4(2-s)/3}}\right]\\
  \le {}& C\left[\left(\frac{h(A, \Phi, n,\alpha, \beta)}{\sqrt{n}(D\varepsilon_n)^2}\right)^2\sum_{k = 1}^{\infty}\frac{(2^kD\varepsilon_n)^{2s}}{2^{4k}} + \left(\frac{ (\sigma + \Phi)}{\sqrt{n}D\varepsilon_n}\right)^2 + \frac{ \sigma^3 \Phi^3 }{n^{2/3}(D\varepsilon_n)^{4(2-s)/3}}\right]\\
  \le{}&C\left[\frac{h^2(A, \Phi, n,\alpha, \beta)}{n(D\varepsilon_n)^{2(2-s)}} + \left(\frac{ (\sigma + \Phi)}{\sqrt{n}D\varepsilon_n}\right)^2 + \frac{ \sigma^3 \Phi^3 }{n^{2/3}(D\varepsilon_n)^{4(2-s)/3}}\right],
\end{split}
\end{align}
where the last inequality is true as $s \le 1.$  We will find a bound on the tail probability of the LSE and the upper bound on the rate of the LSE by substituting the appropriate $h(A, \Phi, n,\alpha, \beta)$ for each of the cases. We do this below. \\

\noindent\textbf{Case 1:} [$\alpha \in [0,2)$ and $\beta \ge 0$] 
 Because we will allow $A$ and $\Phi$ to depend on $n$ and assume $\sigma$, $\alpha$, and $\beta$ to be constants, we will assume that without loss of generality that $\sigma \le \Phi$.
Substituting~\eqref{eq:Kolt_Pollard} in~\eqref{eq:boundInH}, we have that 
\begin{align}\label{eq:case1vc1}
\begin{split}
 &\mathbb{P}\left(\|\widehat{f} - f_0\| \ge D\varepsilon_n\right)\\
  \le{}&C\left[\frac{\left\{ A^{1/2}\Gamma(1 + \beta/2)(1 - \alpha/2)^{1+\beta/2}(\sigma + \Phi)\right\}^2}{n(D\varepsilon_n)^{2(2-s)}} + \left(\frac{ (\sigma + \Phi)}{\sqrt{n}D\varepsilon_n}\right)^2 + \frac{ \sigma^3 \Phi^3 }{n^{2/3}(D\varepsilon_n)^{4(2-s)/3}}\right]\\
  \le{}&  C \left[\frac{ A\Phi^2}{n(D\varepsilon_n)^{2(2-s)}} + \frac{ \Phi^2}{ n (D\varepsilon_n)^2} + \frac{ \Phi^3 }{n^{2/3}(D\varepsilon_n)^{4(2-s)/3}}\right],
\end{split}
\end{align}
where $C$ is a constant depending only $\alpha, \beta,$ and $ \sigma$.
Now choose, 
\begin{equation}\label{eq:en_vC_case1}
\varepsilon_n:= \max \bigg\{ \left[\frac{A\Phi^2}{n}\right]^{1/2(2-s)}, \frac{\Phi}{\sqrt{n}}, \left[\frac{\Phi^{9/2}}{n}\right]^{1/2(2-s)}\bigg\}.
\end{equation}
With the above choice of $\varepsilon_n$, we have that 
\begin{equation}\label{eq:f-hat_tail_Vc}
\mathbb{P}\left(\|\widehat{f} - f_0\| \ge D\varepsilon_n\right) \le \frac{2 C}{ D^{4(2-s)/3}} + \frac{C}{n^{(1-s)/(2-s)}D^2}.\footnote{If $s < 1/2$, then the tail decays faster than $D^{-2}$ for $D \ll n^{3(1-s)/[(2-4s)(2-s)]}$ and decays like $D^{-2}$ for $D$ larger. Note that $D^{-2}$ is essentially the optimal tail behavior for $\|\widehat{f} - f_0\|$ when the errors only have two moments; see Proposition 1.5 of~\cite{lecue2016performance}.}
\end{equation}
If we fix $A$ and $\Phi$ then it is clear that the rate of convergence of the LSE is no worse than $n^{1/(2(2-s))}$.\\

\noindent\textbf{Case 2:} [$\alpha =2$ and $\beta = 0$] 
 Because we will allow $A$ and $\Phi$ to depend on $n$ and assume $\sigma$ to be a constant, we will assume that without loss of generality that $\sigma \le \Phi$.
Substituting~\eqref{eq:AlphaEqual2a} in~\eqref{eq:boundInH}, we have that 
\begin{align}\label{eq:case1vc2}
\begin{split}
 &\mathbb{P}\left(\|\widehat{f} - f_0\| \ge D\varepsilon_n\right)\\
  \le{}&C\left[\frac{\left\{ A^{1/2}(\sigma + \Phi)\log\left(A^{-1}n\right)\right\}^2}{n(D\varepsilon_n)^{2(2-s)}} + \left(\frac{ (\sigma + \Phi)}{\sqrt{n}D\varepsilon_n}\right)^2 + \frac{ \sigma^3 \Phi^3 }{n^{2/3}(D\varepsilon_n)^{4(2-s)/3}}\right]\\
  \le{}&  C \left[\frac{ A\Phi^2 \log^2\left(A^{-1}n\right) }{n(D\varepsilon_n)^{2(2-s)}} + \frac{ \Phi^2}{ n (D\varepsilon_n)^2} + \frac{ \Phi^3 }{n^{2/3}(D\varepsilon_n)^{4(2-s)/3}}\right],
\end{split}
\end{align}
where $C$ is a constant depending only $\alpha, \beta,$ and $ \sigma$. Then it is easy to see that $\widehat{f}$ satisfies~\eqref{eq:f-hat_tail_Vc} if 
\begin{equation}\label{eq:en_vC_case2}
\varepsilon_n:= \max \bigg\{ \left[\frac{A\Phi^2 \log^2\left(A^{-1}n\right)}{n}\right]^{1/2(2-s)}, \frac{\Phi}{\sqrt{n}}, \left[\frac{\Phi^{9/2}}{n}\right]^{1/2(2-s)}\bigg\}.
\end{equation}
If we fix $A$ and $\Phi$ then it is clear that the rate of convergence of the LSE is no worse than $(\sqrt{n}/\log n)^{1/(2-s)}$.
\\

\noindent\textbf{Case 3:} [$\alpha >2$ and $\beta = 0$] 
 Because we will allow $A$ and $\Phi$ to depend on $n$ and assume $\sigma$ to be a constant, we will assume that without loss of generality that $\sigma \le \Phi$.
Substituting~\eqref{eq:NonDonskera} in~\eqref{eq:boundInH}, we have that 
\begin{align}\label{eq:case1vc}
\begin{split}
 &\mathbb{P}\left(\|\widehat{f} - f_0\| \ge D\varepsilon_n\right)\\
  \le{}&C\left[\frac{\left\{ A^{1/\alpha}(\sigma + \Phi)n^{1/2 - 1/\alpha}\right\}^2}{n(D\varepsilon_n)^{2(2-s)}} + \left(\frac{ (\sigma + \Phi)}{\sqrt{n}D\varepsilon_n}\right)^2 + \frac{ \sigma^3 \Phi^3 }{n^{2/3}(D\varepsilon_n)^{4(2-s)/3}}\right]\\
  \le{}&  C \left[\frac{ A^{2/\alpha}\Phi^2  n^{1-2/\alpha} }{n(D\varepsilon_n)^{2(2-s)}} + \frac{ \Phi^2}{ n (D\varepsilon_n)^2} + \frac{ \Phi^3 }{n^{2/3}(D\varepsilon_n)^{4(2-s)/3}}\right],
\end{split}
\end{align}
where $C$ is a constant depending only $\alpha, \beta,$ and $ \sigma$. Then it is easy to see that $\widehat{f}$ satisfies~\eqref{eq:f-hat_tail_Vc} if 
\begin{equation}\label{eq:en_vC_case_3}
\varepsilon_n:= \max \bigg\{ \left[\frac{ A^{2/\alpha}\Phi^2}{n^{2/\alpha}}\right]^{1/2(2-s)}, \frac{\Phi}{\sqrt{n}}, \left[\frac{\Phi^{9/2}}{n}\right]^{1/2(2-s)}\bigg\}.
\end{equation}
As in the previous cases, if we fix $A$ and $\Phi$, the rate of convergence of the LSE is no worse than $n^{1/\alpha(2-s)}$.

Finally, for each of the above cases $s\le 1$, thus we have that $\varepsilon_n^{-1} \E\|\widehat{f}-f_0\| = O(1)$.

\subsection{Proof of Lemma~\ref{lem:UniformCoveringPhi}}\label{sec:proof_UniformCovering}

The proof is based on the use of symmetrization by Rademacher variables followed by application of the sub-Gaussian maximal inequality given by Corollary 2.2.8 of~\cite{VdVW96} conditionally on $\{(\epsilon_i, X_i), 1\le i\le n\}$.  % \begin{proof}
Recall that
\[
\phi_n(\delta; B) = \mathbb{E}\left[\sup_{\delta/2\le \|f - f_0\| \le \delta}\mathbb{G}_n(T_B(f; \epsilon, X, \delta))\right].
\]
By Symmetrization (Corollary 3.2.2 of~\cite{Gine16}), we get
\[
\phi_n(\delta; B) \le 2\mathbb{E}\left[\sup_{f\in\mathcal{F}_{\delta}}\,\frac{1}{\sqrt{n}}\sum_{i=1}^n R_iT_B(f; \epsilon_i, X_i, \delta)\right],
\]
where $R_1,\ldots,R_n$ are i.i.d.~Rademacher random variables independent of $(\epsilon_1,X_1),\ldots,(\epsilon_n,X_n).$ Now by Lemma A.1 of~\cite{srebro2010optimistic} (also see~\cite[Theorem 3.2]{chatterjee2018matrix}), we have
\begin{equation}\label{eq:Refined-dudley-inequality}
\phi_n(\delta; B) \le \mathbb{E}\left[\inf_{\gamma\ge0}4n^{1/2}\gamma + 10\int_{\gamma}^{\eta_n} \sqrt{\log(\eta, \{T_B(f):f\in\mathcal{F}_{\delta}\}, \|\cdot\|_n)}d\eta\right],
\end{equation}
where
$\|g\|_n^2 := n^{-1}\sum_{i=1}^n g^2(\epsilon_i, X_i),$ for $g:\mathbb{R}\times\rchi\to\mathbb{R}$, and $\eta_n := \sup_{f\in\mathcal{F}_{\delta}}\|T_B(f)\|_n$. It is clear that
\begin{equation}\label{eq:TBIneq}
|T_B(f_1; \epsilon, X) - T_B(f_2; \epsilon, X)| \le (2|\epsilon| + 4M)|f_1(X) - f_2(X)|\mathbf{1}\{U(\epsilon, X; \delta) \le B\}.
\end{equation}
Define a measure $Q$ on $\{X_1, \ldots, X_n\}$ as
\[
Q(\{X_i\}) = \frac{(2|\epsilon_i| + 4M)^2\mathbf{1}\{U(\epsilon_i, X_i; \delta) \le B\}}{\sum_{j=1}^n (2|\epsilon_j| + 4M)^2\mathbf{1}\{U(\epsilon_j,X_j;\delta) \le B\}},\quad 1 \le i\le n.
\]
From inequality~\eqref{eq:TBIneq}, we get
\[
\|T_B(f_1; \epsilon, X) - T_B(f_2; \epsilon, X)\|_{n} \le \|(2|\epsilon| + 4M)\mathbf{1}\{U(\epsilon, X;\delta) \le B\}\|_{n}\|f_1 - f_2\|_{2,Q}.
\]
Thus it follows that for any $\eta > 0,$
\begin{align}\label{eq:6878}
\begin{split}
&\log N(\|(2|\epsilon| + 4M)\mathbf{1}\{U(\epsilon, X;\delta) \le B\}\|_{n}\eta, \{T_B(f)\}, \|\cdot\|_n)\\
&\qquad\le \log N\left({\eta}, \mathcal{F}_{\delta}, L_2(Q)\right) \le A\left(\frac{\eta}{\|F_{\delta}\|_{2,Q}}\right)^{-\alpha}\log^{\beta}\left(\frac{\|F_{\delta}\|_{2,Q}}{\eta}\right).
\end{split}
\end{align}
Hence using the fact $\|U(\epsilon, X; \delta)\mathbf{1}\{U(\epsilon, X; \delta) \le B\}\|_{n} = \|(2|\epsilon|+4M)\mathbf{1}\{U(\epsilon,X,\delta) \le B\}\|_n\|F_{\delta}\|_{2,Q}$,~\eqref{eq:6878} yields
\[
\log N(\|U(\epsilon,X;\delta)\mathbf{1}\{U(\epsilon,X;\delta) \le B\}\|_n\eta, \mathcal{F}_{\delta}, \|\cdot\|_n) \le A\eta^{-\alpha}\log^{\beta}(1/\eta)\quad\mbox{for all}\quad \eta > 0.
\]
Substituting this bound in~\eqref{eq:Refined-dudley-inequality}, we get
\[
\phi_n(\delta;B) \le \mathbb{E}\left[\|U(\epsilon,X;\delta)\mathbf{1}\{U(\epsilon,X;\delta) \le B\}\|_n\inf_{\gamma\ge0}4\sqrt{n}\gamma + 10A^{1/2}\int_{\gamma}^{\Theta_n} \eta^{-\alpha/2}\log^{\beta/2}(1/\eta)d\eta\right],
\]
where $\Theta_n := \sup_{f\in\mathcal{F}_{\delta}}\|T_B(f; \epsilon, \delta)\|_n/\|U(\epsilon, X; \delta)\mathbf{1}\{U(\epsilon, X;\delta) \le B\}\|_n$. Because $\Theta_n \le 1$ and the infimum is a non-random quantity, we obtain
\begin{align*}
\phi_n(\delta;B) &\le \mathbb{E}\left[\|U(\epsilon,X;\delta)\mathbf{1}\{U(\epsilon,X;\delta) \le B\}\|_n\right]\left[\inf_{\gamma\ge0}4\sqrt{n}\gamma + 10A^{1/2}\int_{\gamma}^{1} \eta^{-\alpha/2}\log^{\beta/2}(1/\eta)d\eta\right]\\
% &\le \inf_{\gamma \ge 0}4\sqrt{n}\mathbb{E}[\|U(\epsilon,X;\delta)\mathbf{1}\{U(\epsilon,X;\delta)\le B\}\|_n]\gamma + 10A^{1/2}\int_{\gamma}^1 \eta^{-\alpha}\log^{\beta/2}(1/\eta)d\eta\\
&\le \|U(\epsilon,X;\delta)\|\left[\inf_{\gamma \ge 0}4\sqrt{n}\gamma + 10A^{1/2}\int_{\gamma}^1 \eta^{-\alpha}\log^{\beta/2}(1/\eta)d\eta\right]\\
&\le C(\sigma+\Phi)\Phi\delta^sG_n,
\end{align*}
where $C$ is a universal constant and
\[
G_n := \inf_{\gamma\ge0}\sqrt{n}\gamma + A^{1/2}\int_{\gamma}^1 \eta^{-\alpha/2}\log^{\beta/2}(1/\eta)d\eta.
\]
We now complete the proof by bounding $G_n$ separately in each of the three cases.
\begin{enumerate}
\item[Proof of~\eqref{eq:Kolt_Pollard}:] If $\alpha < 2$ and $\beta > 0$, then we can take $\gamma = 0$ in the infimum of $G(\cdot)$ and using $\int_0^1 \eta^{-\alpha/2}\log^{\beta/2}(1/\eta)d\eta < \infty$, we get
\[
\phi_n(\delta; B) \le CA^{1/2}(\sigma + \Phi)\Phi\delta^s,
\]
for a constant $C > 0$ depending only on $\alpha, \beta$.
\item[{Proof of~\eqref{eq:AlphaEqual2a}:}]
If $\alpha = 2$ and $\beta = 0$, taking $\gamma = (A^{-1}n)^{-1/2}$ yields
\[
G_n \le A^{1/2} + A^{1/2}\int_{(A^{-1}n)^{-1/2}}^1 \eta^{-1}d\eta = 1 + A^{1/2}\log(A^{-1/2}n^{1/2}) \le A^{1/2}\log(A^{-1}n),
\]
and hence $\phi_n(\delta; B) \le CA^{1/2}(\sigma+\Phi)\Phi\delta^s\log(A^{-1}n).$
\item[{Proof of~\eqref{eq:NonDonskera}:}]

 If $\alpha > 2$ and $\beta = 0$, taking $\gamma = (A^{-1}n)^{-1/\alpha}$ yields
\begin{align*}
G_n &\le A^{1/\alpha}n^{1/2 - 1/\alpha} + A^{1/2}\int_{A^{1/\alpha}n^{-1/\alpha}}^1 \eta^{-\alpha/2}d\eta\\
&= A^{1/\alpha}n^{1/2 - 1/\alpha} + \frac{A^{1/2}}{\alpha/2 - 1}[(A^{-1}n)^{1/2 - 1/\alpha} - 1] \le CA^{1/\alpha}n^{1/2 - 1/\alpha},
\end{align*}
for some constant $C > 0$ depending only on $\alpha$. Hence $$\phi_n(\delta;B) \le C A^{1/\alpha}n^{1/2 - 1/\alpha}(\sigma+\Phi)\Phi\delta^s.$$
\end{enumerate}

\section{Proof of Theorem~\ref{thm:Lowerbnd}} % (fold)
\label{sec:proof_of_theorem_thm:lowerbnd}

First consider the case when $s\in \{0,1\}$.
 When $s=1$, let $\widetilde\F:= \{\mathbf{1}_{[a,1]}: a\in [0,1]\}$, in this case~\cite[Example 3]{han2018robustness} show that $\widetilde F_\delta =\mathbf{1}_{[1-\delta^2, 1]}$ (implying that $s=1$). They further show that LSE converges at the parametric rate of $n^{-1/2}$, clearly satisfying~\eqref{eq:lwrbnd}. 

When $s = 0$, fix any $k>1$ and let  
\[\widetilde\F :=\{\sum_{i=1}^k c_i \mathbf{1}_{[x_{i-1},x_i]}: |c_i|\le 1, 0\le x_0< x_1<\ldots<x_k\le 1\}.
\]
Example~4 of~\cite{han2018robustness}  shows that $\widetilde F_\delta\equiv 1$ and hence $s=0$. Moreover, they show that the LSE cannot converge at a rate faster than $n^{-1/4}$ when $\epsilon$ is independent and has roughly two moments\footnote{If $\epsilon$ is Gaussian or bounded, then one can show that the LSE converges at a rate no worse than $\sqrt{\log n/n}$;~\cite[Example 4]{han2018robustness}}.

Now, fix $s \in (0,1)$. The proof below is constructive, i.e., we will give a specific choice of $\widetilde{\F}$ (depending on $s$) and $\xi$ such that the LSE satisfies~\eqref{eq:lwrbnd}. The proof is almost identical to the proof of Proposition~2 in~\cite{han2018robustness}; we use the $\widetilde{\F}$ constructed in Section 5.1.1 of~\cite{han2018robustness}.  However, there are two differences: (1) The chosen error distribution is different; (2) Instead of using the Paley-Zygmund for the lower bound calculations, we use Nagaev inequalities~\cite{nagaev1979large}. These differences lead to different choices of $\delta_1$ (~\cite[Page 33, Section 5.1.3]{han2018robustness}) and $\delta_2$ (~\cite[Page 29, Section 5.1.2]{han2018robustness}) and a different lower bound. In the following, we only detail the parts of the proof that are different. For a full proof, we refer the reader to the proof of Theorem~2 of \cite{han2018robustness}.  To ease the reading, we have kept the notation (almost) identical to those in~\cite{han2018robustness} except they use $\gamma$ instead of $s$ and $\xi$ instead of $\epsilon$. 

 Consider the function class $\widetilde{\mathcal{F}}$ defined in Section 5.1.1 of~\cite{han2018robustness}, with $\gamma = s$. This function class is contained in the set of all intervals of $[0, 1]$ and hence is a VC class (i.e, $\alpha = 0$ and $\beta = 1$ in condition (VC)). For this function class with $f_0 \equiv 0$, the local envelope $\widetilde{F}_{\delta}$ satisfies $\|\widetilde{F}_{\delta}\|_2 \le \sqrt{2}\delta^{s}$; see Lemma 14 of~\cite{han2018robustness} for a proof. 

 Fix $0< \delta_2 = o(1)$. Then $N(\delta_2)$, $\mathcal{I}_l$, and $\mathcal{E}_n$ defined in~\cite{han2018robustness} satisfy: (1) $ N(\delta_2) \asymp \delta_2^{-2(2-2s)};$ (2)  on the set $\mathcal{E}_n$ we have $|\mathcal{I}_l| \asymp n \delta_2^2;$ and (3) $\mathbb{P}(\mathcal{E}_n^c) \lesssim 1/n.$

%  \[  \quad  |\mathcal{I}_l| \asymp n \delta_2^2, \quad \text{ and } \quad
% \mathbb{P}(\mathcal{E}_n^c) \lesssim C\delta_2^{-(2-2s)}\exp\left(-C'\delta_2^{-(2-2s)}\right).
%  \]
% $\mathbb{P}\left(\max_{1\le l \le N(\delta_2)} ||\mathcal{I}_l| - n \delta_2^2| \le \sqrt{n \delta_2^2\log n}+ \log n \right) \ge 1- 1/n$ 
% $sqrt{n \delta_2^2\log n}+ \log n$

Define a symmetric random variable $\xi$ with distribution
\[
\mathbb{P}(|\xi| \ge t) = \frac{\log^2(2)}{t^2\log^2(1 + t)},
\]
for $t \ge 1$. Because $\mathbb{P}(|\xi| \ge 0) = 1$ and $\lim_{t\to\infty}\mathbb{P}(|\xi| \ge t) = 0$.
It is clear that
\[
\mathbb{E}[\xi^2] = 2\int_0^{\infty} t\mathbb{P}(|\xi| \ge t)dt \le \int_0^{\infty} \frac{2\log^22}{t\log^2(1 + t)} dt\lesssim \sum_{n=1}^{\infty} \frac{1}{n\log^2n} < \infty.
\]
Further
\[
\mathbb{E}[|\xi|^{2+\delta}] = (2+\delta)\int_0^{\infty} \frac{t^{1+\delta}}{t^2\log^2(1+t)}dt \gtrsim \int_1^{\infty} \frac{1}{t^{1-\delta}\log^2(1 + t)}dt \gtrsim \int_1^{\infty} \frac{1}{t^{1-\delta/2}} dt= \infty.
\]
Therefore, $\xi$ only has two moments and no more. Let $\xi, \xi_1, \xi_2, \ldots, \xi_n$ are independent and identically distributed random variables. Following inequality (5.2) of~\cite{han2018robustness}, we have 
\begin{equation}\label{eq:52'}
\mathbb{P}\left(\sup_{f\in\widetilde{\mathcal{F}}_{\delta_2}}\sum_{i=1}^n \xi_if(X_i) \ge t_n\right) \ge \mathbb{E}_{X}\left[\mathbb{P}_{\xi}\left(\max_{1\le l\le N(\delta_2)}\sum_{i=1}^{|\mathcal{I}_l|} \xi_i^{(l)} \ge t_n\right)\mathbf{1}\{\mathcal{E}_n\}\right],
\end{equation}
with $N(\delta_2) \in [\delta_2^{-(2-2s)}, 2\delta_2^{-(2-2s)}]$ and $\xi_i^{(l)}$  as defined in~\cite{han2018robustness}.
%  and an event $\mathcal{E}_n$ satisfying
% \[
% \mathbb{P}(\mathcal{E}_n^c) \le C\delta_2^{-(2-2s)}\exp\left(-C'\delta_2^{-(2-2s)}\right).
% \]
% Further by the discussion in the beginning of Proof of Proposition 2, claim (1) in~\cite{han2018robustness}, we have that 
% \[
% |\mathcal{I}_l| ~\asymp~ n \delta_2^2\quad\mbox{and}\quad N(\delta_2) \asymp \delta_2^{-(2-2s)}.
% \]
Set $m_l := |\mathcal{I}_l|$ and $S_l := \sum_{i=1}^{|\mathcal{I}_l|} \xi_i^{(l)}$. The inequality following (1.58)  of~\cite[Page 759]{nagaev1979large} implies that
\begin{equation}\label{eq:Nagaev-implication}
\mathbb{P}\left(\sum_{i=1}^{m_l} \xi_i^{(l)} \ge t\right) \ge \frac{m_l}{2}\mathbb{P}\left(\xi_i^{(l)} \ge 2 t\right) = \frac{m_l\log^2(2)}{8t^2\log^2(1 + 2t)} \quad\text{when}\quad t \ge 2\sqrt{m_l\mbox{Var}(\xi)}.
\end{equation}
We will now use the arguments in \cite[Page 22]{MR1666908} to bound the probability on the right of~\eqref{eq:52'}:
\begin{align}\label{eq:delapenagine}
\begin{split}
\mathbb{P}\left(\max_{1\le l\le N(\delta_2)}S_l \ge t_n\right) &\ge 1- \prod_{i=1}^{N(\delta_2)}\left(1- \mathbb{P}\left(S_l \ge t_n\right)\right)\\
&\ge 1- \exp\left(-\sum_{i=1}^{N(\delta_2)}\mathbb{P}\left(S_l \ge t_n\right)\right)\\
&\ge \frac{\sum_{i=1}^{N(\delta_2)}\mathbb{P}\left(S_l \ge t_n\right)}{1+ \sum_{i=1}^{N(\delta_2)}\mathbb{P}\left(S_l \ge t_n\right)}.
\end{split}
\end{align}
Recalling that $ x/(1+x) \ge 1/2$ if and only if $x\ge1$ and combining~\eqref{eq:Nagaev-implication} and~\eqref{eq:delapenagine}, we have 
\begin{equation}\label{eq:comb_naga_gine}
\mathbb{P}\left(\max_{1\le l\le N(\delta_2)}S_l \ge t_n\right) \ge 1/2  \qquad\text{when}\qquad \frac{\log^2(2)}{8t^2_n\log^2(1 + 2t_n)} \sum_{i=1}^{N(\delta_2)} m_l \ge 1
\end{equation}
Thus, we have that
\[ \mathbb{P}\left(\max_{1\le l\le N(\delta_2)}S_l \ge t_n\right) \ge 1/2 \qquad\text{when}\qquad t_n = c \frac{\sqrt{\sum_{i=1}^{N(\delta_2)} m_l }}{\log (\sum_{l=1}^{N(\delta_2)} m_l)},\]
for some constant $c > 0$.
Note that with $N(\delta_2)$ is diverging to infinity and $C_1 \le m_l/m_{l'} \le C_2$ for all $l, l'$, thus $$\frac{\sqrt{\sum_{l=1}^{N(\delta_2)} m_l}}{\log(\sum_{l=1}^{N(\delta_2)} m_l)} \gg 2\max_{1\le l\le N(\delta_2)}\sqrt{m_l\mbox{Var}(\xi)},$$ and hence~\eqref{eq:Nagaev-implication} is applicable.
Therefore,
\begin{equation}\label{eq:tightlwrbound}
\mathbb{P}\left(\sup_{f\in\widetilde{\mathcal{F}}:\|f\|_2 \le \delta_2}\sum_{i=1}^n \xi_if(X_i) \ge \frac{C' n^{1/2}\delta_2^{s}}{2(2-2s) \log(1/\delta_2)}\right) \ge 1/2.
\end{equation}
Further (5.7) of~\cite{han2018robustness} yields that with probability at least $1 - 1/n$,
\[
\sup_{f\in\widetilde{\mathcal{F}}:\|f\|\le\delta_2}|\mathbb{G}_nf^2| \le C\delta_2\sqrt{\log n},
\]
and hence with probability at least $1/2 - 1/n$, we have
\begin{align*}
\sup_{f\in\widetilde{\mathcal{F}}:\|f\| \le \delta_2}(\mathbb{P}_n - P)(2\xi f - f^2) &\gtrsim \frac{\delta_2^{s} n^{-1/2}}{\log n} - \delta_2\sqrt{\log n}n^{-1/2}\\
&\gtrsim \frac{\delta_2^{s} n^{-1/2}}{\log n}, \text{ when }  \delta_2= o\big((\log n)^{-3/(2(1-s))}\big).
\end{align*} This implies that with probability at least $1/2 - 1/n$,
\begin{align*}
F_n(\delta_2) &\ge \frac{C\delta_2^{s} n^{-1/2}}{\log n} - \delta_2^2 \gtrsim \frac{C\delta_2^{s} n^{-1/2}}{\log n}, \text{ when } \delta_2 \lesssim (n\log^2 n)^{-1/(2(2-s))}.
\end{align*}
Taking $\delta_2 \asymp (n\log^2 n)^{-1/(2(2-s))}$, we have that $F_n(\delta_2) \gtrsim (n\log^2 n)^{-1/(2-s)}$  with probability at least $1/2 - 1/n$.

Now to get an upper bound on $E_n(\delta_1)$ (defined in~\cite[Page 18]{han2018robustness}, observe that by Lemma~\ref{lem:UniformCoveringPhi}, we get that 
\[
\mathbb{E}\left[\sup_{f\in\widetilde{\mathcal{F}}:\|f\|_2 \le \delta_1}\sum_{i=1}^n \xi_if(X_i)\right] \lesssim \delta^{s}n^{1/2} (\|\xi\|_{2}+1),
\]
and following the proof of Proposition 2 (claim (2)) in~\cite{han2018robustness}, we have that
\[
\mathbb{E}\left[\sup_{f\in\widetilde{\mathcal{F}}:\|f\|_2 \le \delta_1}|\mathbb{G}_nf^2|\right] \lesssim \delta_1^{s}.
\]
Therefore,
\[
\mathbb{E}[E_n(\delta_1)] \lesssim \delta_1^{s}n^{1/2}n^{-1} + n^{-1/2}\delta_1^{s} \asymp \delta_1^{s}n^{-1/2}.
\]
Hence with probability at least $3/4$, $E_n(\delta_1) \le C\delta_1^{s}n^{-1/2}$. Taking $\delta_1 \lesssim n^{-1/(2(2-s))}/(\log n)^{2/(s(2-s))}$, we get $\delta_1 \le \delta_2$ and $E_n(\delta_1) < F_n(\delta_2)$ with positive probability. Hence, by Proposition 1 of~\cite{han2018robustness}, there exists a least squares estimator $\widehat{f}$ such that
\[
\|\widehat{f} - f_0\| = \|\widehat{f}\| \gtrsim n^{-1/(2(2-s))}/(\log n)^{2/(s(2-s))}.
\]
\noeqref{eq:lower-bound-main}
\noeqref{eq:TrueTailBound_L2}
\noeqref{eq:ep_n_extra_log}
\noeqref{eq:en_vC_case1}
\noeqref{eq:en_vC_case2}
\noeqref{eq:en_vC_case_3}
% \bibliographystyle{chicago}
% \bibliography{SigNoise}

\end{document}